\documentclass[11pt, oneside]{book}
\usepackage{tocloft}

\usepackage{amsmath,amsthm,amssymb}

\newcommand{\remove}[1]{} 
\newtheorem{thm}{Theorem}[section]
\newtheorem{claim}[thm]{Claim}
\newtheorem{lem}[thm]{Lemma}
\newtheorem{define}[thm]{Definition}
\newtheorem{cor}[thm]{Corollary}

\def\F{{\mathbb{F}}}
\def\Q{{\mathbb{Q}}}
\def\Z{{\mathbb{Z}}}
\def\N{{\mathbb{N}}}
\def\R{{\mathbb{R}}}
\def\K{{\mathbb{K}}}
\def\C{{\mathbb{C}}}
\def\A{{\mathbb{A}}}
\def\P{{\mathbb{P}}}
\def\cP{{\cal P}}
\def\cR{{\cal R}}
\def\L{{\cal L}}
\def\B{{\mathcal B}}
\def\S{{\mathbf S}}
\def\cE{{\mathcal E}}
\def\V{{\mathbf{V}}}
\def\I{{\mathbf{I}}}
\def\bx{{\mathbf x}}
\def\by{{\mathbf y}}
\def\E{{\mathbb E}}
\def\SD{{\mathbf{SD}}}
\def\SC{{\mathbf{SC}}}

\def\half{ \frac{1}{2}}
\newcommand{\ip}[2]{\langle #1,#2 \rangle}
\def\sumN{\sum_{i=1}^n}
\def\_{\,\,\,\,\,}
\def\prob{{\mathbf{Pr}}}
\newcommand{\entropy}[1]{ {\text{H}_{\infty}\left({#1}\right)} }
\def\modulo{\text{mod}}
\def\omm{ \{0,1\} }
\def\id{ \textit{id} }

\def\D{{\partial}}
\def\gap{\textsf{gap}}
\def\sign{\textsf{sign}}
\def\bias{\textsf{bias}}
\def\spar{\textsf{sparse}}
\def\span{\textsf{span}}
\def\Part{\textbf{Part}}
\def\Mon{\textbf{Mon}}
\def\sing{\textbf{sing}}
\def\Und{\textsf{Und}}
\def\Comp{\textsf{Comp}}
\def\rank{\textsf{rank}}
\def\poly{\textsf{poly}}
\def\codim{\textsf{codim}}
\def\cp{\textsf{cp}}
\def\uni{\textsf{Uni}}
\def\ext{\textsf{\bf Ext}}
\def\extt{\textsf{\bf Ext2}}
\def\cros{{\textbf{cr}}}
\def\dist{{\textbf{dist}}}
\def\SL{{\textbf{SL}}}
\def\stab{{\textrm{Stab}}}
\def\Dec{{\textrm{Dec}}}

\def\fplus{{\,+_f\,}}
 
\newcommand{\epsclose}{\stackrel{\epsilon}{\thicksim} }
\newcommand{\eclose}[1]{\stackrel{{#1}}{\thicksim} }
\newcommand{\eps}{\epsilon}

\begin{document}

\title{\textbf{Incidence Theorems and Their Applications}}
\author{Zeev Dvir\thanks{Princeton university, Mathematics and Computer Science Departments. Princeton, NJ, 08540. zdvir@princeton.edu.}}
\date{}   
\maketitle

\chapter*{\centering \begin{normalsize}Abstract\end{normalsize}}
\begin{quotation}
\noindent We survey recent (and not so recent) results concerning arrangements of lines, points and other geometric objects and the applications these results have in theoretical computer science and combinatorics. The  three main types of problems we will discuss are:
\begin{enumerate}
	\item {\bf Counting incidences:} Given a set (or several sets) of geometric objects (lines, points, etc.), what is the maximum number of incidences (or intersections) that can exist between elements in different sets? We will see several results of this type, such as the Szemeredi-Trotter theorem, over the reals and over finite fields and discuss their applications in combinatorics (e.g., in the recent solution of Guth and Katz to Erdos' distance problem) and in computer science (in explicit constructions of multi-source extractors).
	\item {\bf Kakeya type problems:} These problems deal with arrangements of lines that point in different directions. The goal is to try and understand to what extent these lines can overlap one another. We will discuss these questions both over the reals and over finite fields and see how they come up in the theory of randomness-extractors. 
	\item {\bf Sylvester-Gallai type problems:} In this type of problems, one is presented with a configuration of points that contain many `local' dependencies (e.g., three points on a line) and is asked to derive a bound on the dimension of the span of all points. We will discuss several recent results of this type, over various fields, and see their connection to the theory of locally correctable error-correcting codes.
\end{enumerate}

Throughout the different parts of the survey, two types of techniques will make frequent appearance. One is the polynomial method, which uses polynomial interpolation to impose an algebraic structure on the problem at hand. The other recurrent techniques will come from the area of additive combinatorics. 

\end{quotation}
\clearpage

\tableofcontents
\newpage
  
\chapter{Overview}\label{chap-overview}

\newif\ifediting

\ifediting
\documentclass[11pt]{article}

\usepackage{amsmath,amsthm,amssymb}
\newcommand{\remove}[1]{}
\setlength{\topmargin}{0.3in} \setlength{\headheight}{0in}
\setlength{\headsep}{0in} \setlength{\textheight}{8.0in}
\setlength{\topsep}{0.1in} \setlength{\itemsep}{0.0in}
\parskip=0.05in
 \textwidth=6.5in 
\oddsidemargin=0truecm \evensidemargin=0truecm

\newtheorem{thm}{Theorem}[section]
\newtheorem{claim}[thm]{Claim}
\newtheorem{lem}[thm]{Lemma}
\newtheorem{define}[thm]{Definition}
\newtheorem{cor}[thm]{Corollary}
\newtheorem{obs}[thm]{Observation}
\newtheorem{example}[thm]{Example}
\newtheorem{construct}[thm]{Construction}
\newtheorem{conjecture}[thm]{Conjecture}
\newtheorem{THM}{Theorem}
\newtheorem{question}{Question}
\newtheorem{fact}[thm]{Fact}
\newtheorem{prop}[thm]{Proposition}

\def\F{{\mathbb{F}}}
\def\Q{{\mathbb{Q}}}
\def\Z{{\mathbb{Z}}}
\def\N{{\mathbb{N}}}
\def\R{{\mathbb{R}}}
\def\K{{\mathbb{K}}}
\def\C{{\mathbb{C}}}
\def\A{{\mathbb{A}}}
\def\P{{\mathbb{P}}}
\def\cP{{\cal P}}
\def\cS{{\mathcal S}}
\def\cE{{\mathcal E}}
\def\V{{\mathbf{V}}}
\def\I{{\mathbf{I}}}
\def\bx{{\mathbf x}}
\def\by{{\mathbf y}}
\def\E{{\mathbb E}}

\def\half{ \frac{1}{2}}
\newcommand{\ip}[2]{\langle #1,#2 \rangle}
\def\sumN{\sum_{i=1}^n}
\def\_{\,\,\,\,\,}
\def\prob{{\mathbf{Pr}}}
\newcommand{\entropy}[1]{ {\text{H}_{\infty}\left({#1}\right)} }
\def\modulo{\text{mod}}
\def\omm{ \{0,1\} }
\def\id{ \textit{id} }

\def\D{{\partial}}
\def\gap{\textsf{gap}}
\def\sign{\textsf{sign}}
\def\spar{\textsf{sparse}}
\def\span{\textsf{span}}
\def\Part{\textbf{Part}}
\def\Mon{\textbf{Mon}}
\def\sing{\textbf{sing}}
\def\Und{\textsf{Und}}
\def\Comp{\textsf{Comp}}
\def\rank{\textsf{rank}}
\def\poly{\textsf{poly}}
\def\codim{\textsf{codim}}
\def\cp{\textsf{cp}}
\def\uni{\textsf{Uni}}
\def\ext{\textsf{\bf Ext}}
\def\extt{\textsf{\bf Ext2}}

\def\fplus{{\,+_f\,}}

\newcommand{\epsclose}{\stackrel{\epsilon}{\thicksim} }
\newcommand{\eclose}[1]{\stackrel{{#1}}{\thicksim} }
\newcommand{\eps}{\epsilon}
\newcommand{\Anote}[1]{\begin{quote}{\sf Avi's Note:} {\sl{#1}} \end{quote}}

\begin{document}

\title{Incidence Theorems -- Lecture Notes}
\date{}
\maketitle


\fi

Consider a finite set of points, $P$, in some vector space and another set $L$ of lines. An {\em incidence} is a pair $(p,\ell) \in P \times L$ such that $p \in \ell$. There are many types of questions one can ask about the set of incidences and many different conditions one can impose on the corresponding set of points and lines. For example, the Szemeredi-Trotter theorem (which will be discussed at length below) gives an upper bound on the {\em number} of possible incidences. More generally, in this survey we will be interested in a variety of problems and theorems relating to arrangements of lines and points and the surprising applications these theorems have, in theoretical computer science and in combinatorics. The term `incidence theorems' is used in a very broad sense and might include results that could fall under other categories. We will study questions about incidences between lines and points, lines and lines (where an incidence is a pair of intersecting lines), circles and points and more. 

Some of the results we will cover have direct and powerful applications to problems in theoretical computer sciences and combinatorics. One example in combinatorics is the recent solution of Erdos' distance problem by Guth and Katz \cite{GK10}. The problem is to lower bound the number of distinct distances defined by a set of points in the real plane and the solution (which is optimal up to logarithmic factors) uses a clever reduction to a problem on counting incidences of lines \cite{ES10}. 

In theoretical computer science, incidence theorems (mainly over finite fields) have been used in recent years to construct {\em extractors}, which are procedures that  transform {\em weak} sources of randomness (that is, distributions that have some amount of randomness but are not completely uniform) into completely uniform random bits. Extractors have many theoretical applications, ranging from  cryptography to data structures to metric embeddings (to name just a few) and the current state-of-the-art constructions all use incidence theorems in one way or another. The need to understand incidences comes from trying to analyze simple looking  constructions that use basic algebraic operations. For example, how `random' is $X\cdot Y + Z$, when $X,Y,Z$ are three independent random variables each distributed uniformly over a large subset of $\F_p$.

We will see incidence problems over finite fields, over the reals, in low dimension and in high dimension. These changes in field/dimension are pretty drastic and, as a consequence, the ideas appearing in the proofs will be quite diverse. However, two main techniques will make frequent appearance. One is the `polynomial method' which uses polynomial interpolation to try and  `force' an algebraic structure on the problem. The other  recurrent techniques will come from additive combinatorics. These are general tools to argue about sets in Abelian groups and the way they behave under certain group operations. These two techniques are surprisingly flexible and can be applied in many different scenarios and over different fields.

The survey is divided into four chapters, following this overview chapter. The first chapter will be devoted to problems of counting incidences over the real numbers (Szemeredi-Trotter and others) and will contain applications mostly from combinatorics (including the Guth-Katz solution to  Erdos' distance problem). The second chapter will be devoted to the Szemeredi-Trotter theorem over finite fields and its applications to the explicit constructions of multi-source extractors. The third chapter will be devoted to {\em Kakeya} type problems which deal with arrangements of lines pointing in different directions (over finite and infinite fields). The applications in this chapter will be to the construction of another variant of extractors -- seeded extractors. The fourth and final chapter will deal with arrangements of points with many collinear triples. These are related to questions  in theoretical computer science having to do with {\em locally correctable} error correcting codes. More details and definitions relating to each of the aforementioned chapters are given in the next four subsections of this overview which serves as a road map to the various sections.

This survey is aimed at both mathematicians and computer scientists and could serve as a basis for a one semester course. Ideally, each chapter should be read from start to finish (the different chapters are mostly independent of each other). We  only assume familiarity with undergraduate level algebra, including the basics of finite fields and polynomials. 


\paragraph{Notations:} We will use $\lesssim, \gtrsim$ and $\sim$ to denote (in)equality up to multiplicative absolute constants. That is, $X \lesssim Y$ means `there exists an absolute constant $C$ such that $X \leq CY$'. In some places, we opt to use instead the computer science notations of $O(\cdot), \Omega(\cdot)$ and $\theta(\cdot)$ to make some expressions more readable. So $X = O(Y)$ is the same as $X \lesssim Y$, $X = \Omega(Y)$ is the same as $X \gtrsim Y$ and $X  = \theta(Y)$ is the same as $X \sim Y$. This allows us to write, for example, $X = 2^{\Omega(Y)}$ to mean that there exists an absolute constant $C$ such that $x \geq 2^{CY}$.

\paragraph{Sources:} Aside from research papers there were two main sources that were used in the preparation of this survey. The first is a sequence of posts on Terry Tao's blog which cover a large portion of Chapter~\ref{chap-1} (see e.g. \cite{TaoCell}). Ben Green's lecture notes on additive combinatorics \cite{GreenCourse}  were the main source in preparing the Chapter~\ref{chap-2}. Both of these sources were indispensable in preparing this survey and I am grateful to both authors.

\subsection*{Chapter~\ref{chap-1}: Counting incidences over the reals}

Let $P$ be a finite set of points  and $L$ a finite set of lines in $\R^2$.  Let $$ I(P,L) = \{ (p,\ell) \in P \times L \,|\, p \in \ell \}$$ denote the set of {\em incidences} between $P$ and $L$. A basic question we will ask is how big can $I(P,L)$ be. The Szemeredi-Trotter (ST) theorem \cite{ST83} gives the (tight) upper bound of $$ |I(P,L)| \lesssim (|L|\cdot |P|)^{2/3} + |L| + |P|.$$ 

We begin this chapter in Section~\ref{sec-sttheorem} with two different proofs of this theorem. The first proof, presented in Section~\ref{sec-cellpart}, is due to Tao \cite{TaoCell} (based on \cite{CEGS90} and similar to the original proof of \cite{ST83}) and uses the method of cell partitions. The idea is to partition the two dimensional plane into cells, each containing a bounded number of points/lines and to argue about each cell separately. This uses the special `ordered' structure of the real numbers (this proof strategy is also the only one that generalizes to the complex numbers \cite{Toth}). The second proof, presented in Section~\ref{sec-crossing}, is due to Szekely \cite{Szekely} ands  uses the crossing number inequality for planar drawings of graphs and is perhaps the most elegant proof known for this theorem. This proof can also be adapted easily to handle intersections of more complex objects such as curves. We continue in Section~\ref{sec-applications} with some simple applications of the ST theorem to geometric and algebraic problems. These include  proving sum product estimates and counting distances between sets of points.

Sections \ref{sec-elekes} to \ref{sec-gkthm} are devoted to the proof of the Guth-Katz theorem on Erdos' distance counting problem. This theorem, obtained in \cite{GK10}, says that a set of $N$ points in the real plane define at least $\gtrsim N / \log N$ distinct distances. This gives an almost complete answer to an old question of Erdos (the upper bound has a factor of $\sqrt{\log N}$ instead of $\log N$). The tools used in the proof are developed over several sections which  contain several other related results.

In Section~\ref{sec-elekes} we discuss the Elekes-Sharir framework \cite{ES10} which reduces distance counting to a question about incidences of a specific family of lines in $\R^3$, much in the spirit of the ST theorem. Sections~\ref{sec-polynomial} and \ref{sec-polyhs} introduce the two main techniques used in the proof of the Guth-Katz theorem.  In Section~\ref{sec-polynomial} we introduce for the first time one of the main characters of this survey -- the polynomial method. As a first example to the power of this method, we show how it can be used to give a solution to another beautiful geometric conjecture -- the joints conjecture \cite{GK10b}. Here, we have a set of lines in $\R^3$ and want to upper bound the number of {\em joints}, or non-coplanar intersections of three lines or more. In Section~\ref{sec-polyhs} we introduce the second ingredient in the Guth-Katz theorem -- the polynomial Ham-Sandwich theorem. This technique, introduced by Guth in \cite{Guth08}, combines the polynomial method with the method of cell partitions. As an example of how this theorem is used we give a third proof of the ST theorem which was discovered recently \cite{KMS11}. 

Section~\ref{sec-gkthm}, contains a relatively detailed sketch of the  proof of the Guth-Katz theorem (omitting some of the more technical algebraic parts). The main result proved in this section is an incidence theorem upper bounding the number of pairwise intersections in a set of $N$ lines in $\R^3$. If we don't assume anything, $N$ lines can have $\gtrsim N^2$ intersections (an intersection is a pair of lines that intersect). An example is a set of $N/2$ horizontal lines and $N/2$ vertical lines, all lying in the same plane. If we assume, however, that the lines are `truly' in 3 dimensions, in the sense that no large subset of them lies in two dimensions, we can get a better (and tight) bound of $\leq N^{1.5} \log N$. This theorem then implies the bound on distinct distances using the Elekes-Sharir framework.

In the last Section of this chapter, Section~\ref{sec-sumprodgk}, we see yet another beautiful application of the three dimensional incidence theorem of Guth and Katz  obtaining optimal bounds in the flavor of the sum product theorem \cite{IRR11}. 
 
\subsection*{Chapter~\ref{chap-2}: Counting incidences over finite fields}

This chapter deals with the analog of the Szemeredi-Trotter theorem over finite fields and its applications. When we replace the field $\R$ with a finite field $\F_q$ of $q$ elements things become much more tricky and much less is known (in particular there are no tight bounds). Assuming nothing on the field, the best possible upper bound on the number of intersections between $N$ lines and $N$ points is $\sim N^{1.5}$, which is what one gets from only using the fact that two points determine a line (using a simple Cauchy-Schwarz calculation). However, if we assume that $\F_q$ does not contain large sub-fields (as is the case, for example, if $q$ is prime) one can obtain a small improvement of the form $N^{1.5 - \eps}$ for some positive $\eps$, provided $N \ll p^2$. This was shown by Bourgain, Katz and Tao as an application of the sum product theorem over finite fields \cite{BKT04}. The sum product theorem says that, under the same conditions on subfields, for every set $A \subset \F_q$ of size at most $q^{1-\alpha}$ we have $\max\{ |A+A|, |A\cdot A|\} > |A|^{1+ \alpha'}$, where $\alpha'$ depends only on $\alpha$. The set $A+A$ is defined as the set of all elements of the form $a+a'$ with $a,a' \in A$ ($A\cdot A$ is defined in a similar way). 

The proof of the finite field ST theorem is given in Sections \ref{sec-ruzsa} -- \ref{sec-stfinite}.  Section~\ref{sec-ruzsa} describes the machinery called `Ruzsa calculus' -- a set of useful claims for working with sumsets. Section~\ref{sec-growth} proves a theorem about growth of subsets of $\F_p$ (we will only deal with prime fields) which is a main ingredient of the proof of the ST theorem. Section~\ref{sec-bsg} proves the Balog-Szemeredi-Gowers theorem, a crucial tool in this proof and in many other results in additive combinatorics. Finally, Section~\ref{sec-stfinite} puts it all together and proves the final result. We note that, unlike previous expositions (and the original \cite{BKT04}), we opt to first prove the ST theorem and then derive the sum product theorem from it as an application. This choice allows us to derive a slightly more streamlined proof of the ST theorem.

As an application of these results over finite fields we will discuss, in Section~\ref{sec-multi}, the theory of {\em multi-source extractors} coming from theoretical computer science. We will see how to translate the finite field ST theorem into explicit mappings which transform `weak' structured sources of randomness into purely random bits. More precisely, suppose you are given samples from several (at least two) independent random variables and want to use them to output uniform random bits. It is not hard to show that a random function will do the job, but finding {\em explicit} (that is, efficiently computable) constructions is a difficult task. Such constructions have applications in theoretical computer science, in particular in the area of de-randomization, which studies the power of randomized computation vs. deterministic computation.

We will discuss in some detail two representative results in this area: the extractors of Barak, Impagliazzo and Wigderson for several independent blocks \cite{BIW06}, which were the first to introduce the tools of additive combinatorics to this area, and Bourgain's two source extractor \cite{Bour2source}.  Both rely crucially on the finite field Szemeredi-Trotter theorem of \cite{BKT04}.
 
\subsection*{Chapter~\ref{chap-3}: Packing lines in different directions --  Kakeya sets}

This chapter deals with a somewhat different type of theorems that describe  the way lines in different directions can overlap. In Sections~\ref{sec-kakeyareal} and \ref{sec-kakeyafin} we will discuss these questions over the real numbers and over finite fields, respectively. In Section~\ref{sec-mergers} we will discuss applications of the finite field results to problems in theoretical computer science.

A Kakeya set $K \subset \R^n$ is a compact set containing a unit line segment in every direction. These sets can have measure zero. An important open problem is to understand the minimum Minkowski or Hausdorff dimension\footnote{For a definition see Section~\ref{sec-kakeyareal}.} of a Kakeya set. This question reduces in a natural way to a discrete incidence question involving a finite set of lines in many `sufficiently separated' directions. The Kakeya conjecture states that Kakeya sets must have maximal dimension (i.e., have dimension $n$). The conjecture is open in dimensions $n \geq 3$ and was shown to have deep connections with other problems in Analysis, Number Theory and PDE's  (see \cite{Tao01}). 

The most successful line of attack on this conjecture was initiated by Bourgain \cite{Bou99} and later developed by Katz and Tao \cite{KT02} and uses tools from additive combinatorics. In Section~\ref{sec-kakeyareal} we will discuss Kakeya sets over the reals and prove a $\geq (4/7)n$ bound on the Minkowski dimension, which is very close to the best known lower bound of $(0.596...)n$. The underlying additive combinatorics problem that arises in this context is upper bounding the number of differences $a-b$, for pairs $(a,b) \in G \subset A \times B$ in some graph $G$ as a function of the number of sums (or, more generally, linear combinations) on the same graph. We will not discuss the applications of the Euclidean Kakeya conjecture since they are out of scope for this survey (we are focusing on applications in discrete mathematics and computer science). Even though we will not directly use additive combinatorics results developed in Chapter~\ref{chap-2}, they will be in the background and will provide intuition as to what is going on.

Over a finite field $\F_q$ a Kakeya set is a set containing a line in every direction (a line will contain $q$ points). It was conjectured by Wolff \cite{Wolff99} that the minimum size of a Kakeya set is at least $C_n \cdot q^n$ for some constant $C_n$ depending only on $n$. We will see the proof of this conjecture (obtained by the author in \cite{Dvir08}) which uses  the polynomial method.  An application of this result, described in Section~\ref{sec-mergers}, is a construction of {\em seeded extractors}, which are explicit mappings that transform a `weak' random source into a close-to-uniform distribution with the aid of a short random `seed' (since there is a single source, the extractor must use a seed). A specific question that arises in this setting is the following: Suppose Alice and Bob each pick a point $X,Y \in \F_q^n$ ($X$ for Alice, $Y$ for Bob). Consider the random variable $Z$ computed by picking a random point on the line through $X,Y$. If both Alice and Bob pick their points independently at random then it is easy to see that $Z$ will also be random. But what happens when Bob picks his points $Y$ to be some function $Y = F(X)$? Using the connection to the Kakeya conjecture one can show that, in this case, $Z$ is still sufficiently random in the sense that it cannot hit any small set with high probability.

\subsection*{Chapter~\ref{chap-4}: From local to global --  Sylvester-Gallai type theorems}


The Sylvester-Gallai (SG) theorem says that, in a finite set of points in $\R^n$, not all on the same line, there exists a line intersecting exactly two of the points. In other words, if for every two points $u,v$ in the set, the line through $u,v$ contains a third point in the set, then all points are on the same line. Besides being a natural incidence theorem, one can also look at this theorem as converting local geometric information (collinear triples) into global upper bounds on the dimension (i.e., putting all points on a single line, which is one dimensional). We will see several generalizations of this theorem, obtained in \cite{BDWY11}, in various settings. For example, assume that for every point $u$ in a set of $N$ points there are at least $N/100$ other points $v$ such that the line through $u,v$ contains a third point. We will see in this case that the points all lie on an affine subspace of dimension bounded by a constant. The proof technique here is different than what we have seen so far and will rely on convex optimization techniques among other things. These results will be described in Section~\ref{sec-sgreal} with the main technical tool, a rank lower bound for design matrices, proved in Section~\ref{sec-rankdesign}.
 
In Section~\ref{sec-sgfinite} we will  consider this type of questions over a finite field and see how the bounds are weaker in this case. In particular, under the same assumption as above (with $N/100$) the best possible upper bound on the dimension will be $\lesssim \log_q(N))$, where $q$ is the characteristic of the field \cite{BDSS11}. Here, we will again rely on tools from additive combinatorics and will use results proved in Chapter~\ref{chap-2}.
 
In Section~\ref{sec-lcc} we will see how this type of  questions arise naturally in computer science applications involving error correcting codes which are `locally correctable'. A (linear){\em Locally-Correctable-Code} (LCC) is a (linear) error correcting code in which each symbol of a possible corrupted codeword can be corrected by looking at only a few other locations (in the same corrupted codeword). Such codes are very different than `regular' error correcting codes (in which decoding is usually done in one shot for all symbols) and have interesting applications in complexity theory\footnote{They are also very much related to Locally {\em Decodable} Codes (LDCs) which are discussed at length in the survey \cite{Y_now}.}.

\ifediting
\bibliographystyle{alpha}
 \bibliography{incidence}

\end{document}
\fi

\chapter{Counting Incidences Over the Reals}\label{chap-1}

\newif\ifediting

\ifediting
\documentclass[11pt]{article}

\usepackage{amsmath,amsthm,amssymb}
\newcommand{\remove}[1]{}
\setlength{\topmargin}{0.3in} \setlength{\headheight}{0in}
\setlength{\headsep}{0in} \setlength{\textheight}{8.0in}
\setlength{\topsep}{0.1in} \setlength{\itemsep}{0.0in}
\parskip=0.05in
 \textwidth=6.5in 
\oddsidemargin=0truecm \evensidemargin=0truecm

\newtheorem{thm}{Theorem}[section]
\newtheorem{claim}[thm]{Claim}
\newtheorem{lem}[thm]{Lemma}
\newtheorem{define}[thm]{Definition}
\newtheorem{cor}[thm]{Corollary}
\newtheorem{obs}[thm]{Observation}
\newtheorem{example}[thm]{Example}
\newtheorem{construct}[thm]{Construction}
\newtheorem{conjecture}[thm]{Conjecture}
\newtheorem{THM}{Theorem}
\newtheorem{question}{Question}
\newtheorem{fact}[thm]{Fact}
\newtheorem{prop}[thm]{Proposition}

\def\F{{\mathbb{F}}}
\def\Q{{\mathbb{Q}}}
\def\Z{{\mathbb{Z}}}
\def\N{{\mathbb{N}}}
\def\R{{\mathbb{R}}}
\def\K{{\mathbb{K}}}
\def\C{{\mathbb{C}}}
\def\A{{\mathbb{A}}}
\def\P{{\mathbb{P}}}
\def\cP{{\cal P}}
\def\cS{{\mathcal S}}
\def\cE{{\mathcal E}}
\def\V{{\mathbf{V}}}
\def\I{{\mathbf{I}}}
\def\bx{{\mathbf x}}
\def\by{{\mathbf y}}
\def\E{{\mathbb E}}

\def\half{ \frac{1}{2}}
\newcommand{\ip}[2]{\langle #1,#2 \rangle}
\def\sumN{\sum_{i=1}^n}
\def\_{\,\,\,\,\,}
\def\prob{{\mathbf{Pr}}}
\newcommand{\entropy}[1]{ {\text{H}_{\infty}\left({#1}\right)} }
\def\modulo{\text{mod}}
\def\omm{ \{0,1\} }
\def\id{ \textit{id} }

\def\D{{\partial}}
\def\gap{\textsf{gap}}
\def\sign{\textsf{sign}}
\def\spar{\textsf{sparse}}
\def\span{\textsf{span}}
\def\Part{\textbf{Part}}
\def\Mon{\textbf{Mon}}
\def\sing{\textbf{sing}}
\def\Und{\textsf{Und}}
\def\Comp{\textsf{Comp}}
\def\rank{\textsf{rank}}
\def\poly{\textsf{poly}}
\def\codim{\textsf{codim}}
\def\cp{\textsf{cp}}
\def\uni{\textsf{Uni}}
\def\ext{\textsf{\bf Ext}}
\def\extt{\textsf{\bf Ext2}}

\def\fplus{{\,+_f\,}}

\newcommand{\epsclose}{\stackrel{\epsilon}{\thicksim} }
\newcommand{\eclose}[1]{\stackrel{{#1}}{\thicksim} }
\newcommand{\eps}{\epsilon}
\newcommand{\Anote}[1]{\begin{quote}{\sf Avi's Note:} {\sl{#1}} \end{quote}}

\begin{document}

\title{Incidence Theorems -- Lecture Notes}
\date{}
\maketitle


\fi

\section{The Szemeredi-Trotter theorem}\label{sec-sttheorem}

Let $L$ be a finite set of points in $\R^2$ and let $P$ be a finite set of points in $\R^2$. We define $$ I(P,L) = \{ (p,\ell) \in P\times L \,|\, p \in \ell\}$$ to be the set of incidences between $P$ and $L$. We will prove the following result of Szemeredi and Trotter \cite{ST83}.

\begin{thm}[ST theorem]
Under the above notations we have $$|I(P,L)| \lesssim (|P|\cdot |L|)^{2/3} + |L| + |P|.$$
\end{thm}
We will use $\lesssim, \gtrsim$ and $\sim$ to denote (in)equality up to multiplicative absolute constants.

The following example shows that this bound is tight. Let $L$ be the set of $N \sim M^3$ lines of the form $\{ (x,y) \in \R^2 \,|\, y = ax + b, \}$ with $a\in [M], b \in [M^2] $. Let $P = \{ (x,y) \in \R^2 \,|\, x \in [M], y \in [2M^2]\}$ be a set of $N \sim M^3$ points. Observe that each line $\ell \in L$ intersects $P$ in $ \geq M$ points (for each $x \in [M]$, $y =ax+ b \leq 2M^2$). This gives a total of $M^4 \sim N^{4/3}$ incidences. 

As a step towards proving the ST theorem we prove the following  claim which gives an `easy' bound on the number of incidences. It is `easy' not just because it has a simple proof but also because it only uses the fact that every two points define a single line and every pair of lines can intersect in at most one point (these facts hold over any field). The proof of the claim will use the Cauchy-Schwarz inequality which says that $$ \left( \sum_{i=1}^k a_i \cdot b_i \right)^2 \leq \left( \sum_{i=1}^k a_i^2 \right)\cdot \left( \sum_{i=1}^k b_i^2 \right)$$ whenever $a_i,b_i$ are positive real numbers. 
\begin{claim}\label{cla-CS}
Let $P,L$ be as above. Then we have the following two bounds:
$$ I(P,L) \lesssim |P| \cdot |L|^{1/2} + |L| $$ and
$$ I(P,L) \lesssim |L| \cdot |P|^{1/2} + |P|.$$
\end{claim}
\begin{proof}
We will only prove the first assertion (the second one follows using a similar argument or by duality). The only geometric property used is that through every two points passes only one line. First, observe that 
\begin{equation}\label{eq-P2L}
	|I(P,L)| \leq |P|^2 + |L|.
\end{equation}
To see this, count first the lines that have at most one point in $P $ on them. These lines contribute at most $|L|$ incidences. The rest of the lines have at least two points in $P$ on each line. The total number of incidences on these lines is at most $|P|^2$ since otherwise there would be a point $p \in P$ that lies on $> |P|$ lines and each of these lines must have one additional point on it and so there are more than $|P|$ points -- a contradiction.

We now bound the number of incidences. We use $1_{p \in \ell}$ to denote the indicator function which is equal to $1$ if $p \in \ell$ and equal to zero otherwise.
\begin{eqnarray}
|I(P,L)|^2 &=& \left(\sum_{\ell \in L} \sum_{p \in P} 1_{p \in \ell}\right)^2 \\
&\leq& |L| \cdot  \sum_{\ell \in L} \left( \sum_{p \in P} 1_{p \in \ell}\right)^2 \text{(Cauchy Schwarz)}\\
&=& |L| \cdot \sum_{p_1,p_2 \in P} \sum_{\ell \in L} 1_{p_1 \in \ell}\cdot 1_{p_2 \in \ell} \\
&\leq& |L| \cdot \left( |I(P,L)| + |P|^2 \right) \\
&\leq& |L|^2 + 2|L| \cdot |P|^2,
\end{eqnarray}
which implies the bound.
\end{proof}

\subsection{Proof using cell partitions}\label{sec-cellpart}

The first proof of the ST theorem we will see uses the idea of cell partitions and is perhaps the most direct of the three proofs we will encounter. The proof we will see is due to Tao \cite{TaoCell} (based loosely on \cite{CEGS90}) and is similar in spirit to the original proof of Szemeredi and Trotter. The idea is to use the properties of the real plane to partition it into small regions such that each region will intersect a small fraction of the lines in our set $L$. This allows to `amplify' the easy bound (Claim~\ref{cla-CS}) to a stronger (indeed, optimal) bound by applying it to separated smaller instances of the problem.   

\begin{lem}\label{lem-partition}
For every $r \geq 1$ there exists a set of $O(r)$ lines plus some additional line segments not containing any point in $P$ that partition $\R^2$ into at most $O(r^2)$ regions (convex open sets) such that the interior of each region is incident to at most $O(|L|/r)$ lines in $L$.
\end{lem}

We will sketch the proof of this lemma later. Before that, let's see how it implies the ST theorem: First we can assume w.l.o.g that $$ |L|^{1/2} \ll |P| \ll |L|^2$$ (we use $A \ll B$ to mean that $A \leq c \cdot B$ for some sufficiently small constant $c$). If not, then the bound in the ST theorem follows from Claim~\ref{cla-CS}. We will apply Lemma~\ref{lem-partition} with some $r$ to be chosen later. Let $R$ be the set of lines defining the partition (recall that there are some additional line segments not counted in $R$ that do not contain points in $P$). For each cell $C$ we apply Claim~\ref{cla-CS} to bound the number of incidences in this cell (the cell does not include the boundary). We get that a cell $C$ can have at most $O( |P \cap C| \cdot (|L|/r)^{1/2}  + |L|/r)$ incidences. Summing over all cells we get that $$ |I(P,L)| \leq |I(P,L \cap R)| + O( |P||L|^{1/2}/r^{1/2} + |L|r) + O(|L|r)$$ where the first term counts the incidences of point with lines in $R \cap L$, the second term counts incidences in the open cells and the third term counts the incidences of lines not in $R$ with points in the cell boundary (each line not in $R$ has at most $r$ incidences with points on $R$). Setting $$ r \sim |P|^{2/3}/|L|^{1/3}$$ we get that $$ |I(P,L)| \lesssim |I(P,L \cap R)| + |P|^{2/3}|L|^{2/3}.$$ Since $|P| \ll |L|^2$ we get that $r \leq |L|/10$ and so, we can repeat the same argument on $P, L \cap R$ obtaining a geometric sum that only adds up to a constant. This completes the proof of the ST theorem.

\subsubsection*{Proof of the cell partition lemma}
We only sketch the proof. The proof will be probabilistic. We will pick a random set of the lines in $L$ to be the set $R$ (plus some additional segments) and will argue that it satisfies the lemma with positive probability (this will imply that a good choice exists). This type of argument is common in combinatorics and is usually referred to as the `probabilistic method'. We will make two simplifying assumptions: one is that at most two lines pass through a point (this can be removed by a limiting argument). The second is that there are no vertical lines in $L$ and that no point in $P$ is on a vertical line through the intersection of two lines in $L$ (this can be removed by a random rotation).

The particular procedure we will use to pick the partition is the following: first we take each line $\ell \in L$ to be in $R$ with probability $r/|L|$. This will give us $O(r)$ lines with high probability (say, at least $0.99$). This set of lines can create at most $O(r^2)$ cells. Then, we `fix' the partition so that each cell has a bounded (at most 4) number of line segments bordering it. This `fix' is done by adding vertical line segments through every point that is adjacent to a cell with more than 4 border segments (the number 4 is not important, it can be any constant). These extra line segments  are not in $L$ and, by our `random rotation' assumption, do not hit any point in $P$. One can verify that adding these segments does not increase the number of cells above $O(r^2)$ (there are at most $O(r^2)$ initial `corners' to fix). 

Having described the probabilistic construction we turn to analyze the probability of a cell having too many lines passing through it. Consider a cell with $4$ border segments. Each line passing through the cell must intersect at least one of these bordering segments. If there are more than $M$ lines in the cell than one segment must have at least $M/2$ lines in $L$ passing through it. Since all of these lines were not chosen in the partition we get that this event (for this specific segment) happens with probability at most $$ (1 - r/|L|)^{M/2}.$$ Taking $M$ to be roughly $\sim 100 |L|\log|L|/r$ we get that this probability is at most $|L|^{-100}$. Therefore, we can bound the union of all `bad' events of this form (i.e, of a particular segment or a line in $L$ containing a series of $M$ lines not chosen in $R$) as the product of the number of events times $|L|^{-100}$. Since the number of bad events is much smaller than $|L|^{100}$ we get that there exists a partition with a bound of $O(|L|\log|L| /r)$. A more careful argument can get rid of the logarithmic factor by arguing that (a) the number of `bad' cells is very small and (b) we can use induction on this smaller set to get the required partition. 

This proof seems messy but is actually much cleaner than the original partition proof of Szemeredi and Trotter (which was deterministic). Next, we will see a much simpler proof of ST that does not use cell partition (later on we will see a third proof that uses a very different kind of cell partition using polynomials).

\subsection{Proof using the crossing number inequality}\label{sec-crossing}
Next, we will see a different, very elegant, proof of the ST theorem due to Szekely \cite{Szekely} based on the powerful crossing number inequality \cite{ACNS,Leighton}. We will consider undirected graphs $G = (V,E)$ on a finite set $V$ of vertices and with a set $E \subset V \times V$ of edges. A {\em drawing} of a graph is a placing of the vertices in the real plane $\R^2$ with simple curves connecting two vertices if there is an edge between them (we omit the `formal' definition since this is a very intuitive notion). For a drawing $D$ of $G$ we denote by $\cros(D)$ the number of `crossings' or intersections of edges in the drawing. The crossing number of $G$, denoted $\cros(G)$ is the minimum over all drawings $D$ of $G$ of $\cros(D)$. Thus, a graph is planar if it has a crossing number of zero.

A useful tool when talking about planar graphs is Euler's formula. Given a planar drawing $D$ of a connected graph $G = (V,E)$ we have the following equality
\begin{equation}
	|V| - |E| + |F| = 2,
\end{equation}
where $F$ is the set of faces of the drawing (including the unbounded face). The proof is a very simple induction on $|F|$. If there is one face then the graph is a tree and so $|V| = |E|+1$. If there are more faces then we can remove a single edge and decrease the number of faces by one.

The proof of the ST theorem will use a powerful inequality called the crossing number inequality. This inequality gives a strong lower bound on $\cros(G)$ given the number of edges in $G$. As a preliminary step we shall prove a weaker bound (which we will amplify later).
\begin{claim}\label{cla-weak}
Let $G$ be a graph. Then $\cros(G) > |E| - 3|V|.$
\end{claim}
\begin{proof}
W.l.o.g we can assume $G$ is connected and with at least 3 vertices. It is easy to check that, if $G$ is planar then $3|F| \leq 2|E|$ (draw two points on either side of an edge and count them once by going over all edges and another  by going over all faces, using the fact that each face has at least 3 edges adjacent to it).  Plugging this into Euler's formula we get that, for planar graphs $$ 2 \leq |V| - (1/3)|E|$$ or $|E| < 3|V|$. If the claim was false, we could remove at most $\cros(G) \leq |E| - 3|V|$ edges and obtain a planar drawing of $G$. The new graph will have at least $|E| - (|E| - 3|V|) = 3|V|$ edges -- a contradiction.
\end{proof}

This is clearly not a very good bound and some simple examples demonstrate. To get the final bound we will apply Claim~\ref{cla-weak} on a {\em random} vertex induced subgraph and do some expectation analysis. This is a beautiful example of the power of the probabilistic method.

\begin{thm}[Crossing Number Inequality \cite{Leighton,ACNS}]
Let $G$ be a graph. If $|E| \geq 4|V|$ then $$\cros(G) \geq \frac{|E|^3}{64|V|^2}.$$
\end{thm}
\begin{proof}
Let $G' = (V',E')$ be a random vertex induced subgraph with each vertex of $V$ chosen to be in $V'$ independently with probability $p \in [0,1]$ to be chosen later. Taking the expectation of  $$\cros(G') \geq |E'| - 3|V'|$$  we get that $$ p^4 \cdot \cros(G) \geq p^2 \cdot |E| - 3p \cdot |V|.$$ The right hand side is equal to the expectation of the r.h.s of the original inequality by linearity of expectation. The left hand side requires some explanation: consider a single crossing in a drawing of $D$ which has the smallest number of crossings. This crossing will remain after the random restriction with probability $p^4$. Thus, the expected number of crossings will be $p^4\cros(G)$. This is however only an upper bound on the expected $\cros(G')$ since there could be new ways of obtaining even better drawing after we move to $G'$ (but this inequality is the `right' direction so we're fine). Setting $p = 4|V|/|E|$ (which is at most 1 by our assumption) gives the required bound.
\end{proof}

We now prove the ST theorem using this inequality (this proof is by Szekely \cite{Szekely}): Let $P,L$ be our finite sets of points and lines (as above). Put aside lines that have only 1 point on them (this can contribute $|L|$ incidences) so that every line has at least two points. Consider the drawing of the graph whose vertex set is $P$ and two points share an edge if they are (1) on the same line and (2) there is no third point on the line segment connecting them. The number of edges on a  line $\ell \in L$ is $|P_\ell| - 1$, where $P_\ell = P \cap \ell$. The total number of edges is $$ \sum_{\ell} |P_\ell| - 1 \geq (1/2)|I(P,L)|.$$ The crossing number of this graph is at most $|L|^2$ since each crossing is obtained from the intersection of two lines in $|L|$. Plugging all this into the crossing number inequality we get that either $|I(P,L)| \lesssim |P|$ or that $$ |L|^2 \gtrsim \frac{|I(P,L)|^3}{|P|^2}.$$ Putting all of this together we get $|I(P,L)| \lesssim (|P||L|)^{2/3} + |P| + |L|$.

Notice that this proof can be made to work with simple curves instead of lines as long as two curves intersect in at most $O(1)$ points and two points sit on at most $O(1)$ curves together. In particular, a set $P$ of points and a set $C$ of unit circles can have at most $\lesssim (|C||P|)^{2/3} + |C| + |P|$ incidences (we will use this fact later).

\ifediting
\bibliographystyle{alpha}
 \bibliography{incidence}

\end{document}
\fi

\newif\ifediting

\ifediting
\documentclass[11pt]{article}

\usepackage{amsmath,amsthm,amssymb}
\newcommand{\remove}[1]{}
\setlength{\topmargin}{0.3in} \setlength{\headheight}{0in}
\setlength{\headsep}{0in} \setlength{\textheight}{8.0in}
\setlength{\topsep}{0.1in} \setlength{\itemsep}{0.0in}
\parskip=0.05in
 \textwidth=6.5in 
\oddsidemargin=0truecm \evensidemargin=0truecm

\newtheorem{thm}{Theorem}[section]
\newtheorem{claim}[thm]{Claim}
\newtheorem{lem}[thm]{Lemma}
\newtheorem{define}[thm]{Definition}
\newtheorem{cor}[thm]{Corollary}
\newtheorem{obs}[thm]{Observation}
\newtheorem{example}[thm]{Example}
\newtheorem{construct}[thm]{Construction}
\newtheorem{conjecture}[thm]{Conjecture}
\newtheorem{THM}{Theorem}
\newtheorem{question}{Question}
\newtheorem{fact}[thm]{Fact}
\newtheorem{prop}[thm]{Proposition}

\def\F{{\mathbb{F}}}
\def\Q{{\mathbb{Q}}}
\def\Z{{\mathbb{Z}}}
\def\N{{\mathbb{N}}}
\def\R{{\mathbb{R}}}
\def\K{{\mathbb{K}}}
\def\C{{\mathbb{C}}}
\def\A{{\mathbb{A}}}
\def\P{{\mathbb{P}}}
\def\cP{{\cal P}}
\def\cS{{\mathcal S}}
\def\cE{{\mathcal E}}
\def\V{{\mathbf{V}}}
\def\I{{\mathbf{I}}}
\def\bx{{\mathbf x}}
\def\by{{\mathbf y}}
\def\E{{\mathbb E}}

\def\half{ \frac{1}{2}}
\newcommand{\ip}[2]{\langle #1,#2 \rangle}
\def\sumN{\sum_{i=1}^n}
\def\_{\,\,\,\,\,}
\def\prob{{\mathbf{Pr}}}
\newcommand{\entropy}[1]{ {\text{H}_{\infty}\left({#1}\right)} }
\def\modulo{\text{mod}}
\def\omm{ \{0,1\} }
\def\id{ \textit{id} }

\def\D{{\partial}}
\def\gap{\textsf{gap}}
\def\sign{\textsf{sign}}
\def\spar{\textsf{sparse}}
\def\span{\textsf{span}}
\def\Part{\textbf{Part}}
\def\Mon{\textbf{Mon}}
\def\sing{\textbf{sing}}
\def\Und{\textsf{Und}}
\def\Comp{\textsf{Comp}}
\def\rank{\textsf{rank}}
\def\poly{\textsf{poly}}
\def\codim{\textsf{codim}}
\def\cp{\textsf{cp}}
\def\uni{\textsf{Uni}}
\def\ext{\textsf{\bf Ext}}
\def\extt{\textsf{\bf Ext2}}

\def\fplus{{\,+_f\,}}

\newcommand{\epsclose}{\stackrel{\epsilon}{\thicksim} }
\newcommand{\eclose}[1]{\stackrel{{#1}}{\thicksim} }
\newcommand{\eps}{\epsilon}
\newcommand{\Anote}[1]{\begin{quote}{\sf Avi's Note:} {\sl{#1}} \end{quote}}

\begin{document}

\title{Incidence Theorems -- Lecture Notes}
\date{}
\maketitle


\fi

\section{Applications of Szemeredi-Trotter over $\R$}\label{sec-applications}

We mentioned that the crossing number inequality proof of the Szemeredi-Trotter theorem works also for circles of unit distance. In general, the following is also true and very useful (the proof is left as an easy exercise using the crossing number inequality).
\begin{thm} Suppose we have a family $\Gamma$ of simple curves (i.e., that do not intersect themselves) and a set of points $P$ in $\R^2$ such that (1) every two points define at most $C$ curves in $\Gamma$ and (2) every two curves in $\Gamma$ intersect in at most $C$ points then $$|I(\Gamma,P)| \lesssim (|\Gamma|\cdot |P|)^{2/3} + |\Gamma| + |P|,$$ where the hidden constant in the inequality depends only on $C$.
\end{thm}

If we only have a bound on the number of curves passing through $k$ of the points (for some integer $k \geq 2$) the following was shown by Pach and Sharir (and is not known to be tight for $k \geq 3$):
\begin{thm}[Pach-Sharir \cite{PS98}]
Let $\Gamma$  be a family of curves in the plane and let $P$ be a set of points. Suppose that through every $k$ points of $P$ there are at most $C$ curves and that every two curves intersect in at most $C$ points. Then
$$ |I(\Gamma,P)| \lesssim |P|^{k/(2k-1)} \cdot |\Gamma|^{1 - 1/(2k-1)} + |\Gamma| + |P|, $$ where the hidden constant depends only on $C$.
\end{thm}
In particular, this bound can be used for families of algebraic curves of bounded degree. 
 
A simple but useful corollary of the ST theorem is the following. The proof is an easy calculation and is left to the reader.
\begin{cor}\label{cor-rich}
The $P$ and $L$ be sets of points and lines in $\R^2$. For $k > 1$ let $L_k$ denote the set of lines in $L$ that contain at least $k$ elements of $P$. Then,
$$|L_k| \lesssim \frac{|P|^2}{k^3} + \frac{|P|}{k}. $$
\end{cor}

\subsection{Beck's theorem}
A nice theorem that follows from ST using purely combinatorial arguments is Beck's theorem: 
\begin{thm}[Beck's theorem \cite{Beck}]
Let $P$ be a set of points in the plane and let $L$ be the set of lines containing at least 2 points in $P$. Then one of these two cases must hold:
\begin{enumerate}
	\item There exists a line in $L$ that contains $\gtrsim |P|$ points.
	\item $|L| \gtrsim |P|^2$
\end{enumerate}
In other words, if a set of points $P$ defines $\ll |P|^2$ lines than there is a line containing a constant fraction of the points.
\end{thm}
\begin{proof}
Let $n = |P|$. We partition the lines in $L$ into $\sim \log n$ sets $L_j \subset L$, the $j$'th set contains the lines with $\sim 2^j$ points from $P$ (more precisely, the lines that contain between $2^j$ and $2^{j+1}$ points).  Using Corollary~\ref{cor-rich} on each of the $L_j$'s  we get that $$ |L_j| \lesssim \frac{n^2}{2^{3j}} + \frac{n}{2^j}.$$ Take $C$ to be some large constant to be chosen later and let $M = \cup_{C \leq 2^j \leq n/C} L_j$ be the set of lines with `medium' multiplicity. Since each line in $L_j$ contains $\sim 2^{2j}$ pairs of points we get that there are at most $\lesssim \frac{n^2}{2^j} + 2^jn$ pairs of points on lines in $L_j$. Summing over all $j$'s with $C \leq 2^j \leq n/C$ and taking $C$ to be sufficiently large we get that the number of pairs of points on  lines in $M$ is at most $ n^2/100$. Thus, there are two cases: either there is a line with at least $n/C$ points and we are done. Alternatively, there are $\gtrsim n^2$ points on lines that contain at most $C$ points each. In this case we get that $|L| \gtrsim n/C^2$ which is linear in $n$. 
\end{proof}

\subsection{Erdos unit distance and distance counting problems}

Let $A$ be a set of points in $\R^2$. We define $$ \Delta_1(A) = \{(p,q) \in A^2 \,|\, \|p - q\| = 1 \},$$ (i.e., all pairs of Euclidean distance 1). Erdos conjectured (and this is still open) that for all sets $A$ we have $|\Delta_1(A)| \leq C(\eps) \cdot |A|^{1+\eps}$ for all $\eps > 0$. Again, the construction which obtains the maximal number of unit distances is a grid (this time, however, the step size must be chosen carefully using number theoretic properties).

Using the unit-circle version of the ST theorem  we can get a bound of $|\Delta_1(A)| \lesssim |A|^{4/3}$, which is the best bound known. To see this, consider the $|A|$ circles of radius 1 centered at the points of $|A|$. Two circles can intersect in at most two points and every two points define at most two radius one circles that pass through them. Therefore, we can use the ST theorem to bound the number of incidences by $|A|^{4/3}$. In four dimensions the number of unit distances in an arrangement can be as high as $\sim |A|^2$. In three dimensions the answer is not known.

A related question of Erdos is to lower bound the number of distinct distances defined by the pairs in $A$. Let $$\dist(A) = \{ \| p - q\| \,|\, p,q \in A \}.$$  It was conjectured by Erdos that $$ |\dist(A)| \gtrsim \frac{|A|}{(\log |A|)^{1/2}}. $$ Considering $n$ points on a $\sqrt n \times \sqrt n$ integer grid, gives an example showing this bound is tight.  Using the bound on unit distances above (which, by scaling, bounds the maximal number of distances equal to any real number) we immediately get a lower bound of $\gtrsim |A|^{2/3}$ on the number of distinct distances. A result which comes incredibly close to proving Erdos's conjecture  (with $\log$ instead of $\log^{1/2}$) is a recent breakthrough of Guth and Katz which uses a three dimensional incidence theorem in the spirit of the ST theorem (we will see this proof later on).

\subsection{Sum Product theorem over $\R$}
Let $A \subset \R$ be a finite set and define $$ A+A = \{ a+a' \,|\, a,a' \in A \}$$ and $$A \cdot A = \{ a\cdot a' \,|\, a,a' \in A \}$$ to be the sumset and product set of $A$.  If $A$ is an arithmetic progression than we have $|A+A| \lesssim |A|$ and if $A$ is a geometric progression we have $|A\cdot A| \lesssim |A|$. But can $A$ be both? In other words, can we have an `approximate' sub-ring $A$ of the real numbers (one can also ask this for the integers). Using the ST theorem, Elekes \cite{Elekes} proved the following theorem. The same theorem with the constant $5/4$ replaced by some other constant larger than 1 was proved earlier by Erdos and Szemeredi  \cite{ES83}. 
\begin{thm}[Sum-Product Theorem over $\R$]
Let $A \subset \R$ be a finite set. Then $$ \max\{ |A+A|,|A\cdot A| \} \gtrsim |A|^{5/4}.$$
\end{thm}
\begin{proof}
Let $P = (A\cdot A) \times (A+A)$ and let $L$ contain all lines defined by an equation of the form $y = ax+ b$ with $a \in A^{-1}$ and $b \in A$ ($A^{-1}$ is the set of inverses of elements of $A$ (we can discard the zero). Then $|L| = |A|^2$. Each line in $L$ has $\geq |A|$ incidences with $P$ (take all $(x,y)$ on the line with $x = a^{-1}\cdot x'$ for some $x' \in A$) and so we have $$ |A|^3 \lesssim |A\cdot A |^{2/3} \cdot |A+A|^{2/3} \cdot |A|^{4/3}.$$ If both $|A\cdot A|$ and $|A+A|$ were $\ll |A|^{5/4}$ we would get that the r.h.s is smaller than the l.h.s -- a contradiction. 
\end{proof}
 
A more intricate application of the ST theorem can give an improved bound of $ \max\{ |A+A|,|A\cdot A| \} \gtrsim |A|^{14/11}$ \cite{Solymosi05}. The conjectured bound is $\sim |A|^{2-\eps}$ for all $\eps > 0$.

\subsection{Number of points on a convex curve}
   
Let $\gamma$ be a strictly convex curve in the plane contained in the range $[n]\times [n]$. How many integer lattice points can $\gamma$ have? Using the curve version of the ST theorem we can bound this by $\lesssim n^{2/3}$ (this proof is due to Iosevich \cite{Iosevich04}). This bound is tight and the example which matches it is the convex hull of integer points contained in a ball of radius $n$ \cite{BL98}. The argument proceeds as follows:  take the family $\Gamma$ to include all curves obtained from $\gamma$ by translating it by all integer points in $[n] \times [n]$. One has $|\Gamma|= n^2$. We take $P$ to be all integer points that are on some curve from $\Gamma$ so that $|P| \lesssim n^2$. The condition on the number of curves through every two points and the maximum intersection of two curves can be readily verified (here we must used the strict convexity). Thus the number of incidences is at most $n^{8/3}$. Since the curves are all translations of each other they all contain the same number of integer points. Therefore, each one will contain at most $n^{2/3}$ points.

\ifediting
\bibliographystyle{alpha}
 \bibliography{incidence}

\end{document}
\fi

\newif\ifediting

\ifediting
\documentclass[11pt]{article}

\usepackage{amsmath,amsthm,amssymb}
\newcommand{\remove}[1]{}
\setlength{\topmargin}{0.3in} \setlength{\headheight}{0in}
\setlength{\headsep}{0in} \setlength{\textheight}{8.0in}
\setlength{\topsep}{0.1in} \setlength{\itemsep}{0.0in}
\parskip=0.05in
 \textwidth=6.5in 
\oddsidemargin=0truecm \evensidemargin=0truecm

\newtheorem{thm}{Theorem}[section]
\newtheorem{claim}[thm]{Claim}
\newtheorem{lem}[thm]{Lemma}
\newtheorem{define}[thm]{Definition}
\newtheorem{cor}[thm]{Corollary}
\newtheorem{obs}[thm]{Observation}
\newtheorem{example}[thm]{Example}
\newtheorem{construct}[thm]{Construction}
\newtheorem{conjecture}[thm]{Conjecture}
\newtheorem{THM}{Theorem}
\newtheorem{question}{Question}
\newtheorem{fact}[thm]{Fact}
\newtheorem{prop}[thm]{Proposition}

\def\F{{\mathbb{F}}}
\def\Q{{\mathbb{Q}}}
\def\Z{{\mathbb{Z}}}
\def\N{{\mathbb{N}}}
\def\R{{\mathbb{R}}}
\def\K{{\mathbb{K}}}
\def\C{{\mathbb{C}}}
\def\A{{\mathbb{A}}}
\def\P{{\mathbb{P}}}
\def\cP{{\cal P}}
\def\cS{{\mathcal S}}
\def\cE{{\mathcal E}}
\def\V{{\mathbf{V}}}
\def\I{{\mathbf{I}}}
\def\bx{{\mathbf x}}
\def\by{{\mathbf y}}
\def\E{{\mathbb E}}

\def\half{ \frac{1}{2}}
\newcommand{\ip}[2]{\langle #1,#2 \rangle}
\def\sumN{\sum_{i=1}^n}
\def\_{\,\,\,\,\,}
\def\prob{{\mathbf{Pr}}}
\newcommand{\entropy}[1]{ {\text{H}_{\infty}\left({#1}\right)} }
\def\modulo{\text{mod}}
\def\omm{ \{0,1\} }
\def\id{ \textit{id} }

\def\D{{\partial}}
\def\gap{\textsf{gap}}
\def\sign{\textsf{sign}}
\def\spar{\textsf{sparse}}
\def\span{\textsf{span}}
\def\Part{\textbf{Part}}
\def\Mon{\textbf{Mon}}
\def\sing{\textbf{sing}}
\def\Und{\textsf{Und}}
\def\Comp{\textsf{Comp}}
\def\rank{\textsf{rank}}
\def\poly{\textsf{poly}}
\def\codim{\textsf{codim}}
\def\cp{\textsf{cp}}
\def\uni{\textsf{Uni}}
\def\ext{\textsf{\bf Ext}}
\def\extt{\textsf{\bf Ext2}}

\def\fplus{{\,+_f\,}}

\newcommand{\epsclose}{\stackrel{\epsilon}{\thicksim} }
\newcommand{\eclose}[1]{\stackrel{{#1}}{\thicksim} }
\newcommand{\eps}{\epsilon}
\newcommand{\Anote}[1]{\begin{quote}{\sf Avi's Note:} {\sl{#1}} \end{quote}}

\begin{document}

\title{Incidence Theorems -- Lecture Notes}
\date{}
\maketitle


\fi

\section{The Elekes-Sharir framework}\label{sec-elekes}
 
In a recent breakthrough Guth and Katz \cite{GK10} proved that  any finite set of points $P$ in the real plane defines at least $\gtrsim |P|/\log |P|$ distinct distances. This is tight up to a $\sqrt{ \log |P|}$ factor and was conjectured by Erdos \cite{Erdos46}. The proof combines ideas that were developed in previous works and in particular a general framework developed by Elekes and Sharir in \cite{ES10} that gives a `generic' way to approach such problems by reducing them to an incidence problem. (The original paper of Elekes and Sharir reduced the distance counting problem to an incidence problem between higher degree curves but this was simplified by Guth and Katz to give lines instead of curves.)

To begin, observe that a 4-tuple of points $a,b,c,d \in P$ satisfies $\|a - b\| = \|c - d\|$ iff there exists a rigid motion  $T : \R^2 \mapsto \R^2$ (rotation $+$ translation) such that $T(a) = c$ and $T(b) = d$. Let us denote the set of rigid motions by $\cR$. To define a rigid motion we need to specify a translation (which has two independent parameters) and a rotation (one parameter) thus, we can think of $\cR$ as a three dimensional space. Later on we will fix a concrete parametrization of $\cR$ (minus some points) as $\R^3$ but for now it doesn't matter. For each $a,b \in P$ we define the set $$ L_{a,b} = \{ T \in \cR \,|\, T(a) = b \}$$ of rigid motions mapping $a$ to $b$. This set is `one dimensional' since, after specifying the image of $a$ we can only change the rotation parameter. In our concrete parametrization (which we will see shortly) all of the sets $L_{a,b}$ will in fact be lines in $\cR = \R^3$. Let $L = \{ L_{a,b}\,|\, a,b \in P \}$ be the set of $|P|^2$ lines defined by the point set $P$. We would like to bound the number of distances defined by $P$, denoted $d(P)$, as a function of the number of incidences between the lines in $L$. To this end, consider the following set: $$ Q(P) = \{ (a,b,c,d) \in P^4 \,|\, \|a-b\| = \|c - d\| \}.$$ A quick Cauchy-Schwarz calculation shows that $$ d(P) \geq \frac{|P|^4}{|Q(P)|}.$$ On the other hand, since each 4-tuple in $Q(P)$ gives an intersection between two lines in $L$, we have that $$ |Q(P)| \sim |I(L)| = |\{ (\ell,\ell') \in L^2 \,|\, \ell \cap \ell' \neq \emptyset \}|.$$ Thus, it will suffice to give a bound of $\lesssim |P|^3 \cdot \log |P|$ to obtain the bound of Guth-Katz on $d(P)$. In general, $|P|^2$ lines in $\R^3$
 can have many more intersections but, using some special properties of this specific family of lines we will in fact obtain this bound.

We now describe the concrete parametrization of $\cR$ mentioned above. We begin by removing from $\cR$ all translations. It is easy to see that the number of 4-tuples in $Q(P)$ arising from pure translations is at most $|P|^3$ (since every three points determine the fourth one). Now, every map in $\cR$ is a rotation by $\theta \in (0,2\pi)$  (say, to the right) around some fixed point $f = (f_x,f_y)$. If $T(a) = b$ then one sees that the fixed point $f$ must lie on the perpendicular bisector of segment $a,b$. That is, on the line $E_{a,b} = \{ c \,|\, \|c - a\| = \|c - b\| \}$ passing through the mid-point $m = (a-b)/2$ and in direction perpendicular to $a-b$. Our parametrization $\rho : \cR\setminus \{ \text{translations}\} \mapsto \R^3$ will be defined as $$\rho(T) = (f_x,f_y,(1/2)\cot(\theta/2)).$$ We now show that
\begin{claim}
For each $a,b \in P$ we have that $\rho(L_{a,b})$ is a line in $\R^3$.
\end{claim}
\begin{proof}
It will help to draw a picture at this point with  the two points $a,b$, the line $E_{a,b}$ and the fixed point $f$ on this line. We will consider the triangle formed by the three points $a,f$ and $m = (a-b)/2$ (the point on $E_{a,b}$ that is directly between $a$ and $b$). This is a right angled triangle with an angle of $\theta/2$ between the segments $fa$ and $fm$.  Thus, $(1/2)\cot(\theta/2) = \frac{\|f-m\|}{\|a-b\|}$. Let $d = (d_x,d_y)$ be a unit vector parallel to $E_{a,b}$. We thus have that $f = m + \|f-m\|\cdot d$ (or with a minus sign). Setting $t = \|f-m\|$, this shows that $$ \rho(L_{a,b}) = \left\{ (m_x,m_y,0) + t \cdot (d_x,d_y,\|a-b\|^{-1}) \,\,|\,\, t\in \R \right\}$$ which is a line.
\end{proof}

Let $N = |P|^2$. We have $N$ lines in $\R^3$ and want to show that there are at most $\sim N^{1.5} \log N$ incidences. This is clearly not true for an arbitrary set of $N$ lines. A trivial example where the number of incidences is $N^2$ is when all lines pass through a single point. Another example is when the lines are all in a single plane inside $\R^3$. If no two lines are parallel we would have $\sim N^2$ incidences (every pair intersects). Surprisingly enough, there is only one more type of examples with $\sim N^2$ incidences: doubly ruled surfaces. Take for example the set $S = \{ (x,y,xy) \,|\, x,y \in \R\} $. This set is `ruled' by two families of lines: the lines of the form $\{ (x,y,xy) \,|\, x \in \R\}$ and the lines of the form $\{ (x,y,xy) \,|\, y \in \R\}$. If we take $N/2$ lines from each family we will get $\sim N^2$ intersections. In other words, the set $S$ contains a `grid' of lines (horizontal and vertical) such that every horizontal lines intersects every vertical line. In general a doubly ruled surface is defined as a surface in which every point has two lines passing through it. A singly ruled surface is a surface in which every point has at least one line passing through it. It is known that the only non-planar doubly ruled surfaces (up to linear change of coordinates) is the one we just saw and the surface $\{ (x,y,z) \,|\, x^2 + y^2 - z^2 =1 \}$. There are no non-planar triply ruled surfaces. 

Guth and Katz proved the following:
\begin{thm}[Guth-Katz]\label{thm-GK1}
Let $L$ be a set of $N$ lines in $\R^3$ such that no more than $\sqrt N$ lines intersect at a single point and no plane or doubly ruled surface contains more than $\sqrt N$ lines. Then the number of incidences of lines in $L$, $|I(L)|$, is at most $\lesssim N^{1.5} \cdot \log N$.
\end{thm}
The bound in the theorem is tight (even with the logarithmic factor) and clearly the conditions cannot be relaxed. Luckily, the lines $L_{a,b}$ defined by a point set $P$ in the above manner satisfy the conditions of the theorem and so this theorem implies the bound on the number of distinct distances. An example of a set of lines matching the bound in the theorem is as follows: Take an $S \times S$ grid in the plane $z=0$ and another identical grid in the plane $z=1$ and pass a line through every two points, one on each grid. The number of lines is thus $N = S^4$ and a careful calculation shows that the number of incidences is $\sim S^6 \log S$ (see the appendix in Guth and Katz's original paper for the proof).

For each $k$ let $I_{\geq k}(L)$ denote the set of points that have at least $k$ lines in $L$ passing through them. Define $I_{=k}(L)$ in a similar manner requiring that there are exactly $k$ lines through the point. Theorem~\ref{thm-GK1} will follow from the following lemma (which is also tight using the same construction as above).
\begin{lem}\label{lem-GK1}
Let $L$ be as in Theorem~\ref{thm-GK1}.  Then for every $k \geq 2$, $$|I_{\geq k}(L)| \lesssim \frac{N^{1.5}}{k^2}.$$ 
\end{lem}

Let us see how this Lemma implies Theorem~\ref{thm-GK1}. 
\begin{eqnarray*}
|I(L)| &=& \sum_{k=2}^{\sqrt N} |I_{=k}(L)| \cdot k^2 \\
&=& \sum_{k} \left(|I_{\geq k}(L)| - |I_{\geq k+1}(L)| \right) \cdot k^2\\
&\lesssim& \sum_{k} |I_{\geq k}(L)|\cdot k \\
&\lesssim& N^{1.5} \cdot \sum_{k}(1/k) \lesssim N^{1.5} \cdot \log N.
\end{eqnarray*}

The case $k=2$ and $k \geq 3$ of the Lemma are proved in \cite{GK10} using  different arguments (the case $k=3$ was proved earlier in \cite{EKS11}). Interestingly, the case $k \geq 3$  does not require the condition on doubly ruled surfaces and can be proven without it.  

We  still need to show that the lines $L_{a,b}$ coming from $P$ satisfy the conditions of Theorem~\ref{thm-GK1} (we omit the mapping $\rho$ at this point to save on notations). To see that at most $\sqrt N = |P|$ lines meet at a point observe that, if not, we could find two lines $L_{a,b}$ and $L_{a,b'}$ that intersect. This clearly cannot happen since this would imply that there is a rigid motion $T$ mapping $a$ to $b$ and also mapping $a$ to $b'$. Let us consider now the maximum number of lines in a plane. Let $L_a = \{L_{a,b} \,|\, b \in P \}$ and observe that the lines in $L_a$ are disjoint. Moreover, by the parametrization of the lines, we have that all lines in $L_a$ are of different directions. Thus, every plane can contain at most one line from $L_a$. Thus, there could be at most $|P| = \sqrt N$ lines in a single plane. The condition regarding doubly ruled surfaces is more delicate and can be found in the Guth-Katz paper. 

In the next few sections we will develop the necessary machinery for proving Lemma~\ref{lem-GK1}. Since the proof is quite lengthy we will omit the proofs of some claims having to do with doubly ruled surfaces that are used in the proof of the $k=2$ case.  As a `warmup' to the full proof we will see two proofs of related theorems which use this machinery in a slightly simpler way. One of these will be yet another proof of the Szemeredi-Trotter theorem, this time using the polynomial ham sandwich theorem. The other will be the proof of the joints Conjecture which uses the polynomial method.

\ifediting
\bibliographystyle{alpha}
 \bibliography{incidence}

\end{document}
\fi

\newif\ifediting

\ifediting
\documentclass[11pt]{article}

\usepackage{amsmath,amsthm,amssymb}
\newcommand{\remove}[1]{}
\setlength{\topmargin}{0.3in} \setlength{\headheight}{0in}
\setlength{\headsep}{0in} \setlength{\textheight}{8.0in}
\setlength{\topsep}{0.1in} \setlength{\itemsep}{0.0in}
\parskip=0.05in
 \textwidth=6.5in 
\oddsidemargin=0truecm \evensidemargin=0truecm

\newtheorem{thm}{Theorem}[section]
\newtheorem{claim}[thm]{Claim}
\newtheorem{lem}[thm]{Lemma}
\newtheorem{define}[thm]{Definition}
\newtheorem{cor}[thm]{Corollary}
\newtheorem{obs}[thm]{Observation}
\newtheorem{example}[thm]{Example}
\newtheorem{construct}[thm]{Construction}
\newtheorem{conjecture}[thm]{Conjecture}
\newtheorem{THM}{Theorem}
\newtheorem{question}{Question}
\newtheorem{fact}[thm]{Fact}
\newtheorem{prop}[thm]{Proposition}

\def\F{{\mathbb{F}}}
\def\Q{{\mathbb{Q}}}
\def\Z{{\mathbb{Z}}}
\def\N{{\mathbb{N}}}
\def\R{{\mathbb{R}}}
\def\K{{\mathbb{K}}}
\def\C{{\mathbb{C}}}
\def\A{{\mathbb{A}}}
\def\P{{\mathbb{P}}}
\def\cP{{\cal P}}
\def\cS{{\mathcal S}}
\def\cE{{\mathcal E}}
\def\V{{\mathbf{V}}}
\def\I{{\mathbf{I}}}
\def\bx{{\mathbf x}}
\def\by{{\mathbf y}}
\def\E{{\mathbb E}}

\def\half{ \frac{1}{2}}
\newcommand{\ip}[2]{\langle #1,#2 \rangle}
\def\sumN{\sum_{i=1}^n}
\def\_{\,\,\,\,\,}
\def\prob{{\mathbf{Pr}}}
\newcommand{\entropy}[1]{ {\text{H}_{\infty}\left({#1}\right)} }
\def\modulo{\text{mod}}
\def\omm{ \{0,1\} }
\def\id{ \textit{id} }

\def\D{{\partial}}
\def\gap{\textsf{gap}}
\def\sign{\textsf{sign}}
\def\spar{\textsf{sparse}}
\def\span{\textsf{span}}
\def\Part{\textbf{Part}}
\def\Mon{\textbf{Mon}}
\def\sing{\textbf{sing}}
\def\Und{\textsf{Und}}
\def\Comp{\textsf{Comp}}
\def\rank{\textsf{rank}}
\def\poly{\textsf{poly}}
\def\codim{\textsf{codim}}
\def\cp{\textsf{cp}}
\def\uni{\textsf{Uni}}
\def\ext{\textsf{\bf Ext}}
\def\extt{\textsf{\bf Ext2}}

\def\fplus{{\,+_f\,}}

\newcommand{\epsclose}{\stackrel{\epsilon}{\thicksim} }
\newcommand{\eclose}[1]{\stackrel{{#1}}{\thicksim} }
\newcommand{\eps}{\epsilon}
\newcommand{\Anote}[1]{\begin{quote}{\sf Avi's Note:} {\sl{#1}} \end{quote}}

\begin{document}

\title{Incidence Theorems -- Lecture Notes}
\date{}
\maketitle


\fi

\section{The Polynomial Method and the joints Conjecture}\label{sec-polynomial}

The polynomial method is used to impose an algebraic structure on a geometric problem. The basic ingredient in this method is the following simple claim which holds over any field. Notice that the phrase `non zero polynomial' used in the claim refers to a polynomial with at last one non zero coefficient (over a finite field such a polynomial might still evaluate to zero everywhere).
\begin{claim}
Let $P \subset \F^n$ be a finite set, with $\F$ some field. If $|P| < {n + d \choose d}$ then there exists a non zero polynomial $g \in \F[x_1,\ldots,x_n]$ of degree $\leq d$ such that $g(p)=0$ for all $p \in P$. 
\end{claim}
\begin{proof}
Each constraint of the form $g(p)=0$ is a homogenous linear equation in the coefficients of $g$. The number of coefficients for a polynomial of degree at most $d$ in $n$ variables is exactly ${n +d \choose d}$ and so there must be a non zero solution.
\end{proof}
The second component of the polynomial method is given by bounding the maximum number of zeros a polynomial can have. In the univariate case, this is given by the following well-known fact. Later, we will see a variant of this claim for polynomials with more variables.

\begin{claim}
A non zero univariate polynomial $g(x)$ over a field $\F$ can have at most $\deg(g)$ zeros.
\end{claim}

To illustrate the power of the polynomial  method we will see how it gives a simple proof of the joints conjecture. To this end we must first prove some rather easy claims about restrictions of polynomials. Let $g \in \F[x_1,\ldots,x_n]$ be a degree $d$ polynomial over a field $\F$. Let $S \subset \F^n$ be any affine subspace of dimension  $k$. We can restrict $g$ to $S$ in the following way: write $S$ as the image of a degree one mapping $\phi : \F^k \mapsto \F^n$ so that $$ S =  \{ \phi(t_1,\ldots,t_k) \,|\, t_1,\ldots,t_k \in \F \}.$$ The restriction of $g$ to $S$ is the polynomial $h(t_1,\ldots,t_k) = g(\phi(t_1,\ldots,t_k))$. This depends in general on the particular choice of $\phi$ but for our purposes all $\phi$'s will be the same (we can pick any $\phi$). A basic and useful fact is that $\deg(h) \leq \deg(g)$ for any restriction $h$ of $g$. 

Suppose now that $\ell$ is a line in $\F^n$ and write $\ell$ as $\ell = \{ a + tb \,|\, t \in \F\}$ for some $a,b \in \F^n$. Restricting $g$ to $\ell$ we get a polynomial $h(t) = g(a+tb)$. It will be useful to prove some properties of this restriction. In particular, we would like to understand some of its coefficients. The constant coefficient is the value of $h$ at zero and is simply $h(0) = g(a)$. On the other hand, the coefficient of highest degree $t^d$ will be $g_d(b)$, where $g_d(x_1,\ldots,x_n)$ is the highest degree homogenous component of $g$ (that is, the sum of all monomial of maximal degree in $g$). Another coefficient we will try to understand is that of $t$. To this end we will use the partial derivatives of $g$. Recall that $\partial g / \partial x_i$ is a polynomial in $\F[x_1,\ldots,x_n]$ obtained by taking the derivative of $g$ as a polynomial in $x_i$ (with coefficients in $\F[x_j, j \neq i]$). This is defined over any field but notice some weird things can happen if $\F$ has positive characteristic (e.g, the partial derivative of $x^p$ is zero over $\F_p$ even though this is a non constant polynomial). The gradient of $g$ is the vector of polynomials $$ \nabla(g) = (\partial g / \partial x_1, \ldots, \partial g / \partial x_n ). $$ This vector has geometric meaning which we will not discuss here. Algebraically, we have that the coefficient of $t$ in the restriction $h(t) = g(a + tb)$ to a line is exactly $\ip{\nabla(g)(a)}{b}$. To see this, observe that the coefficient of $t$ is obtained by taking the derivative (w.r.t the single variable $t$) and then evaluating the derivative at $t=0$. Using the chain rule for functions of several variables we get that the derivative of $h(t)$ is $\ip{ \nabla(g)(a+tb)}{b}$ and so the claim follows.

\subsection{The joints problem}

 Let $L$ be a set of lines in $\R^3$. A `joint' w.r.t the arrangement $L$ is a point $p \in \R^3$ through which pass at least three, non coplanar, lines. The basic question  one can  ask is `what is the maximal number of joints possible in an arrangement of $N$ lines'. This problem, raised in \cite{Chazelle90} in relation to computer graphics algorithms, was answered completely by Guth and Katz \cite{GK10b} using a clever application of the polynomial method. This result followed a long line of papers proving incremental results using various techniques, quite different than the polynomial method (see \cite{GK10b} for a list of references). This proof of the joints conjecture by Guth and Katz was the first case where the polynomial method was used directly to argue about problems in Euclidean space (in contrast to finite fields where it was more common). Later in \cite{GK10}, ideas from this proof were used in part of the proof of the Erdos-Distance problem bound. 

A simple lower bound on the number of joints is obtained from taking  $L$ to be the union of the following three sets, each containing $N^2$ lines:
\[ L_{xy} = \{ (i,j,t) | t \in \R\},  \,\, i,j \in [N] \]
\[ L_{yz} = \{ (t,i,j) | t \in \R\},  \,\, i,j \in [N] \]
\[ L_{zx} = \{ (i,t,j) | t \in \R\},  \,\, i,j \in [N] \]
In other words, the set $L$ contains $\sim N^2$ lines in a three dimensional grid. It is easy to check that every point in $[N]^3$  is a joint and so we have that the number of joints can be as large as $|L|^{3/2}$. Not surprisingly (at this point), this is tight.

\begin{thm}[Guth Katz \cite{GK10b}]
Let $L$ be a set of lines in $\R^3$. Then $L$ defines at most $\lesssim |L|^{3/2}$ joints.
\end{thm}
\begin{proof}
The proof given here is a simplified proof found by Kaplan, Sharir and Shustin  \cite{KSS10} who also generalized it to higher dimensions.

Let $J$ be the set of joints and suppose that $|J| > C|L|^{3/2}$ for some large constant $C$ to be chosen later. We can throw away all lines in $L$ that have fewer than $|J|/2|L|$ joints on them. This can decrease the size of $J$ by at most a half.

Let $g(x,y,z)$ be a non zero polynomial with real coefficients of minimal degree that vanishes on the set $J$. We saw in previous sections that $\deg(g) \lesssim |J|^{1/3}$ (since a polynomial of this degree will have $\gtrsim |J|$ coefficients). 
 
Since each of the lines in $L$ contains $|J|/2|L| \geq \deg(g)$ joints (here we take the constant $C$ to be large enough) we get that $g$ must vanish identically on each line in $L$. To see this, observe that the restriction of $g(x,y,z)$ to a line is also a degree $\leq \deg(g)$ polynomial and so, if it is not identically zero, it can have at most $\deg(g) $ zeros. Thus, we have moved from knowing that $g$ vanishes on all joints to knowing that $g$ vanishes on all lines!

Consider a joint $p \in \R^3$ and let $\ell_1,\ell_2,\ell_3 \in L$ be three non coplanar lines passing through $p$. We can find three linearly independent vectors $u_1,u_2,u_3 \in \R^3$ such that for all $i \in [3]$ we have $\ell_i = \{ p + tu_i | t \in \R\}$. Since $g$ vanishes on these three lines we get that for $i \in [3]$, the polynomial $ h_i(t) = g(p+tu_i)$ is identically zero. This means that all the coefficients of $h_i(t)$ are zero and in particular the coefficient of $t$ which is, as we saw, equal to  $\ip{u_i}{\nabla(g)(p)}$. Since the $u_i$'s are linearly independent, we get that $\nabla(g)(p) = 0$ for all joints $p \in J$. One can now check that, over the reals, a non zero polynomial $g$ has at least one non zero partial derivative of degree strictly less that the degree of $g$. Therefore, we get a contradiction since we assumed that $g$ is a minimal degree polynomial vanishing on $J$. 
\end{proof}

It is not hard to generalize this proof to finite fields using the fact that a polynomial $g \in \F[x_1,\ldots,x_n]$ all of whose partial derivatives are zero must be of the form $f(x)^p$ for some other polynomial $f$, where $p$ is equal to the characteristic of the field. For other generalizations, including to algebraic curves, see \cite{KSS10}.

\ifediting
\bibliographystyle{alpha}
 \bibliography{incidence}

\end{document}
\fi

\newif\ifediting

\ifediting
\documentclass[11pt]{article}

\usepackage{amsmath,amsthm,amssymb}
\newcommand{\remove}[1]{}
\setlength{\topmargin}{0.3in} \setlength{\headheight}{0in}
\setlength{\headsep}{0in} \setlength{\textheight}{8.0in}
\setlength{\topsep}{0.1in} \setlength{\itemsep}{0.0in}
\parskip=0.05in
 \textwidth=6.5in 
\oddsidemargin=0truecm \evensidemargin=0truecm

\newtheorem{thm}{Theorem}[section]
\newtheorem{claim}[thm]{Claim}
\newtheorem{lem}[thm]{Lemma}
\newtheorem{define}[thm]{Definition}
\newtheorem{cor}[thm]{Corollary}
\newtheorem{obs}[thm]{Observation}
\newtheorem{example}[thm]{Example}
\newtheorem{construct}[thm]{Construction}
\newtheorem{conjecture}[thm]{Conjecture}
\newtheorem{THM}{Theorem}
\newtheorem{question}{Question}
\newtheorem{fact}[thm]{Fact}
\newtheorem{prop}[thm]{Proposition}

\def\F{{\mathbb{F}}}
\def\Q{{\mathbb{Q}}}
\def\Z{{\mathbb{Z}}}
\def\N{{\mathbb{N}}}
\def\R{{\mathbb{R}}}
\def\K{{\mathbb{K}}}
\def\C{{\mathbb{C}}}
\def\A{{\mathbb{A}}}
\def\P{{\mathbb{P}}}
\def\cP{{\cal P}}
\def\cS{{\mathcal S}}
\def\cE{{\mathcal E}}
\def\V{{\mathbf{V}}}
\def\I{{\mathbf{I}}}
\def\bx{{\mathbf x}}
\def\by{{\mathbf y}}
\def\E{{\mathbb E}}

\def\half{ \frac{1}{2}}
\newcommand{\ip}[2]{\langle #1,#2 \rangle}
\def\sumN{\sum_{i=1}^n}
\def\_{\,\,\,\,\,}
\def\prob{{\mathbf{Pr}}}
\newcommand{\entropy}[1]{ {\text{H}_{\infty}\left({#1}\right)} }
\def\modulo{\text{mod}}
\def\omm{ \{0,1\} }
\def\id{ \textit{id} }

\def\D{{\partial}}
\def\gap{\textsf{gap}}
\def\sign{\textsf{sign}}
\def\spar{\textsf{sparse}}
\def\span{\textsf{span}}
\def\Part{\textbf{Part}}
\def\Mon{\textbf{Mon}}
\def\sing{\textbf{sing}}
\def\Und{\textsf{Und}}
\def\Comp{\textsf{Comp}}
\def\rank{\textsf{rank}}
\def\poly{\textsf{poly}}
\def\codim{\textsf{codim}}
\def\cp{\textsf{cp}}
\def\uni{\textsf{Uni}}
\def\ext{\textsf{\bf Ext}}
\def\extt{\textsf{\bf Ext2}}

\def\fplus{{\,+_f\,}}

\newcommand{\epsclose}{\stackrel{\epsilon}{\thicksim} }
\newcommand{\eclose}[1]{\stackrel{{#1}}{\thicksim} }
\newcommand{\eps}{\epsilon}
\newcommand{\Anote}[1]{\begin{quote}{\sf Avi's Note:} {\sl{#1}} \end{quote}}

\begin{document}

\title{Incidence Theorems -- Lecture Notes}
\date{}
\maketitle


\fi

\section{The Polynomial ham sandwich theorem}\label{sec-polyhs}

One of the main ingredients in the proof of the Guth-Katz theorem is an ingenious use of the polynomial ham sandwich theorem, originally proved by Stone and Tukey  \cite{ST42}. This is a completely different use of polynomials than the one we saw for the joints problem and combines both algebra and topology.  We will  demonstrate its usefulness by seeing how it can be used to give yet another  proof of the Szemeredi-Trotter theorem in two dimensions. This proof will be both `intuitive' (not `magical' like the crossing number inequality proof) and simple (without the messy technicalities of the cell partition proof we saw). 

The folklore ham sandwich theorem states that every three bounded open sets in $\R^3$ can be simultaneously bisected using a single plane. If we identify the three sets with two slices of bread and a slice of ham we see the practical significance of this theorem. More generally, we have:
\begin{thm}[\cite{ST42}]
Let $U_1,\ldots,U_n \subset \R^n$ be bounded open sets. Then there exists a hyperplane  $H = \{ x \in \R^n \,|\, h(x) = 0\}$ (with $h(x)$ a degree one polynomial in $n$ variables) such that for each $i \in [n]$ the two sets $U_i \cap H^+= \{x \in U_i \,|\, h(x) > 0 \}$ and $U_i \cap H^- = \{x \in U_i \,|\,  h(x) < 0 \}$  have equal volume. In this case we say that $H$ {\em bisects} each of the $U_i$'s.
\end{thm}
   
This can be thought of as extending the basic fact that for every $n$ points there is a hyperplane in $\R^n$ that passes through all of them. The proof uses the Borsuk-Ulam theorem from topology, whose proof can be found in \cite{MatBU}.
\begin{thm}[Borsuk-Ulam \cite{Borsuk}]
Let $S^n \subset \R^{n+1}$ be the $n$-dimensional unit sphere and let $f : S^n \mapsto \R^n$ be a continuous map such that $f(-x) = -f(x)$ for all $x \in S^n$ (such a map is called antipodal). Then there exists $x$ such that $f(x) = 0$.
\end{thm}

\begin{proof}[Proof of the ham-sandwich theorem:]
Each hyperplane corresponds to some degree one polynomial $h(x) = h_0 + h_1x_1 + \ldots + h_nx_n $. Since we only care about the sign of $h$ we can scale so that the coefficients form a unit vector $v_h = (h_0,h_1,\ldots,h_n)$. We define a function $f : S^n \mapsto \R^n$ as follows $$f(v_h) = \left( |H^+ \cap U_i| - |H^- \cap U_i|\right)_{i \in [n]}.$$ It is clear that $f$ is continuous and that $f(-x) = -f(x)$. Thus, there exists a zero of $f$ and we are done.
\end{proof}

In their original paper, Stone and Tukey also observed that if we want to bisect more sets, we can do it as long as we have enough degrees of freedom. One way to allow for more degrees of freedom is to replace a hyperplane with a {\em hypersurface}. 

A hypersurface is a set $H = \{ x \in \R^n \,|\, h(x) = 0\}$ where now $h$ can be a polynomial of arbitrary degree $d$. The degree of $H$ is defined to be the degree of its defining polynomial (we will abuse this definition a bit and say that a hypersurface has degree $d$ if it has degree \textbf{bounded} by $d$). Recall that, if we have $t < {n+d \choose d}$ points in $\R^n$ than we can find, by interpolation, a non zero degree $d$ polynomial that is zero on all of them. For the problem of bisecting open sets the same holds: If the number of sets is smaller than ${n+d \choose d}$  we can find a degree $d$ polynomial that bisects all of the sets.

\begin{thm}[Polynomial ham sandwich (PHS)]
Let $U_1,\ldots,U_t \in \R^n$ be bounded open sets with $t < {n+d \choose d}$. Then there exists a degree $d$ hypersurface $H$ that bisects each of the sets $U_i, i \in [t]$.
\end{thm}
\begin{proof}
The proof is identical to the degree one proof. Identify each degree $d$ hypersurface with its (unit) vector of coefficients and apply the Borsuk-Ulam theorem on the function $f$ mapping to the differences. 
\end{proof}

\subsection{Cell partition using polynomials}

The PHS theorem gives a particularly nice way to partition $\R^n$ into cells. In addition to having a `balanced' partition (as we had in the cell partition method we saw earlier) we will have the additional useful property that the boundaries of the partition are defined using a low degree polynomial. The use of the PHS for cell partition originated in a paper of Guth \cite{Guth08} on the multilinear Kakeya problem.

The first step for obtaining the cell partition theorem is to take the PHS to the `limit' and replace the open sets with discrete sets. If $S \subset \R^n$ is a finite set and $H$ is a hypersurface, we say that $H$ bisects $S$ if both sets $S \cap H^-$ and $S \cap H^+$ have size at most $|S|/2$. Notice that this definition  allows for an arbitrary number of points from $S$ to belong to the set $H$ itself.
 
\begin{lem}[Discrete PHS]\label{lem-dphs}
Let $S_1,\ldots,S_t \subset \R^n$ be $t$ finite sets of points with $t < {n+d \choose d}$. Then, there exists a degree $d$ hypersurface $H$ that bisects each of the sets $S_i, i \in [t]$.
\end{lem}
\begin{proof}
Consider $\eps$-neighborhoods $U_1,\ldots,U_t$ of the sets $S_1,\ldots,S_t$ and apply the PHS on this family of open sets obtaining a bisecting hypersurface $H_\eps$. Taking $\eps$ to zero and using the compactness of the unit sphere we get that there is  sequence of bisecting hypersurfaces converging to some degree $d$ hypersurface $h$. If one of the sets $S_i \cap H^+$ or $S_i \cap H^-$ has size larger than $|S_i|/2$ we could find a h.s $h_\eps$ that does not bisect the $\eps$-neighborhood of $S_i$.
\end{proof}

Notice that, if $n$, the dimension, is fixed and the number of sets $t$ grows, we always have a degree $d = O_n( t^{1/n})$ polynomial that bisects $t$ sets. In particular, over $\R^2$, a family of $t$ discrete sets can be bisected by a degree $\sim \sqrt t$ h.s.

We will now use the discrete PHS to get our final cell partition theorem.  We will only need this theorem over $\R^2$ and $\R^3$ but will state it over $\R^n$ for all $n$ (it will help to think of $n$ as a fixed constant and of $t$ as growing to infinity). 

\begin{thm}[Polynomial Cell Partition]
Let $S \subset \R^n$ be a finite set and let $t \geq 1$. Then, there exists a decomposition of $\R^n$ into $O(t)$ cells (open sets) such that each cell has boundary in a hypersurface $H$  of  degree $d = O_n( t^{1/n})$ and each cell contains at most $|S|/t$ points from $|S|$. Notice that the cells do not have to be connected.
\end{thm}
\begin{proof}
We will apply the discrete PHS iteratively to obtain finer and finer partitions. Initially, we get a h.s $H_1$ of degree $d_1 \leq O_n( 1^{1/n})$ that bisects the single set $|S|$ into two sets of size at most $|S|/2$
 each (plus some points on the boundary). Applying the discrete PHS again on these two sets we obtain a degree $d_2 = O_n( 2^{1/n})$ h.s $H_2$ that bisects both sets. This gives four cells (with boundary in the h.s $H_1 \cup H_2$ wich has degree at most $d_1 + d_2$ since its defined by the product of polynomials defining each h.s) with at most $|S|/4$ points in each. Continuing in this fashion $\ell = \log_2 t$ times we obtain $\ell$ hypersurfaces $H_1,\ldots,H_\ell$ with $H_j$ having degree  $O_n(2^{j/n})$ and such that their union $H = \cup_{j \in [\ell]} H_j$, gives a partition into cells containing at most $|S|/t$ points each. The degree of $H$ is bounded by the sum $$ \sum_{j=1}^{\ell} O_n( 2^{j/n}) = O_n( t^{1/n}). $$	
\end{proof}

\subsection{Szemeredi-Trotter using polynomials} 
Using the Polynomial Cell Partition theorem, we either get a `balanced' partition of most of the points into disjoint cells or there is a large number of points that lies on a low degree hypersurface (hence, they possess an algebraic structure). Kaplan, Matousek and Sharir \cite{KMS11} used this argument to give another proof of the Szemeredi-Trotter theorem which we will now see. 
 
Before we can start the proof we need to prove a very simple algebraic claim that we will use in the `algebraic' case (when many points are on the h.s). In general, the polynomial method always requires some additional algebraic claims that depend on the specific problem (e.g., for the joints problem we needed to look at the coefficients of the restriction to lines and express them using the gradient). In this case we can prove what we need in a few lines. In other cases, we will rely on more powerful theorems from algebraic geometry.

\begin{claim}\label{cla-lines}
Let $H \subset \R^2$ be a h.s of degree $d$. Then 
\begin{enumerate}
	\item For every line $\ell \subset \R^2$ not contained in $H$ we have $|\ell \cap H| \leq d$.
	\item There are at most $d$ lines contained in $H$.
\end{enumerate}
\end{claim}
\begin{proof}
Let $h(x,y)$ be a polynomial of degree $\leq d$ defining $H$. Let $\ell$ be a line and consider the restriction of $h$ to $\ell$. As we already discussed, this is a univariate polynomial of degree at most $\deg(h)$ and so, if it is not identically zero, it can have at most $\deg(h)$ zeros. 

To prove 2, suppose that there were $d+1$ lines $\ell_1,\ldots,\ell_{d+1}$ distinct lines contained in $H$. A generic line\footnote{We use the term `generic line' to mean `any line outside some set of measure zero'. More accurately, if we parametrize lines as vectors of coefficients defining them, a generic line is any line outside some fixed set of zeros of some polynomial. Over the reals one can simply take a `perturbed' line. We can also use the word generic for other objects such as hypersurfaces, sets of points etc. with the same meaning.}  $\ell$ will (a) not be contained in $H$ and (b) will intersect each of the $d+1$ lines. This will contradict 1 and so 2. is proved as well.
\end{proof}

Later we will see a more general statement of this form known as  Bezout's theorem.


We can now give the proof of the Szemeredi-Trotter theorem using polynomial cell partition. Let $P,L$ be our sets of points/lines in $\R^2$. As in previous proofs we will use the Cauchy Schwarz bound(s):
$$ |I(P,L)| \lesssim |P||L|^{1/2} + |L|$$ and
$$ |I(P,L)| \lesssim |L||P|^{1/2} + |P|.$$ We can also assume that $|P|^{1/2} \ll |L| \ll |P|^2$ (otherwise the theorem follows from the CS bound).

We will apply the polynomial cell partition theorem with parameter $t$ to be chosen later to obtain a hypersurface $H$ of degree $d = O(\sqrt t)$ and a family of disjoint cells $C_1,\ldots,C_t$ such that $\R^2$ is the disjoint union of the cells and of $H$. For each $i \in [t]$ let $P_i$ denote the set $P \cap C_i$ and let $L_i$ denote the set of lines in $L$ that intersect the cell $C_i$. We also define $P_0$ to be the set of points in $P \cap H$ and $L_0$ to be the set of lines in $L$ that intersect $H$. 

We thus have:
$$ |I(P,L)| = |I(P_0,L_0)| +  \sum_{i=1}^t |I(P_i,L_i)|.$$	
We will use different arguments to bound each of the terms on the r.h.s. The sum of incidences for $i > 0$ can be bounded as follows. We have for each $i \in [t]$, $|P_i| \leq |P|/t$. So,  applying the CS bound on each cell we obtain
\begin{equation}\label{eq-CS}
|I(P_i,L_i)| \lesssim (|P|/t)\cdot |L_i|^{1/2} + |L_i|.
\end{equation}
Since each line in $L_i$ is not contained in $H$ we have, by Claim~\ref{cla-lines}, that it can intersect at most $d = O(\sqrt t)$ cells (since it must cross $H$ when it moves from one cell to another). This gives a bound
$$ \sum_{i=1}^t |L_i| \lesssim t^{1/2}|L|, $$ which, using Cauchy-Schwarz, gives
$$ \sum_{i=1}^t |L_i|^{1/2} \lesssim t^{3/4}|L|^{1/2}.$$
Combining the above we get
\begin{equation}\label{eq-sumt}
\sum_{i=1}^t|I(P_i,L_i)| \lesssim (|P|/t)\cdot t^{3/4}|L|^{1/2} + t^{1/2}|L| = t^{-1/4}|P||L|^{1/2} + t^{1/2}|L|.
\end{equation}

To bound $|I(P_0,L_0)|$ we first split $L_0$ into two sets: $L_0'$ containing lines that are contained in $H$ and $L_0''$ containing lines that are not contained in $H$ (but intersect it at some point). By Claim~\ref{cla-lines} we have $|I(P_0,L_0'')| \leq t^{1/2}|L|$ since each line in $L_0''$ can intersect $H$ in at most $d \sim \sqrt t$ points. We also have $|L_0'| \leq d \sim t^{1/2}$ and so, using the CS bound we have 
$$ |I(P_0,L_0')| \lesssim t^{1/2}|P|^{1/2} + |P| \lesssim t^{1/2}|L| + |P|.$$

Combining the above  we get
$$ |I(P,L)| \lesssim t^{-1/4}|P||L|^{1/2} + t^{1/2}|L| + |P|. $$ Setting $$t \sim \frac{|P|^{4/3}}{|L|^{2/3}}$$ gives the Szemeredi-Trotter theorem. We have $t \geq 1$ since $|L| \ll |P|^2$. This completes the proof \qed.

\ifediting
\bibliographystyle{alpha}
 \bibliography{incidence}

\end{document}
\fi

\newif\ifediting

\ifediting
\documentclass[11pt]{article}

\usepackage{amsmath,amsthm,amssymb}
\newcommand{\remove}[1]{}
\setlength{\topmargin}{0.3in} \setlength{\headheight}{0in}
\setlength{\headsep}{0in} \setlength{\textheight}{8.0in}
\setlength{\topsep}{0.1in} \setlength{\itemsep}{0.0in}
\parskip=0.05in
 \textwidth=6.5in 
\oddsidemargin=0truecm \evensidemargin=0truecm

\newtheorem{thm}{Theorem}[section]
\newtheorem{claim}[thm]{Claim}
\newtheorem{lem}[thm]{Lemma}
\newtheorem{define}[thm]{Definition}
\newtheorem{cor}[thm]{Corollary}
\newtheorem{obs}[thm]{Observation}
\newtheorem{example}[thm]{Example}
\newtheorem{construct}[thm]{Construction}
\newtheorem{conjecture}[thm]{Conjecture}
\newtheorem{THM}{Theorem}
\newtheorem{question}{Question}
\newtheorem{fact}[thm]{Fact}
\newtheorem{prop}[thm]{Proposition}

\def\F{{\mathbb{F}}}
\def\Q{{\mathbb{Q}}}
\def\Z{{\mathbb{Z}}}
\def\N{{\mathbb{N}}}
\def\R{{\mathbb{R}}}
\def\K{{\mathbb{K}}}
\def\C{{\mathbb{C}}}
\def\A{{\mathbb{A}}}
\def\P{{\mathbb{P}}}
\def\cP{{\cal P}}
\def\cS{{\mathcal S}}
\def\cE{{\mathcal E}}
\def\V{{\mathbf{V}}}
\def\I{{\mathbf{I}}}
\def\bx{{\mathbf x}}
\def\by{{\mathbf y}}
\def\E{{\mathbb E}}

\def\half{ \frac{1}{2}}
\newcommand{\ip}[2]{\langle #1,#2 \rangle}
\def\sumN{\sum_{i=1}^n}
\def\_{\,\,\,\,\,}
\def\prob{{\mathbf{Pr}}}
\newcommand{\entropy}[1]{ {\text{H}_{\infty}\left({#1}\right)} }
\def\modulo{\text{mod}}
\def\omm{ \{0,1\} }
\def\id{ \textit{id} }

\def\D{{\partial}}
\def\gap{\textsf{gap}}
\def\sign{\textsf{sign}}
\def\spar{\textsf{sparse}}
\def\span{\textsf{span}}
\def\Part{\textbf{Part}}
\def\Mon{\textbf{Mon}}
\def\sing{\textbf{sing}}
\def\Und{\textsf{Und}}
\def\Comp{\textsf{Comp}}
\def\rank{\textsf{rank}}
\def\poly{\textsf{poly}}
\def\codim{\textsf{codim}}
\def\cp{\textsf{cp}}
\def\uni{\textsf{Uni}}
\def\ext{\textsf{\bf Ext}}
\def\extt{\textsf{\bf Ext2}}

\def\fplus{{\,+_f\,}}

\newcommand{\epsclose}{\stackrel{\epsilon}{\thicksim} }
\newcommand{\eclose}[1]{\stackrel{{#1}}{\thicksim} }
\newcommand{\eps}{\epsilon}
\newcommand{\Anote}[1]{\begin{quote}{\sf Avi's Note:} {\sl{#1}} \end{quote}}

\begin{document}

\title{Incidence Theorems -- Lecture Notes}
\date{}
\maketitle


\fi

\section{The Guth-Katz incidence theorem for lines in $\R^3$}\label{sec-gkthm}
We have now developed enough machinery and intuition to start discussing the proof of the Guth-Katz theorem regarding incidences of lines in $\R^3$. Recall that the statement we are trying to prove is:
\begin{thm}[\cite{GK10}]\label{thm-GK}
Let $L$ be a set of $N^2$ lines in $\R^3$ such that no more than $ N$ lines intersect at a single point and no plane or doubly ruled surface contains more than $ N$ lines. Then the number of incidences of lines in $L$, $|I(L)|$, is at most $\lesssim N^3 \cdot \log N$. 
\end{thm}
Also recall that we argued that this Theorem will follow from the following estimate on the sets $I_{\geq k}(L)$ of points that have at least $k$ lines in $L$ passing through them:
\begin{lem}\label{lem-GK}
Let $L$ be as in Theorem~\ref{thm-GK}.  Then for every $k \geq 2$, $$|I_{\geq k}(L)| \lesssim \frac{N^{3}}{k^2}.$$ 
\end{lem}
We will prove Lemma~\ref{lem-GK} first for $k \geq 3$ and then for $k=2$ (using different arguments). Since the statement is asymptotic we can actually separate into the two cases when $k$ is either larger than some big constant $C$ or smaller than $C$ (the case $k < C$ will only use the fact that at least two lines meet at a point).

\subsection{The $k \geq 3$ case}
This case of the lemma does not require any conditions on doubly ruled surfaces and so we only assume that no plane contains more than $N$ lines in $L$. We can also assume that $k < N$ since each point has at most $N$ lines through it. 

The high level idea is as follows: Using the polynomial cell partition theorem, we can partition the points in $I_{\geq k}(L)$ into cells whose boundary is a low degree surface. We then separate into two cases: the cellular case and the algebraic case. The cellular case is when most of the points are inside the interior of the cells. In this case we will use the `weak' three dimensional Szemeredi-Trotter theorem (meaning, the ST theorem one gets from projecting everything to a generic plane) in each cell and sum up the resulting bounds. This case is very similar to the cell partition proof of the ST theorem we saw. The second, and harder, case is when most points are on the algebraic surface. The proof in this case is similar  to the proof of the joints conjecture with the added difficulty that some intersections are planar. In the algebraic case we will argue using a degree argument that the surface must contain many of the lines in $L$ (those lines with many points on them). We will then use the assumption $k \geq 3$ to argue that these lines are `special' in some concrete sense and that a surface that contains too many `special' lines must contain a plane and this plane must contain many of the lines (contradicting our assumption). In this last part of the proof we will also need to distinguish between points that have 3 non-coplanar lines through them and points through which there are 3 planar lines.

Since the full proof requires some careful book-keeping we will make some simplifying assumption along the way. These will usually be benign and can be removed easily by simple averaging arguments. To begin,  we assume the following two `regularity' assumptions: 
\begin{enumerate}
	\item Every point in $I_{\geq k}(L)$ has at most $2k$ lines in $L$ passing through it. (To remove this assumption we need to argue about each interval $[2^i,2^{i+1}]$ and sum the results).
	\item Each line in $L$ is incident to at least $\gtrsim \frac{|I_{\geq k}(L)|k}{|L|}$ lines. This is the `average' number of lines incident to a point and so many points will have at least some fraction of this number of lines passing through them. 
\end{enumerate}

Let $S = I_{\geq k}(L)$ be the set whose size we wish to bound. Suppose $|S| \geq C \cdot \frac{N^3}{k^2}$ for some large constant $C$ to be specified later. We will use the cell partition lemma obtained from the polynomial ham sandwich theorem (stated here for  $\R^3$):
\begin{thm}[Polynomial Cell Partition]
Let $S \subset \R^3$ be a finite set and let $t \geq 1$. Then, there exists a decomposition of $\R^3$ into $\lesssim t$ cells (open sets) such that each cell has boundary in a hypersurface $H$  of  degree $d \lesssim t^{1/3}$ and s.t each cell contains at most $|S|/t$ points.
\end{thm}
We will apply this theorem and choose the parameter $t$ so that the hypersurface $H$ will be of degree $d = \lceil 3\cdot (N/k) \rceil$. We can assume $k \ll N$ since otherwise the bound we are trying to prove is trivial. This guarantees that $d \geq 1$. This means that each cell contains at most $ \lesssim |S|/d^3$ points and each line passes through at most $d$ cells (since crossing between cells means intersecting $H$ and a line not contained in $H$ can intersect it in at most $\deg(H)$ points as we already saw). Let $S_H = S \cap H$ and let $S_C = S \setminus S_H$. Clearly, one of these sets will have size at least $|S|/2$. We begin with the case $|S_C| \geq |S|/2$.

\subsubsection{The cellular case}
Assume $|S_C| \geq |S|/2$. We will use the following easy corollary of the Szemeredi-Trotter theorem. We already used this corollary in dual form (with lines replaced with points) when we proved Beck's theorem. Even though we proved this bound in the plane $\R^2$, the same statement holds in three dimensions using a generic projection to a plane (which preserves intersections). 
\begin{cor}\label{cor-rich2}
Let $P$ and $L$ be sets of points and lines in $\R^3$. For $k > 1$ let $P_k$ denote the set of points in $P$ that have at least $k$ lines passing through them. Then,
$$|P_k| \lesssim \frac{|L|^2}{k^3} + \frac{|L|}{k}. $$
\end{cor}
Let $L_i$ denote the set of lines in $L$ that pass through the $i$'th cell. Applying Corollary~\ref{cor-rich2} to each cell we get
\begin{equation}\label{eq-sumi}
\frac{|S|}{2} \leq |S_C| \leq \sum_i \left( \frac{|L_i|^2}{k^3} + \frac{|L_i|}{k}\right)
\end{equation}
Observe that $\sum_i |L_i| \leq d\cdot |L|$ since each lines passes through at most $d$ cells. Also, from our first regularity assumption (each point has at most $2k$ lines passing through it), we get that $\max_i |L_i| \leq \frac{|S|2k}{d^3}$. This implies that $$\sum_i |L_i|^2 \leq \max_i |L_i| \cdot \sum_i |L_i| \leq \frac{|S|\cdot |L| \cdot 2k}{d^2} = \frac{|S|\cdot N^2 \cdot 2k}{d^2}. $$ Using these bounds in (\ref{eq-sumi}) we obtain
\begin{eqnarray*}
\frac{|S|}{2} &\leq&  \frac{|S|\cdot N^2 \cdot 2}{k^2\cdot d^2} + \frac{d \cdot N^2}{k}\\
&=& \frac{2|S|}{9} + 3\cdot \frac{N^3}{k^2} \ll \frac{|S|}{2},
\end{eqnarray*}
where the last inequality used the assumption that $|S| \gg N^3/k^2$. This is a contradiction and so we conclude that we must be in the algebraic case.

\subsubsection{The algebraic case}
 Observing the proof of the cellular case we see that we could also get a contradiction if $|S_C| \geq |S|/100$ or any other small constant. This only requires taking $d = D \cdot (N/k)$ for a larger constant $D$ (instead of $D=3$). This means that we can also get $|S_H| \geq (1-\eps)|S|$ for any constant $\eps$. Taking $\eps$ small enough and doing some averaging arguments (removing points not on $H$) we can actually reduce to the case where all points in $S$ are in $S_H$ so from now on we make this further simplifying assumption.

Thus, the situation is as follows: We have a set of points $S$  with $|S| \gg N^3/k^2$ such that all points lie on a hypersurface $H$ of degree $d  \lesssim  N/k$ and such that through every point in $S$ there are $k \geq 3$ lines in some set $L$ of $N^2$ lines. The assumption $k \geq 3$ will come into play now when we analyze algebraic properties of $H$ at the points in $S$. 

Choosing $C$ large enough (so that $|S| > CN^3/k^2$) and using our second regularity assumption (that each lines has many points on it) we get that each line in $L$ contains at least $$\gtrsim \frac{|S|\cdot k}{|L|} \geq 10d$$ points in $S \subset H$ on it (the constant $10$ will be important later). This means, in particular, that all the lines in $L$ are completely contained in $H$. Thus, each point in $S$ has three lines in $H$ passing through it. We will separate these points into `critical' points, through which there are three non-coplanar lines (as in the joints problem), and to `flat' points, which are non critical and through which there are three planar lines. We can also define `critical' lines to be those lines that contain at least $5d$ critical points and `flat' lines that contain at least $5d$ flat points. Since each line has at least $10d$ points on it we have that each line is either critical or flat. We separate again into two cases depending on whether at least half the lines are critical or at least half are flat lines. 

\subsubsection*{At least $N^2/2$ critical lines}

Recall the proof of the joints conjecture. We saw that if a surface $H = \{ h(x,y,z) = 0 \}$ has three non-coplanar lines passing through a point $p \in H$ then the gradient $$\nabla_h = (\partial h / \partial x, \partial h / \partial y, \partial h / \partial z ),$$ which is composed of three polynomials of degree $\leq \deg(h)$, must vanish at $p$. Let $f_1,f_2,f_3 \in \R[x,y,z]$ denote the three components of the gradient (so that $f_1 = \partial h / \partial x$ etc..). Since they must vanish on all critical points, and since each critical line contains $> 5d$ critical points, we have that $f_1,f_2,f_3$ must also vanish on all critical lines. Thus, the hypersurface $H$ shares many lines with each of the three hypersurfaces $F_i = \{f_i(x,y,z) = 0 \}$. We would like to say that this cannot happen. For an arbitrary pair of hypersurfaces $G_1 = \{ g_1(x,y,z)=0\}$ and $G_2 = \{ g_2(x,y,z)=0\}$ there can be no bound on the number of lines contained in both since the two polynomials $g_1$ and $g_2$ might share a factor $r(x,y,z)$ (that is, $r$ divides both) such that the hypersurface $R = \{ r(x,y,z)=0 \}$ contains many lines (perhaps an infinite number of lines). Fortunately, this is the only case where this can happen. We will now prove this fact using the following classical result known as Bezout's theorem.
\begin{thm}[Bezout]
Let $f(x,y), g(x,y) \in \R[x,y]$ be two polynomial without a common factor. Then, $f$ and $g$ have at most $\deg(f)\deg(g)$ common roots.
\end{thm}
The proof of this result is not hard but requires some discussion of resultants which is beyond the scope of this survey. Using Bezout's theorem we can prove the following claim.
\begin{claim}\label{cla-sharelines}
Let $G_1 = \{ g_1(x,y,z)=0\}$ and $G_2 = \{ g_2(x,y,z)=0\}$ be two hypersurfaces so that $g_1$ and $g_2$ do not have a common factor. Then $G_1 \cap G_2$ can contain at most $\deg(g_1)\cdot \deg(g_2)$ lines.
\end{claim}
\begin{proof}
Take a generic plane $A \subset \R^3$ and consider the restrictions $\hat g_1(u,v)$ and $\hat g_2(u,v)$ of $g_1,g_2$ to this plane ($u,v$ are new variables parametrizing the plane). It is not hard to show that, since $g_1,g_2$ do not share a factor, the restrictions to a generic plane will also not have a common factor. This means that $\hat g_1,\hat g_2$ can have at most $\deg(g_1)\deg(g_2)$ common roots. But a generic plane will intersect each of the lines contained in $G_1 \cap G_2$ in distinct points and so each line will add a common zero to $\hat g_1,\hat g_2$. This shows that the number of lines is bounded by $\deg(g_1)\deg(g_2)$.
\end{proof}

To use Claim~\ref{cla-sharelines} we need to argue that $h$ (the polynomial defining $H$) does not have  a factor in common with  its partial derivatives $f_1,f_2,f_3$. This, however, can  happen if $h$ has a repeated factor. Namely, if we factor $h$ into its irreducible components $h = \prod_{j} p_j(x,y,z)^{\alpha_j}$ then one of the $\alpha_j$'s is at least $2$. Such a polynomial will share a factor  $p_j$ with each of its partial derivatives. In our case, since we are only interested in the set of zeros of $h$ (and do not mind reducing the degree) we can assume w.l.o.g that $h$ has no repeated factors (i.e., is square-free). For square-free polynomials, one can easily show that $h$ does not share a factor with at least one of the partial derivatives. This means that, using Claim~\ref{cla-sharelines}, there are at most $ d^2 \ll N^2/2$ lines -- a contradiction. This brings us to the last case:

\subsubsection*{There are at least $N^2/2$ flat lines}
This case is similar but will require us to use the assumption that no plane contains more than $N$ lines. We saw that critical points are common zeros of some family of three low degree polynomials and that this family of polynomials cannot have a common factor with $h$. For flat points a similar statement holds but with $9$ polynomials. Specifically:
\begin{claim}\label{cla-flat}
There are $9$ polynomials $\pi_1,\ldots,\pi_9 \in \R[x,y,z]$ of degree at most $3d$ each such that:
\begin{enumerate}
	\item Each flat point is a zero of all $9$ polynomials $\pi_i$.
	\item If $h(x,y,z)$ is `plane-free' (i.e., if no irreducible factor of $h$ is of degree one) then $h$ does not share a factor with at least one of the polynomials $\pi_j$.
\end{enumerate}
\end{claim}
We will not prove this claim here and just say that these $9$ mysterious polynomials are the `second fundamental form' of the surface $H$ and are related to the second order derivatives of $h$ (or more precisely, to the quadratic approximation of $H$ at a point). Since the degrees of the $\pi_j$'s are bounded by $3d$ and since there are at least $5d$ flat points on each flat line we get that all flat lines are contained in the $9$ hypersurfaces $\Pi_j = \{ \pi_j(x,y,z) = 0 \}$. 

We now write $h(x,y,z) = h_p(x,y,z) \cdot h_n(x,y,z)$, where $h_p$ contains all the `planar' (degree one)  irreducible components of $h$ and $h_n$ is `plane-free'. Thus, the hypersurface $H_p = \{ h_p(x,y,z) = 0 \}$ is the union of all planes contained in $H$. Using Claim~\ref{cla-flat} and Claim~\ref{cla-sharelines} we get that $H$ and $H_n$ can share at most $5d^2 \ll N^2$ lines and so there are many ($\gtrsim N^2$)  lines contained in $H_p$ (clearly, each line in $H$ must be contained in one of the irreducible components). By a pigeonhole argument, and using the fact hat $H_p$ can have at most $\deg(H_n) \leq d$ degree one components, we get that at least $\gtrsim N^2/d > N$ lines are in some plane. This is a contradiction and so the proof of Lemma~\ref{lem-GK} is complete for the case $k \geq 3$.

\subsection{The $k=2$ case}

We now have to prove Lemma~\ref{lem-GK} in the case $k=2$. Here we will eventually use the fact that at most $N$ lines are in a doubly-ruled surface. Recall that a doubly ruled surface is a surface in which every point has two lines passing through it. A singly ruled surface is a surface in which every point has at least one line passing through it. We are not assuming anything on the number of lines contained in a singly ruled surface. It is known that there are only two examples (up to isomorphism) of doubly ruled surfaces, both of degree two. In order for us to not stray too far off our topic, we will state some facts about doubly (and singly) ruled surfaces without proof. The interested reader can find the missing proofs (or pointers to them) in \cite{GK10}.

Let us begin with the proof and denote again by $S = I_{\geq 2}(L)$ the set of points of intersection of at least two lines in $L$. We will want to prove that $|S| \lesssim N^3$ so, for contradiction, suppose $|S| > C\cdot N^3$ for some large constant $C$ to be chosen later. 

The proof uses again the polynomial method. This time, unlike the $k \geq 3$ case, we will jump straight to the `algebraic case' and find a polynomial that vanishes on all lines in $L$. We already saw that a degree $\lesssim t^{1/3}$ polynomial can be found that vanishes on a given set of $t$ points in $\R^3$. We wish to prove a similar statement with lines instead of points. The next claim does just that (we state it only for $\R^3$ but a similar claim holds for any dimension and any field).
\begin{claim}[Simple Interpolation]\label{cla-interpol1}
Let $\ell_1,\ldots,\ell_t$ be $t$ lines in $\R^3$. Then there exists a non zero polynomial of degree $\leq 10\cdot t^{1/2}$ that vanishes on all of the lines $\ell_i$ (i.e., the restriction of the polynomial to each of the lines is identically zero).
\end{claim}
\begin{proof}
A polynomial $f(x,y,z)$ of degree $10\cdot t^{1/2}$ has ${10t^{1/2} + 3 \choose 3} > 10 t^{1.5}$ coefficients. Each constraint of the form $f|_{\ell_i} \equiv 0$ ($f$ vanishes on $\ell_i$) gives at most $\deg(f)+1 \leq 10t^{1/2}+1$ homogenous linear equations in the coefficients of $f$ (each coming from the vanishing of one of the coefficients of the univariate restriction to the line $\ell_i$). Thus, we have enough coefficients to satisfy all the constraints in a non-trivial way.
\end{proof}
Notice that this proof is very general and can work in many other settings (with lines replaced with almost any algebraic object you can think of).

This claim is indeed simple but not very useful in our case. Applying it directly on the set $L$ will give us a polynomial of degree $d \leq 10N$ that vanishes on all lines in $L$. For reasons that will become clear later, we will actually need a polynomial of much smaller degree (a small constant times $N$) to vanish on $L$. 
Fortunately, $L$ is not an arbitrary set of lines (which would make this task impossible) but a set with many intersections. Since $|S| > CN^3$ we know that a constant fraction of the  lines have at least $CN/10$ points of intersection on them. Replacing $C$ with $C/10$ and throwing away a small fraction of the lines we can assume w.l.o.g that each line in $L$ has at least $C\cdot N$  distinct points of intersection on it, where $C$ is some large constant to be chosen later. Using this additional structure we can prove the following improved version of the interpolation claim.
\begin{claim}[Interpolation using incidences]\label{cla-interpol2}
Suppose $C$ is a large enough constant. Let $L$ be a set of $N^2$ lines in $\R^3$ such that each line in $L$ intersects at least $CN$ other lines in distinct points. Then, there exists a non zero polynomial of degree $d \lesssim N/\sqrt{C}$ that vanishes on all lines in $L$.
\end{claim}
\begin{proof}
Take a random subset $L'$ of $L$ by picking each line with probability $1/C$. With high probability, each line in $L$ (no typo, this is the original set $L$) will still intersect at least $N/2$  lines in $L'$. Use Claim~\ref{cla-interpol1} to find a polynomial $f(x,y,z)$ of degree $10\sqrt{|L'|} \lesssim N/\sqrt{C}$ that vanishes on $L'$. We now observe that $f$ must vanish also on $L$ since the restriction of $f$ to each line in $L$ has at least $N/2 > \deg(f)$ zeros (when $C$ is large enough).
\end{proof}

Let $f(x,y,z)$ be a polynomial of degree $d \leq N/\sqrt{C}$ given by the last claim such that $f$ vanishes on all lines in $L$. Let $F$ be the hypersurface defined by $f$. Write $f = \prod_{i}f_i(x,y,z)$ as a product of irreducible polynomials $f_i$ and recall that, w.l.o.g, $f$ is square free and so no $f_i$ is repeated. Thus, $F$ is the union of the hypersurfaces $F_i$ defined by the different $f_i$'s. If we denote by $d_i$ the degree of $f_i$ we have that $d = \sum_i d_i$. We now partition the irreducible factors into 4 groups. Let $f_{pl}$ be the product of all $f_i$'s that are of degree one (corresponding to $F_i$'s that are planes). Let $f_{dr}$ be the product of the doubly-ruled irreducible components. Let $f_{sr}$ denote the product of the singly-ruled components and Let $f_{nr}$ be the product of the non-ruled components. Also define $F_{pl},F_{dr},F_{sr}$ and $F_{nr}$ to be the hypersurfaces defined by these four polynomials. 

Since each line $\ell \in L$ is contained in $F$, it must be contained in one of the irreducible factors of $F$. The set of incidences $S$ can be partitioned into incidences between lines in different factors and to incidences of lines that are in the same factor. A line $\ell$ in a factor $f_i$ can intersect lines in factors not containing $\ell$ in at most $d$ points. This follows from the basic fact (that we proved in previous sections) that a line can have at most $\deg(H)$ intersections with a hypersurface $H$ not containing it. Here we use this fact for $H$ being the hypersurface defined by the product of the factors of $f$ that do not contain the line $\ell$ (which contains all lines that do not share a factor with $\ell$). This means that the total number of incidences between lines in different factors is bounded by $|L|d \lesssim N^3$.

We thus have to consider only incidences between lines that are in the same factor. Consider first intersections of lines in the factors of $F_{pl}$. Since there are at most $N$ lines in each plane, we have at most $N^2$ intersections inside each plane. Since there are at most $d \leq N$ planar factors in $F_{pl}$, the total number of incidences of this kind is bounded by $\lesssim N^3$. The same argument precisely works also for intersections between lines in $F_{dr}$ using the assumption that every doubly ruled surface contains at most $N$ lines. 

We now consider intersections of lines in  $F_{sr}$. We will use the following fact from the theory of ruled surfaces:
\begin{claim}
Let $S \subset \R^3$ be a singly-ruled surface. Then, every line in $S$, with the exception of at most two lines, can intersect at most $\deg(S)$ other lines in $S$.
\end{claim}
In other words, if there are 3 lines in $S$ that have more than $\deg(S)$ intersection with other lines in $S$, than $S$ must be doubly ruled. Using this claim, we can bound the number of incidences in $F_{sr}$ by $\lesssim N^3$ as follows. Each singly ruled factor $F_i \subset F_{sr}$ can have  two `exceptional' lines which can have at most $|L| = N^2$ incidences each. This sums up to $|\{\text{components of } F_{sr} \}|\cdot 2N^2 \leq dN^2 \leq N^3$. Each `non exceptional' line in a factor $F_i$ of $F_{sr}$ can contribute at most $\deg(F_i) = d_i \leq d$ intersections and so the total is again bounded by $|L|d \leq N^3$. 

We are left with the task of bounding the number of intersections of lines in $F_{nr}$. Suppose that there are more than $AN^3$ such intersections, where $A$ is some large constant to be chosen later. We will use the following claim without proof:
\begin{claim}\label{cla-nonruled}
A non-ruled surface $S \subset \R^3$ can contain at most $\deg(S)^2$ lines.
\end{claim}
This means that $F_{nr}$ can contain at most $d^2  \leq N^2/C$ lines. Notice that we can pick $A$ (which controls the number of incidences) and $C$ to be as large as we want and so we can  argue by induction on  the problem of bounding the number of incidences of lines in $F_{nr}$. That is, we can assume Lemma~\ref{lem-GK} (for $k=2$) holds for $(N-1)^2$ lines and then use this assumption on the lines in $F_{nr}$. This requires some careful choice of constants but can be carried out (we will not do this here). A delicate point is that we must satisfy the assumption that at most $N-1$ lines in $F_{nr}$ are in any plane or doubly ruled surface. This can be achieved using the following argument: As long as there is a plane or doubly ruled surface containing more than $N-1$ lines in $L' = L \cap F_{nr}$, remove these lines from $L'$. This can repeat at most $\lesssim N$ times and so the intersections between lines we have removed are at most $\lesssim N^3$. The remaining lines are not contained in any of the planes of doubly ruled surfaces we removed and so have at most $\lesssim N^3$ intersections with removed lines (using the same argument we used before). This will result in a small decrease in the constants that can be ignored since we can choose $C$ and $A$ to be sufficiently large constants.

\ifediting
\bibliographystyle{alpha}
 \bibliography{incidence}

\end{document}
\fi

\newif\ifediting

\ifediting
\documentclass[11pt]{article}

\usepackage{amsmath,amsthm,amssymb}
\newcommand{\remove}[1]{}
\setlength{\topmargin}{0.3in} \setlength{\headheight}{0in}
\setlength{\headsep}{0in} \setlength{\textheight}{8.0in}
\setlength{\topsep}{0.1in} \setlength{\itemsep}{0.0in}
\parskip=0.05in
 \textwidth=6.5in 
\oddsidemargin=0truecm \evensidemargin=0truecm

\newtheorem{thm}{Theorem}[section]
\newtheorem{claim}[thm]{Claim}
\newtheorem{lem}[thm]{Lemma}
\newtheorem{define}[thm]{Definition}
\newtheorem{cor}[thm]{Corollary}
\newtheorem{obs}[thm]{Observation}
\newtheorem{example}[thm]{Example}
\newtheorem{construct}[thm]{Construction}
\newtheorem{conjecture}[thm]{Conjecture}
\newtheorem{THM}{Theorem}
\newtheorem{question}{Question}
\newtheorem{fact}[thm]{Fact}
\newtheorem{prop}[thm]{Proposition}

\def\F{{\mathbb{F}}}
\def\Q{{\mathbb{Q}}}
\def\Z{{\mathbb{Z}}}
\def\N{{\mathbb{N}}}
\def\R{{\mathbb{R}}}
\def\K{{\mathbb{K}}}
\def\C{{\mathbb{C}}}
\def\A{{\mathbb{A}}}
\def\P{{\mathbb{P}}}
\def\cP{{\cal P}}
\def\cS{{\mathcal S}}
\def\cE{{\mathcal E}}
\def\V{{\mathbf{V}}}
\def\I{{\mathbf{I}}}
\def\bx{{\mathbf x}}
\def\by{{\mathbf y}}
\def\E{{\mathbb E}}

\def\half{ \frac{1}{2}}
\newcommand{\ip}[2]{\langle #1,#2 \rangle}
\def\sumN{\sum_{i=1}^n}
\def\_{\,\,\,\,\,}
\def\prob{{\mathbf{Pr}}}
\newcommand{\entropy}[1]{ {\text{H}_{\infty}\left({#1}\right)} }
\def\modulo{\text{mod}}
\def\omm{ \{0,1\} }
\def\id{ \textit{id} }

\def\D{{\partial}}
\def\gap{\textsf{gap}}
\def\sign{\textsf{sign}}
\def\spar{\textsf{sparse}}
\def\span{\textsf{span}}
\def\Part{\textbf{Part}}
\def\Mon{\textbf{Mon}}
\def\sing{\textbf{sing}}
\def\Und{\textsf{Und}}
\def\Comp{\textsf{Comp}}
\def\rank{\textsf{rank}}
\def\poly{\textsf{poly}}
\def\codim{\textsf{codim}}
\def\cp{\textsf{cp}}
\def\uni{\textsf{Uni}}
\def\ext{\textsf{\bf Ext}}
\def\extt{\textsf{\bf Ext2}}

\def\fplus{{\,+_f\,}}

\newcommand{\epsclose}{\stackrel{\epsilon}{\thicksim} }
\newcommand{\eclose}[1]{\stackrel{{#1}}{\thicksim} }
\newcommand{\eps}{\epsilon}
\newcommand{\Anote}[1]{\begin{quote}{\sf Avi's Note:} {\sl{#1}} \end{quote}}

\begin{document}

\title{Incidence Theorems -- Lecture Notes}
\date{}
\maketitle


\fi

\section{Application of the Guth-Katz  bound to sum-product estimates}\label{sec-sumprodgk}

We saw how the Szemeredi Trotter theorem can be used to derive the sum product theorem over the reals showing that for every set $A \subset \R$ one of the sets $A+A,A\cdot A$ has size at least $|A|^{5/4}$. We will now see how the three dimensional incidence theorem proved by Guth and Katz can be used to give a similar estimate. Recall the GK bound:
\begin{thm}[Guth-Katz]\label{thm-GK3}
Let $L$ be a set of $N^2$ lines in $\R^3$ such that no more than $ N$ lines intersect at a single point and no plane or doubly ruled surface contains more than $ N$ lines. Then the number of incidences of lines in $L$, $|I(L)|$, is at most $\lesssim N^3 \cdot \log N$.
\end{thm}

Recently, Iosevich, Roche-Newton and Rudnev \cite{IRR11} used this theorem to prove:
\begin{thm}\label{thm-sumprod2}
Let $A \subset \R$ be any set. Then $$ |A\cdot A - A\cdot A| \gtrsim \frac{|A|^2}{\log |A|}. $$
\end{thm}
(The same result also holds with the minus sign replaced by a plus). Notice that, unlike the previous sum-product estimates we saw, this bound is tight up to the logarithmic factor. This follows by taking $A$ to be an arithmetic progression.

To get a feeling why such a bound should follow from the result of Guth-Katz observe that the GK bound for the distinct distances problem (obtained from Theorem~\ref{thm-GK3} in a black-box way) automatically gives that, for a set $A \subset \R$ we have $$|\{ (a-b)^2 + (c-d)^2 \,|\, a,b,c,d \in \A \}| \gtrsim \frac{|A|^2}{\log |A|}.$$ This is obtained by taking the set of points $P = A \times A$ in the plane and counting the distances defined by this set. This bound is also tight (up to logarithmic factors) for an arithmetic progression. To argue about $A\cdot A - A\cdot A$, however, we will need to use a slightly different reduction. Recall that the reduction from the distance problem to the incidence bound was obtained by considering the group of distance preserving linear mappings acting on the plane. To prove Theorem~\ref{thm-sumprod2} we will need to consider mappings that preserve determinants (or areas). This is the group $\SL_2(\R)$ of $2 \times 2$ matrices with determinant equal to 1. Again, this is a three dimensional group and the trick will be to identify it with $\R^3$ in a way that the mappings sending a point $p$ to a point $q$ form  a line in $\R^3$ (this will actually be simpler to show in this case).

For two vectors $v = (a,b)$ and $u = (c,d)$ we denote $\det(v,u) = ad-bc$. Observe that for four vectors $v,u,v',u'$, no two of which are multiples of each other, we have that $\det(v,u)=\det(v',u')$ iff there exists a map $T \in \SL_2(\R)$ that sends $v$ to $v'$ and $u$ to $u'$. One direction is obvious. To see the other direction let $v,u,v',u'$ be as above and observe that there is a unique $T$ sending $v$ to $v'$ and $u$ to $u'$. We now have that $\det(v',u') = \det(T)\cdot \det(v,u)$ and so $\det(T)=1$ as required.

Fix two vectors $v = (a,b),v' = (c,d) \in \R^2$ that are not multiples of each other (i.e., $\det(v,v') \neq 0$). Let $L_{v,v'} \subset \SL_2(\R)$ be the set of mappings with $L(v)=v'$. We wish to understand how this set looks like. If $v=(1,0)$ and $v' = (0,1)$ this is easy: 
\[ L_{(1,0),(0,1)} = \left\{ \left(\begin{array}{cc}
0 & -1 \\
1 & t \end{array} \right)\,\,,\,\, t \in \R \right\}. \]
For general $v = (a,b),v' = (c,d)$ we need to conjugate by the matrix taking $v,v'$ to $(1,0), (0,1)$. This gives the set:
\[ \frac{1}{ad-bc}\cdot \left(\begin{array}{cc}
a & c \\
b & d \end{array} \right) \cdot
\left(\begin{array}{cc}
0 & -1 \\
1 & t \end{array} \right) \cdot
\left(\begin{array}{cc}
d & -c \\
-b & a \end{array} \right)\,\,, t \in \R
\]
Which gives the line
\[ \frac{1}{ad-bc}\cdot
\left(\begin{array}{cc}
cd +ab - bct & -c^2 - a^2 + act \\
d^2 + b^2 - bdt & -cd -ab +adt \end{array} \right)\,\,, t \in \R.
\]
The lines $L_{v,v'}$ are contained in the three dimensional surface $H = \{(x_1,x_2,x_3,x_4)|  x_1x_4 - x_2x_3 = 1\}$ which lives in $\R^4$. We can project this surface to $\R^3$ using the projection $(x_1,x_2,x_3)$. This projection is one-to-one as long as $x_1 \neq 0$. Using a generic rotation around zero, we can assume that this projection preserves the structure of the finite set of lines we will be interested in (i.e., those coming from $v,v'$ with both coordinates in $A$). 

Let $P = A \times A$ and let $L = \{L_{v,v'}\, | \, v,v' \in P\}$ be our set of $|P|^2$ lines. Following the  Elekes-Sharir framework we define the set
$$ Q(P) = \{ (v,u,v',u') \in P^4 \,|\, \det(v,u) = \det(v',u') \}.$$ Applying Cauchy-Schwarz we get that $$ |A\cdot A - A\cdot A| = |\{ \det(v,u)\, |\, v,u \in P\}| \geq \frac{|P|^4}{|Q(P)|}.$$ On the other hand, since each 4-tuple in $Q(P)$ gives an intersection between two lines in $L$, we have that $$ |Q(P)| \sim |I(L)| = |\{ (\ell,\ell') \in L^2 \,|\, \ell \cap \ell' \neq \emptyset \}|.$$ Thus, it will suffice to give a bound of $\lesssim |P|^3 \cdot \log |P|$. This bound will follow from Theorem~\ref{thm-GK3} if we can argue that the set $L$ satisfies the conditions of the theorem. As before, the condition on at most $|P|$ lines through a single point follows from the fact that no mapping can map a single point to two distinct points. The two conditions on planes and doubly ruled surfaces can be verified from the explicit description of the lines in $L$ given above.

To prove the same statement for $A\cdot A + A \cdot A$ observe that the size of $Q(P)$ is the same in this case (since $ad-bc=a'd'-b'c'$ iff $ad+b'c'=bc+a'd'$).

\ifediting
\bibliographystyle{alpha}
 \bibliography{incidence}

\end{document}
\fi

\chapter{Counting Incidences Over Finite Fields}\label{chap-2}


\newif\ifediting

\ifediting
\documentclass[11pt]{article}

\usepackage{amsmath,amsthm,amssymb}
\newcommand{\remove}[1]{}
\setlength{\topmargin}{0.3in} \setlength{\headheight}{0in}
\setlength{\headsep}{0in} \setlength{\textheight}{8.0in}
\setlength{\topsep}{0.1in} \setlength{\itemsep}{0.0in}
\parskip=0.05in
 \textwidth=6.5in 
\oddsidemargin=0truecm \evensidemargin=0truecm

\newtheorem{thm}{Theorem}[section]
\newtheorem{claim}[thm]{Claim}
\newtheorem{lem}[thm]{Lemma}
\newtheorem{define}[thm]{Definition}
\newtheorem{cor}[thm]{Corollary}
\newtheorem{obs}[thm]{Observation}
\newtheorem{example}[thm]{Example}
\newtheorem{construct}[thm]{Construction}
\newtheorem{conjecture}[thm]{Conjecture}
\newtheorem{THM}{Theorem}
\newtheorem{question}{Question}
\newtheorem{fact}[thm]{Fact}
\newtheorem{prop}[thm]{Proposition}

\def\F{{\mathbb{F}}}
\def\Q{{\mathbb{Q}}}
\def\Z{{\mathbb{Z}}}
\def\N{{\mathbb{N}}}
\def\R{{\mathbb{R}}}
\def\K{{\mathbb{K}}}
\def\C{{\mathbb{C}}}
\def\A{{\mathbb{A}}}
\def\P{{\mathbb{P}}}
\def\cP{{\cal P}}
\def\cS{{\mathcal S}}
\def\cE{{\mathcal E}}
\def\V{{\mathbf{V}}}
\def\I{{\mathbf{I}}}
\def\bx{{\mathbf x}}
\def\by{{\mathbf y}}
\def\E{{\mathbb E}}

\def\half{ \frac{1}{2}}
\newcommand{\ip}[2]{\langle #1,#2 \rangle}
\def\sumN{\sum_{i=1}^n}
\def\_{\,\,\,\,\,}
\def\prob{{\mathbf{Pr}}}
\newcommand{\entropy}[1]{ {\text{H}_{\infty}\left({#1}\right)} }
\def\modulo{\text{mod}}
\def\omm{ \{0,1\} }
\def\id{ \textit{id} }

\def\D{{\partial}}
\def\gap{\textsf{gap}}
\def\sign{\textsf{sign}}
\def\spar{\textsf{sparse}}
\def\span{\textsf{span}}
\def\Part{\textbf{Part}}
\def\Mon{\textbf{Mon}}
\def\sing{\textbf{sing}}
\def\Und{\textsf{Und}}
\def\Comp{\textsf{Comp}}
\def\rank{\textsf{rank}}
\def\poly{\textsf{poly}}
\def\codim{\textsf{codim}}
\def\cp{\textsf{cp}}
\def\uni{\textsf{Uni}}
\def\ext{\textsf{\bf Ext}}
\def\extt{\textsf{\bf Ext2}}

\def\fplus{{\,+_f\,}}

\newcommand{\epsclose}{\stackrel{\epsilon}{\thicksim} }
\newcommand{\eclose}[1]{\stackrel{{#1}}{\thicksim} }
\newcommand{\eps}{\epsilon}
\newcommand{\Anote}[1]{\begin{quote}{\sf Avi's Note:} {\sl{#1}} \end{quote}}

\begin{document}

\title{Incidence Theorems -- Lecture Notes}
\date{}
\maketitle


\fi

\section{Ruzsa calculus}\label{sec-ruzsa}

We begin  developing the necessary machinery for proving a Szemeredi-Trotter type result over prime finite fields. This result, due to Bourgain, Katz and Tao \cite{BKT04}, says that a set of $N$ lines and $N$ points in $\F_p^2$ can have at most $O(N^{1.5 - \eps})$ incidences, where $\eps$ is some positive real constant and $N$ is not too large. We will discuss the precise statement of this theorem in Section~\ref{sec-stfinite} after we have developed some machinery from additive combinatorics in this and the following two sections. 

The first ingredient we will need is Ruzsa calculus \cite{Ruzsa96}. This is a set of small claims about additive structure in arbitrary abelian groups. The usefulness of this calculus will become clear in the following sections.

 Let $G$ be an abelian group and let $A,B \subset G$ be subsets. We already defined the sets $A+B$ and $A-B$ of sums/differences of elements of $A$. We can define $k\cdot A$ to be the set $A+A+A+\ldots +A$, $k$-times. Be careful not to confuse this with the set $\{ka \,|\, a \in A\}$ which always has size bounded by $|A|$. We will generally only work with finite subsets of $G$  and so will omit the word `finite' in all of our definitions/claims. We use the Cartesian product notation $A \times B = \{ (a,b) \,|\, a \in A, b \in B \}$. 
 
We begin with a simple, yet powerful, lemma known as Ruzsa's triangle inequality.
\begin{lem}[Triangle Inequality]\label{lem-triangle}
Let $A,B,C \subset G$. Then
$$ |A|\cdot |B - C| \leq |A-B|\cdot |A-C|. $$
\end{lem}
\begin{proof}
We can define two functions $f : B-C \mapsto B$ and $g : B-C \mapsto C$ such that For every $v \in B-C$  it holds that $f(v) - g(v) = v$ (these are not defined uniquely, just pick some pair of values with difference $v$). Consider the mapping $$\phi : A \times (B-C) \mapsto (A-B)\times (A-C)$$ given by $$ \phi(a,v) = (a - f(v),a - g(v)).$$ Observe that $\phi$ is injective since we can recover $v$ from the difference between the two coordinates of the output. This implies the required bound on the set sizes.
\end{proof}

To explain the name of this lemma, consider the `Ruzsa Distance' between two sets: $$ d(A,B) = \log \left( \frac{|A-B|}{\sqrt{|A||B|}}\right).$$ The lemma just proved shows that $d(B,C) \leq d(A,B) + d(A,C)$ which justifies calling this function `distance' (though it is not a real distance function since $d(A,A)$ might be non zero).

\begin{thm}[Ruzsa calculus]\label{thm-Ruzsacalc}
There exists an absolute constant $c$ such that the following holds. Let $A,B,C \subset G$ be such that $|A|=|B|=|C|=N$.
\begin{enumerate}
	\item If $|A+B| \leq K\cdot N$ then $|A-B| \leq K^c \cdot N$.
	\item If $|A+B| \leq K\cdot N$ then $|A+A| \leq K^c \cdot N$.
	\item If $|A+B|,|A+C| \leq K\cdot N$ then $|B+C| \leq K^c \cdot N$.
	\item If $|A+B| \leq K\cdot N$, $|C+C| \leq K \cdot N$ and $|A \cap C| \geq K^{-1} \cdot N$ then $|C+B| \leq K^c \cdot N$.
	\item If $|A+B|,|A+C| \leq K\cdot N$ then $|A+B+C| \leq K^c \cdot N$.
	\item If $|A+A| \leq K\cdot N$ then for all non-negative integers $k,\ell$ there exists $c(k,\ell)$ such that $|k \cdot A- \ell \cdot A| \leq K^{c(k,\ell)} \cdot N$, where $c(k,\ell)$ does not depend on the group $G$ or on the set $A$. 
\end{enumerate}
\end{thm}
\begin{proof}
Throughout the proof we will (ab)use the constant $c$ freely and treat it as a `generic' constant that can change from one line to another (all inequalities will eventually work if we pick $c$ large enough). A cleaner way to do this would be to define $A \lesssim B$ to mean $ A \leq K^c B$ for some absolute constant $c$.

We start with some useful notations (some of which will be familiar). Let $$Q(A,B) = \{ (a,b,a',b') \in A \times B \times A \times B \,|\, a+b = a'+b'\}$$ and let $S(x) = \{ (a,b) \in A \times B \,|\, a+b = x \}$ and 	$R(x) = \{ (a,b) \in A \times B \,|\, a-b = x \}$. Then $$ |Q(A,B)| = \sum_{x}|S(x)|^2 = \sum_{x}|R(x)|^2.$$ Recall that, using Cauchy-Schwarz, we get that $$ |Q(A,B)| = \sum_{x}|S(x)|^2 \geq \frac{ \left( \sum_{x}|S(x)|\right)^2}{|A+B|} = \frac{|A|^2 |B|^2}{|A+B|}.$$ Using the fact that $Q(A,B) = Q(A,-B)$ we also get that $$ |Q(A,B)| \geq  \frac{|A|^2 |B|^2}{|A-B|}.$$

\paragraph{1.} From $$ |A||B| \cdot \max_{x}|R(x)| \geq \sum_{x}|R(x)|^2 \geq \frac{|A|^2 |B|^2}{|A+B|} $$ we get $$ \max_{x}|R(x)| \geq \frac{|A| |B|}{|A+B|}. $$ Let $x_0$ be such that $|R(x_0)| \geq \frac{|A| |B|}{|A+B|} $. We define the map $$ \phi : R(x_0) \times (A-B) \mapsto (A+B) \times (A+B)$$ to be $$ \phi((a,b),v) = (f(v)+b, g(v) + a),$$ where $f,g$ are fixed functions on $A-B$ such that $f(v) \in A, g(v) \in B$ and $f(v)-g(v)=v$. Again, we can check that $\phi$ is injective, which gives $$ \frac{|A| |B|}{|A+B|} \leq |R(x_0)| \leq \frac{|A+B|^2}{|A-B|}.$$ Plugging in the bound $|A+B| \leq KN$ we obtain the bound $|A-B| \leq K^3N$.

\paragraph{2.} Using the triangle inequality (Lemma~\ref{lem-triangle}) and 1. we get $$ |B||A-A| \leq |B-A||B-A| \leq K^c \cdot N^2. $$

\paragraph{3.} Similarly, 
\begin{eqnarray*}
(N/K^c)|B+C| \leq |A||B-C| \leq \\ |A-B||A-C| \leq K^c |A+B||A+C| \leq K^c N^2.
\end{eqnarray*}

\paragraph{4.} Using the triangle inequality and the previously proved parts
\begin{eqnarray*}
(N/K)|B-C| \leq |A\cap C||B-C| \leq \\ |A\cap C-B||A\cap C-C| \leq |A-B||C-C| \leq K^c N^2.
\end{eqnarray*}

\paragraph{5.} Here we need to do some work. The main step is to find a set $S$ of size roughly $N$ such that $|S+(A+B)| \leq K^c N$. Then we will have a bound on  $|(A+B)+(A+B)|$, which will become a bound on $|A+B+C|$ using 4. and the fact that some shift of $C$ has large intersection with $A$ (and so with $A+B$). It is a good exercise to stop reading now and try to fill in the rest of the proof.

We will take $S$ to be $$S = \{ x \in G\, |\, |S(x)| \geq N/10K \}.$$ Since $\sum_{x} |S(x)|^2 \geq N^4/|A+B| \geq N^3/K$ we must have $|S| \geq N/10K$. Observe that each element $x+(a+b) \in S + (A+B)$ has at least $N/10K$ distinct representations as a sum of the form $$x+a+b = (a_i + a) + (b_i + b)$$ with $a,a_i \in A,\,\,b,b_i\in B$ s.t $a_i + b_i = x$. This means that $$ |S+ (A+B)| \leq (10K/N)|A+A||B+B| \leq K^c N.$$ We now use part 2. to obtain  $$|(A+B) + (A+B)| \leq K^c N.$$ 

Since $|A-C| \leq K^c N$ there exists an element $x \in G$ that can be written in at least $N/K^c$ ways as a difference $x = a-c$ with $a\in A, c \in C$. This means that $|(C+x) \cap A| \geq N/K^c$. This implies $|(C+x) \cap (A+B)| \geq N/K^c$. We also know that $|(C+x)+(C+x)| \leq |C+C| \leq K^c N$ (since, from 3, $|B+C| \leq K^c \cdot N$) and so, using 4, we get $$ | (C+x) + (A+B) | \leq K^c N.$$ which gives $|A+B+C| \leq K^c N$.

\paragraph{6.} This follows immediately from a repeated application of 5.

\end{proof}

\ifediting
\bibliographystyle{alpha}
 \bibliography{incidence}

\end{document}
\fi

\newif\ifediting

\ifediting
\documentclass[11pt]{article}

\usepackage{amsmath,amsthm,amssymb}
\newcommand{\remove}[1]{}
\setlength{\topmargin}{0.3in} \setlength{\headheight}{0in}
\setlength{\headsep}{0in} \setlength{\textheight}{8.0in}
\setlength{\topsep}{0.1in} \setlength{\itemsep}{0.0in}
\parskip=0.05in
 \textwidth=6.5in 
\oddsidemargin=0truecm \evensidemargin=0truecm

\newtheorem{thm}{Theorem}[section]
\newtheorem{claim}[thm]{Claim}
\newtheorem{lem}[thm]{Lemma}
\newtheorem{define}[thm]{Definition}
\newtheorem{cor}[thm]{Corollary}
\newtheorem{obs}[thm]{Observation}
\newtheorem{example}[thm]{Example}
\newtheorem{construct}[thm]{Construction}
\newtheorem{conjecture}[thm]{Conjecture}
\newtheorem{THM}{Theorem}
\newtheorem{question}{Question}
\newtheorem{fact}[thm]{Fact}
\newtheorem{prop}[thm]{Proposition}

\def\F{{\mathbb{F}}}
\def\Q{{\mathbb{Q}}}
\def\Z{{\mathbb{Z}}}
\def\N{{\mathbb{N}}}
\def\R{{\mathbb{R}}}
\def\K{{\mathbb{K}}}
\def\C{{\mathbb{C}}}
\def\A{{\mathbb{A}}}
\def\P{{\mathbb{P}}}
\def\cP{{\cal P}}
\def\cS{{\mathcal S}}
\def\cE{{\mathcal E}}
\def\V{{\mathbf{V}}}
\def\I{{\mathbf{I}}}
\def\bx{{\mathbf x}}
\def\by{{\mathbf y}}
\def\E{{\mathbb E}}

\def\half{ \frac{1}{2}}
\newcommand{\ip}[2]{\langle #1,#2 \rangle}
\def\sumN{\sum_{i=1}^n}
\def\_{\,\,\,\,\,}
\def\prob{{\mathbf{Pr}}}
\newcommand{\entropy}[1]{ {\text{H}_{\infty}\left({#1}\right)} }
\def\modulo{\text{mod}}
\def\omm{ \{0,1\} }
\def\id{ \textit{id} }

\def\D{{\partial}}
\def\gap{\textsf{gap}}
\def\sign{\textsf{sign}}
\def\spar{\textsf{sparse}}
\def\span{\textsf{span}}
\def\Part{\textbf{Part}}
\def\Mon{\textbf{Mon}}
\def\sing{\textbf{sing}}
\def\Und{\textsf{Und}}
\def\Comp{\textsf{Comp}}
\def\rank{\textsf{rank}}
\def\poly{\textsf{poly}}
\def\codim{\textsf{codim}}
\def\cp{\textsf{cp}}
\def\uni{\textsf{Uni}}
\def\ext{\textsf{\bf Ext}}
\def\extt{\textsf{\bf Ext2}}

\def\fplus{{\,+_f\,}}

\newcommand{\epsclose}{\stackrel{\epsilon}{\thicksim} }
\newcommand{\eclose}[1]{\stackrel{{#1}}{\thicksim} }
\newcommand{\eps}{\epsilon}
\newcommand{\Anote}[1]{\begin{quote}{\sf Avi's Note:} {\sl{#1}} \end{quote}}

\begin{document}

\title{Incidence Theorems -- Lecture Notes}
\date{}
\maketitle


\fi

\section{Growth in $\F_p$}\label{sec-growth}

Let $\F = \F_p$ be a finite field of prime cardinality. Recall that such fields do not contain any subfields. For simplicity, we will talk about prime fields but all of our arguments can be extended to fields not containing large subfields. Our goal in this section will be to show that for a set $A \subset \F$ and for almost all values of $\lambda \in \F$ we have $|A + \lambda A | \gg |A|$. Here we denote by $\lambda A = \{ \lambda a \,|\, a \in A \}$ (do not confuse this with $k \cdot A$ used for iterated sums of elements in $A$). Later we will need to develop a `distributional' variant of the same statement.

Our first step is to show that there is at least {\em one} good value of $\lambda$.
\begin{lem}\label{lem-onegood}
Let $A \subset \F$ then there exists $\lambda \in \F$ such that $$ |A+\lambda A| \geq \frac{1}{2} \cdot \min\{ |A|^2, p \}.$$
\end{lem}
\begin{proof}
We will use the familiar notation $Q(A,B)$ for the set of quadruples $a+b=a'+b'$ with $a,a' \in A$ and $b,b' \in B$. Recall also that $$ |Q(A,B)| \geq \frac{|A|^2 |B|^2}{|A+B|}.$$ Summing over $\lambda$ we get
$$ \sum_{\lambda \neq 0} |Q(A,\lambda A)| = \left|\left\{ (a_1,a_2,a_3,a_4,\lambda) \in A^4 \times \F^* \,|\, a_1 + \lambda a_2 = a_3 + \lambda a_4. \right\}\right|.$$ The solutions with $a_1=a_3$ and $a_2 = a_4$ contribute at most $|A|^2 (p-1)$ to this sum (since $\lambda$ can be anything). The solutions with $(a_1,a_2) \neq (a_3,a_4)$ determine a unique $\lambda$ and so contribute at most $|A|^2 (|A|^2- 1)$. Over all we have $$ \sum_{\lambda \neq 0} |Q(A,\lambda A)| \leq |A|^2 (p-1) + |A|^2 (|A|^2- 1). $$ This means that there exists $\lambda_0$  such that $|Q(A,\lambda_0 A)| \leq |A|^2 + |A|^4/(p-1)$. This implies the required bound on $|A + \lambda_0 A|$.
\end{proof}

Let $$ \stab_K(A) = \{ \lambda \in \F^* \,|\, |A+\lambda A| \leq K|A|\}.$$ The next lemma shows that the set $\stab_K(A)$ behaves somewhat similarly to a sub field. 
\begin{lem}\label{lem-stab}
There exists an absolute constant $c$ such that:
\begin{enumerate}
	\item If $\lambda \in \stab_K(A)$ then $-\lambda, 1/\lambda$ are in $\stab_{K^c}(A)$.
	\item If $\lambda_1,\lambda_2 \in \stab_K(A)$ then $\lambda_1\lambda_2, \lambda_1+\lambda_2$ are in $\stab_{K^c}(A)$.
\end{enumerate}
\end{lem}
\begin{proof} The proof is an immediate application of Ruzsa calculus (see last section):
\noindent{1.} The claim about $-\lambda$ follows from Ruzsa calculus. The claim about $1/\lambda$ is trivial since $|A+\lambda A| = |(1/\lambda)A + A| $.
\noindent{2.} Using Ruzsa calculus we have $$|A+(\lambda_1  + \lambda_2 )A| \leq |A + \lambda_1 A + \lambda_2 A| \leq K^c |A|.$$ To argue about the product observe that $|A+ \lambda_1 A| \leq K |A|$ and $|A + (1/\lambda_2)A| \leq K |A|$ so, by Ruzsa, we have $$|A + (\lambda_1 \lambda_2)A| = |\lambda_1 A + (1/\lambda_2) A| \leq K^c |A|.$$
\end{proof}

We would like to argue that, if $\stab_K(A)$ is large, then for some $c$, $\stab_{K^c}(A)$ contains all of $\F$ (contradicting Lemma~\ref{lem-onegood}). This will be obtained by the following lemma.
\begin{lem}\label{lem-growth}
Let $A \subset \F$ then $$ |3 \cdot A^2 - 3 \cdot A^2 | \geq \frac{1}{2} \min \{ |A|^2, p \},$$ where $$3 \cdot A^2 - 3 \cdot A^2 = A \cdot A + A \cdot A + A \cdot A - A \cdot A - A \cdot A - A \cdot A .$$
\end{lem}
\begin{proof}
Observe that, if $\lambda \not\in \frac{A-A}{A-A}$ then $|A + \lambda A| = |A|^2$. We divide the proof into two cases:
\paragraph{Case 1: $\frac{A-A}{A-A} \neq \F$.} In this case there must exist $\lambda \in \frac{A-A}{A-A}$ such that $\lambda+1 \not\in \frac{A-A}{A-A}$ (here we use the particular structure of the field $\F_p$). This implies $|A + (\lambda+1)A| = |A|^2$. Write $\lambda = \frac{a_1-a_3}{a_2-a_4}$. We have $$ (a_2 - a_4)(A+ (\lambda+1)A) \subset (a_2-a_4)A + (a_1-a_3 + a_2 - a_4)A \subset 3 \cdot A^2 - 3 \cdot A^2.$$ And since the size of the set on the left is $|A|^2$ we are done.

\paragraph{Case 2: $\frac{A-A}{A-A} = \F$.} Then, from Lemma~\ref{lem-onegood} we have that there exists $\lambda \in \frac{A-A}{A-A}$ such that $|A + \lambda A| \geq \frac{1}{2}\min \{ |A|^2, p \}.$ Write $\lambda = \frac{a_1-a_3}{a_2-a_4}$. Then $$ (A+\lambda A)(a_2 - a_4) \subset 3 \cdot A^2 - 3 \cdot A^2$$ and the bound follows also in this case. 
\end{proof}

An immediate corollary of Lemma~\ref{lem-growth} is the following
\begin{cor}\label{cor-iterate}
Let $A \subset \F$ be of size $p^\delta$. Then $|k \cdot A^k - k \cdot A^k| = \F$ for some $k = k(\delta)$.
\end{cor}
\begin{proof}
Applying Lemma~\ref{lem-growth} gets us all the way up to size $p/2$. To make the final jump observe that, if $|A| > p/2$ ($p$ is odd!) then $A \cap (x-A) \neq \emptyset$ for all $x \in \F$. This means that $A+A = \F$ and so one more addition will finish the job.
\end{proof}

Combining the above   we get to our goal:
\begin{thm}\label{thm-growth}
Let $A,T \subset \F$ with  $p^\alpha \leq |A| \leq p^{1-\alpha}$ and $|T| \geq p^\beta$. Then there exists $\lambda \in T$ such that $|A + \lambda A| \geq |A|^{1 + c(\alpha,\beta)}, $ where $c(\alpha,\beta)>0$ is a constant depending only on $\alpha$ and $\beta$.
\end{thm}
\begin{proof}
We let $K = |A|^{c(\alpha,\beta)}$. If the theorem is not true then $T \subset \stab_K(A)$ which implies that taking some $k = k(\alpha,\beta)$ sums and products (as in Corollary~\ref{cor-iterate}) we will have $\F = \stab_{K'}(A)$ with $K' \leq  K^{k(\alpha,\beta)}$. Picking $c(\alpha,\beta)$ small enough we get that for all $\lambda \in \F$, $|A + \lambda A| \leq K' |A| \ll \min\{ |A|^2, p \}.$ This contradicts Lemma~\ref{lem-onegood}.
\end{proof}

Our next goal will be to prove a `distributional' version of Theorem~\ref{thm-growth}.

\ifediting
\bibliographystyle{alpha}
 \bibliography{incidence}

\end{document}
\fi

\newif\ifediting

\ifediting
\documentclass[11pt]{article}

\usepackage{amsmath,amsthm,amssymb}
\newcommand{\remove}[1]{}
\setlength{\topmargin}{0.3in} \setlength{\headheight}{0in}
\setlength{\headsep}{0in} \setlength{\textheight}{8.0in}
\setlength{\topsep}{0.1in} \setlength{\itemsep}{0.0in}
\parskip=0.05in
 \textwidth=6.5in 
\oddsidemargin=0truecm \evensidemargin=0truecm

\newtheorem{thm}{Theorem}[section]
\newtheorem{claim}[thm]{Claim}
\newtheorem{lem}[thm]{Lemma}
\newtheorem{define}[thm]{Definition}
\newtheorem{cor}[thm]{Corollary}
\newtheorem{obs}[thm]{Observation}
\newtheorem{example}[thm]{Example}
\newtheorem{construct}[thm]{Construction}
\newtheorem{conjecture}[thm]{Conjecture}
\newtheorem{THM}{Theorem}
\newtheorem{question}{Question}
\newtheorem{fact}[thm]{Fact}
\newtheorem{prop}[thm]{Proposition}

\def\F{{\mathbb{F}}}
\def\Q{{\mathbb{Q}}}
\def\Z{{\mathbb{Z}}}
\def\N{{\mathbb{N}}}
\def\R{{\mathbb{R}}}
\def\K{{\mathbb{K}}}
\def\C{{\mathbb{C}}}
\def\A{{\mathbb{A}}}
\def\P{{\mathbb{P}}}
\def\cP{{\cal P}}
\def\cS{{\mathcal S}}
\def\cE{{\mathcal E}}
\def\V{{\mathbf{V}}}
\def\I{{\mathbf{I}}}
\def\bx{{\mathbf x}}
\def\by{{\mathbf y}}
\def\E{{\mathbb E}}

\def\half{ \frac{1}{2}}
\newcommand{\ip}[2]{\langle #1,#2 \rangle}
\def\sumN{\sum_{i=1}^n}
\def\_{\,\,\,\,\,}
\def\prob{{\mathbf{Pr}}}
\newcommand{\entropy}[1]{ {\text{H}_{\infty}\left({#1}\right)} }
\def\modulo{\text{mod}}
\def\omm{ \{0,1\} }
\def\id{ \textit{id} }

\def\D{{\partial}}
\def\gap{\textsf{gap}}
\def\sign{\textsf{sign}}
\def\spar{\textsf{sparse}}
\def\span{\textsf{span}}
\def\Part{\textbf{Part}}
\def\Mon{\textbf{Mon}}
\def\sing{\textbf{sing}}
\def\Und{\textsf{Und}}
\def\Comp{\textsf{Comp}}
\def\rank{\textsf{rank}}
\def\poly{\textsf{poly}}
\def\codim{\textsf{codim}}
\def\cp{\textsf{cp}}
\def\uni{\textsf{Uni}}
\def\ext{\textsf{\bf Ext}}
\def\extt{\textsf{\bf Ext2}}

\def\fplus{{\,+_f\,}}

\newcommand{\epsclose}{\stackrel{\epsilon}{\thicksim} }
\newcommand{\eclose}[1]{\stackrel{{#1}}{\thicksim} }
\newcommand{\eps}{\epsilon}
\newcommand{\Anote}[1]{\begin{quote}{\sf Avi's Note:} {\sl{#1}} \end{quote}}

\begin{document}

\title{Incidence Theorems -- Lecture Notes}
\date{}
\maketitle


\fi

\section{The Balog-Szemeredi-Gowers theorem}\label{sec-bsg}
In the previous section we showed:
\begin{thm}\label{thm-growth2}
Let $A,T \subset \F$ with  $p^\alpha \leq |A| \leq p^{1-\alpha}$ and $|T| \geq p^\beta$. Then there exists $\lambda \in T$ such that $|A + \lambda A| \geq |A|^{1 + c(\alpha,\beta)}, $ where $c(\alpha,\beta)>0$ is a constant depending only on $\alpha$ and $\beta$.
\end{thm}

Let's try to see why this is the kind of result we could hope to use to prove the ST theorem and why it is not really strong enough. We will demonstrate this by considering a very special case of a line/point arrangement. Let $P$ and $L$ be sets of $N$ points and $N$ lines in $\F^2$ with $N \sim p^\delta$ for some `nice' $\delta$ (say, between $1/4$ and $7/4$). Suppose also that $P = P_x \times P_y$ with $|P_x|, |P_y| \lesssim N^{1/2+\eps}$  and that the lines in $L$ are given by equations of the form $Y = aX+b$ with $a,b \in A$ and $|A| \leq N^{1/2+\eps}$. Then, if $|I(P,L)| > N^{3/2-\eps}$, then, for at least $N^{1-\eps}$ lines $Y = aX + b$ there will be at least $N^{1/2-\eps}$ values of $x \in P_x$ for which $ax + b \in P_y$ (this is the `typical' value required to obtain $N^{3/2-\eps}$ intersections). This means that $A + xA$ is `small' (contained in $P_y$) for many values of $x$. This would contradict Theorem~\ref{thm-growth2} if the information was complete (i.e., if we knew that for {\em all} $a,b \in A$ and $x \in P_x$, $a+xa \in P_y$). However, the information is given in an incomplete form, as a quantitative incidence bound, and so we need a stronger version of Theorem~\ref{thm-growth2} that can handle such information.

The idea is to work with $Q(A,B)$ instead of $|A+B|$. Recall that $$ |A+B| \geq \frac{|A|^2|B|^2}{|Q(A,B)|}.$$ Thus, if we define $$ E(A,B) =  \frac{|A|^2|B|^2}{|Q(A,B)|}$$ we have $\max\{ |A|, |B|\} \leq E(A,B) \leq |A+B|.$ We will call the quantity $E(A,B)$ the {\em additive energy} (or just {\em energy}) of $A+B$, as it relates to the $\ell_2$ norm of the distribution obtained by sampling $a,b$ independently in $A,B$ and summing them\footnote{We normalize $E(A,B)$ so that it is in the same scale as $|A+B|$. In some texts other scalings are used, e.g. in \cite{GreenCourse} the scaling is so that $E(A,B)$ is in the range $[0,1]$. In some places, the term additive energy is used for the quantity $Q(A,B)$.   }. As an example, consider an arithmetic progression $A$ of size $N$ and notice that, in this case, both $|A+A|$ and $E(A,A)$ are bounded by $\lesssim N$. Now, let $B$ be a set of size $N$ with $|B+B| = (1/2)|B|(|B|-1)$ (i.e., a set with no dependencies). Here we also have $|B+B| \sim E(B,B) \sim |B|^2.$ However, taking $C = A \cup B$ we get that $|C+C| \gtrsim |B|^2 \gtrsim |C|^2$ but $E(C,C) \lesssim N$ (because $|Q(C,C)| \geq |Q(A,A)|$). That is, the energy $E(A,B)$ can capture information about sufficiently large subsets of $A$ (or $B$) that do not grow in addition -- information that is not captured by $|A+B|$. A partial converse to this statement is given by the following important result known as the Balog-Szemeredi-Gowers Theorem (or BSG for short).

\begin{thm}[\cite{BalSze,Gow}]  \label{thm-bsg}
Let $A,B \subset G$ be sets of size $N$ in an abelian group $G$. Suppose that $E(A,B) \leq KN$. Then, there exist subsets $A' \subset A$ and $B' \subset B$ with $|A'|,|B'| \geq N/K^c$ and with $|A'+B'| \leq K^c N$. Here, $c >0$ is some absolute constant.
\end{thm}

This theorem will allow us (with some work) to derive an additive energy version of Theorem~\ref{thm-growth2} with $|A+\lambda A|$ replaced by $E(A,\lambda A)$. The BSG theorem will actually follow from a relatively generic graph theoretic lemma which we now state.

\begin{lem}\label{lem-3path}
Let $H \subset V \times U$ be a bipartite graph with $|V|=|U|=N$. Suppose $|H| \geq \alpha N^2$ (the number of edges). Then, there are subsets $V' \subset V$ and $U' \subset U$ with $|V'|,|U'| \geq \alpha^c N$ and such that for all $v \in V'$, $u \in U'$ there are at least $\alpha^c N^2$ paths of length three between $v$ and $u$. 
\end{lem}

Before proving Lemma~\ref{lem-3path}, let's see how it implies the BSG theorem.
\paragraph{Proof of Theorem~\ref{thm-bsg}:} Suppose $E(A,B) \leq KN$. Then $|Q(A,B)| \geq N^3/K.$ Let $P$ be the set of values $x$ with $|R(x)| = |\{ (a,b) \in A \times B \,|\, a-b=x \}| \geq N/2K.$ This is the set of `popular differences' and can be seen to have size at least $|P| \geq N/2K$ (we saw this argument last time). Consider the graph $H \subset A \times B$ whose edges corresponds to pairs $(a,b)$ with $a-b \in P$ and label each such edge $(a,b)$ with the value $a-b$. From the definition of $P$ we have that $H$ has at least $|P|(N/2K) \geq (1/4K^2) \cdot N^2$ edges and so, applying Lemma~\ref{lem-3path} with $\alpha = 1/4K^2$, we have subsets $A' \subset A$ and $B' \subset B$
 with $|A'|,|B'| \geq N/K^c$ and such that, for all $(a,b) \in A' \times B'$ there are at least $N^2/K^c$ paths of length three between $a$ and $b$ in the graph $H$. Consider such a path $a \rightarrow b' \rightarrow a' \rightarrow b$. Writing $$ a-b = (a-b') - (a'-b') + (a'-b) $$ and using the fact that all three differences in this sum are popular, we see that $a-b$ can be written in at least $ N^2/K^c$ distinct ways as $a-b = x_1 - x_2 + x_3$ with $x_1,x_2,x_3 \in P$. This implies $$ |A'-B'| \leq \frac{K^c|P|^3}{N^2} \lesssim  K^{c'} N$$ for some other absolute constant $c' > 0$. Going from $|A' - B'|$ to $|A' + B'|$ is possible using Ruzsa calculus and loses another constant. \qed

To prove Lemma~\ref{lem-3path} we first prove a simpler lemma on paths of length two. We will denote by $\Gamma(S)$ the set of neighbors of $S$ in the graph $H$, where $S$ is some subset of the vertices of either $V$ or $U$.
\begin{lem}\label{lem-2path}
Let $H \subset V \times U$ be as in Lemma~\ref{lem-3path}. Then, for every $\eps>0$ there exists a set $V' \subset V$ with $|V'| \geq (\alpha/\sqrt{2})N$ and s.t $$|\{ (v_1,v_2) \in V' \times V' \,|\, |\Gamma(v_1) \cap \Gamma(v_2)| \leq (\eps\alpha^2/2)N \}| \leq \eps|V'|^2.$$ In other words, for all but an $\eps$ fraction of the pairs of vertices in $V'$, the pair will have at least $(\eps \alpha^2/2)N$ common neighbors (or paths of length two).
\end{lem}
\begin{proof}
The proof uses a clever trick introduced by Gowers which combines a `somewhat' random choice of the set $V'$. Picking the set $V'$ completely at random does not seem to work. The idea is to chose $V'$ as the set of neighbors $\Gamma(u)$ of a random vertex $u \in U$. This makes sense, since a pair $(v_1,v_2)$ with few common neighbors are less likely to be in $V'$ than a pair that has many common neighbors. Lets see the calculation.

Denote the set of `bad' pairs by $$ S = \{ (v_1,v_2) \in V \times V \,|\, |\Gamma(v_1) \cap \Gamma(v_2)| \leq (\eps\alpha^2/2)N \}. $$ For each $u \in U$ let $S_u = S \cap \Gamma(u)$ denote the set of bad pairs among the neighbors of $u$. Suppose we pick $u$ at random and consider first the expectation of $|\Gamma(u)|^2$ (the total number of pairs among neighbors of $u$). Using Cauchy-Schwarz we have:
\begin{equation*}\label{eq-Gammau}
\E_u [ |\Gamma(u)|^2] \geq \left( \E_u [ |\Gamma(u)|]\right)^2 = \alpha^2 N^2.
\end{equation*}
We also have
\begin{eqnarray*}
	\E_u[|S_u|] &=& \E_u\left[  \sum_{v_1,v_2} 1_{v_1,v_2 \in S} \cdot 1_{(v_1,u) \in H}\cdot 1_{(v_2,u) \in H}\right] \\
	&=& \sum_{v_1,v_2} 1_{v_1,v_2 \in S} \cdot \frac{|\Gamma(v_1) \cap \Gamma(v_2)|}{N} \\
	&\leq& N^2 \cdot (\eps \alpha^2/2).
\end{eqnarray*}
Combining the two bounds we get
$$ \E_u\left[ \eps|\Gamma(u)|^2 - |S_u| \right] \geq (\eps \alpha^2 /2) N^2. $$
This implies  $|\Gamma(u)| \geq (\alpha/\sqrt 2)N$ and $|S_u| \leq \eps |\Gamma(u)|^2$ as was required.
\end{proof}

\paragraph{Proof of Lemma~\ref{lem-3path}:} We will omit some of the detailed calculations (which can be easily filled in). By throwing away a small fraction of the vertices we can reduce to the case where the minimum degree of a vertex is at least $(\alpha/2)N$. Let $V' \subset V$ be given by Lemma~\ref{lem-2path} so that $|V'| \gtrsim \alpha N$ and such that for all but $\eps |V'|^2$ pairs $(v_1,v_2) \in V'$ we have $|\Gamma(v_1) \cap \Gamma(v_2)| \gtrsim \eps \alpha^2 N$ (we will pick $\eps$ later). Notice that there might be some vertices $v_1 \in V'$ for which there are many (even all) vertices $v_2 \in V'$ that have few common neighbors with them. We can, however, find a subset $V'' \subset V'$ with $|V''| \sim |V'|$ and such that for all $v_1 \in V''$ there are at most $2\eps|V'|$ vertices $v_2 \in V'$ with $|\Gamma(v_1) \cap \Gamma(v_2) | \lesssim \eps \alpha^2 N$. Next, we can find a subset $U' \subset U$ with $|U'| \gtrsim \alpha^2 N$ such that each $u \in U'$ has at least $10 \eps |V'|$ neighbors in $V''$ (here we need to choose $\eps$ sufficiently small, but still polynomial in $\alpha$). This part uses the fact that the minimum degree is large and so there is a quadratic number of edges leaving $V''$. Now, fix $u \in U'$ and $v \in V''$.  We  will build many paths of length three between $u$ and $v$ as follows: Start with $u$ and move to a neighbor of $u$ in $V''$. There are at least $10\eps |V'|$ options for this step and at most $2\eps|V'|$ of them will have less than $\lesssim \eps \alpha^2 N$ common neighbors with $v$ (this is how we defined $V''$). This means that we can complete this path in at least $\gtrsim \eps \alpha^2 N$ ways. This gives $\gtrsim \alpha^c N^2$ distinct paths of length three.

\subsection{Energy version of growth in $\F_p$}

We will now use the BSG theorem to derive an energy  version of Theorem~\ref{thm-growth2}. The proof will follow from the following result (due to Bourgain \cite{BouML}). 
\begin{thm}\label{thm-bourgain}
Let $\F = \F_p$ with $p$ prime. Let $A \subset \F$ and $T \subset \F^*$. Suppose that for all $\lambda \in T$ we have $E(A,\lambda A) \leq K|A|$. Then, there exist $A' \subset A$ and $T' \subset xT$ (for some $x \in \F^*$) such that $|A'| \geq |A|/K^c$, $|T'| \geq |T|/K^c$ and with $|A' + \lambda A'| \leq K^c |A'|$ for all $\lambda \in T'$.
\end{thm}
\begin{proof}
Using the BSG theorem (Theorem~\ref{thm-bsg}) we can find sets $X_\lambda,Y_\lambda \subset A$ for each $\lambda \in T$ such that $|X_\lambda|,|Y_\lambda| \geq |A|/K^c$ and such that $|X_\lambda + \lambda Y_\lambda| \leq K^c |A|$ for all $\lambda \in T$. We will want to somehow `paste' many of these sets together. We start with a simple claim.
\begin{claim}\label{cla-subsets}
Let $S_1,\ldots,S_k \subset S$ be finite sets with $|S_i| \geq \delta|S|$ for all $i \in [k]$. Then, there exists $i \in [k]$ such that $$|\{ j \in [k] \,|\, |S_i \cap S_j| \geq (\delta^2/2)|S| \}| \geq (\delta^2/2)k.$$
\end{claim}
\begin{proof}
Observe that
\begin{eqnarray*}
\sum_{i,j} |S_i \cap S_j| &=& \sum_{i,j} \sum_{x \in |S|} 1_{x \in S_i} \cdot 1_{x \in S_j} \\
&=& \sum_x |\{ i \,|\, x \in S_i\}|^2 \\
&\geq& \frac{1}{|S|} \cdot \left( \sum_{x \in S}|\{i \,|\, x \in S_i \}| \right)^2 \\
&=& \frac{ \left( \sum_i |S_i| \right)^2}{|S|} \geq k^2 \delta^2 |S|.
\end{eqnarray*}
If we take $i \in [k]$ such that $$ \sum_j |S_i \cap S_j| \geq k \delta^2 |S|$$ we will get the required property.
\end{proof}

Using the claim we can find some $\lambda_0 \in T$ and a subset $T' \subset T$ with $|T'| \gtrsim |T|/K^c$ (remember our convention to `reuse' the constant $c$) such that for all $\lambda \in T'$ we have $|X_{\lambda_0} \cap X_\lambda|, |Y_{\lambda_0} \cap Y_\lambda| \gtrsim |A|/K^c$. Notice that, to get this to work, we need to apply the claim with $S = A \times A$ and the family of sets $S_\lambda = X_\lambda \times Y_\lambda$. We find $\lambda_0$ such that $S_{\lambda_0}$ has intersection at least $|A|^2/K^c$ with  $S_\lambda$ for all $\lambda$ in some large set $T'$ and then argue about the intersections of the projections $X_{\lambda_0}, Y_{\lambda_0}$.

In what follows we will use Ruzsa calculus (RC) very freely (without stating each time exactly what part we are using) and the reader is advised to recall the different claims involved. We will use the notation $X \equiv Y$ to mean $|X+Y| \lesssim K^c |A|$ for some absolute constant $c$ (this notation is only for this proof). We know that $X_\lambda \equiv \lambda Y_\lambda$ for all $\lambda \in T$. Thus $X_\lambda \equiv X_\lambda$ for all $\lambda$ and in particular $X_{\lambda_0} \equiv X_{\lambda_0}$. Using RC and the fact that $|X_{\lambda_0} \cap X_\lambda|$ is large for all $\lambda \in T'$ we get that $X_{\lambda_0} \equiv \lambda Y_{\lambda}$ for all $\lambda \in T'$. In the same way, since $Y_\lambda \cap Y_{\lambda_0}$ is large, we get that $X_{\lambda_0} \equiv \lambda Y_{\lambda_0}$ for all $\lambda \in T'$. Using the triangle inequality, we get $\lambda_0 Y_{\lambda_0} \equiv \lambda Y_{\lambda_0}$ for all $\lambda \in T'$ which is the same as $Y_{\lambda_0} \equiv \frac{\lambda}{\lambda_0} Y_{\lambda_0}.$ Set $T'' = (1/\lambda_0)T'$ and $A' = Y_{\lambda_0}$ and the theorem follows.
\end{proof}
 
Combining Theorem~\ref{thm-bourgain} with Theorem~\ref{thm-growth2} we immediately get to our previously described goal:
\begin{cor}[Growth in energy]\label{cor-entropy}
	Let $A,T \subset \F$ with  $p^\alpha \leq |A| \leq p^{1-\alpha}$ and $|T| \geq p^\beta$. Then there exists $\lambda \in T$ such that $E(A,\lambda A) \geq |A|^{1 + c(\alpha,\beta)}, $ where $c(\alpha,\beta)>0$ is a constant depending only on $\alpha$ and $\beta$.
\end{cor}

\ifediting
\bibliographystyle{alpha}
 \bibliography{incidence}

\end{document}
\fi

\newif\ifediting

\ifediting
\documentclass[11pt]{article}

\usepackage{amsmath,amsthm,amssymb}
\newcommand{\remove}[1]{}
\setlength{\topmargin}{0.3in} \setlength{\headheight}{0in}
\setlength{\headsep}{0in} \setlength{\textheight}{8.0in}
\setlength{\topsep}{0.1in} \setlength{\itemsep}{0.0in}
\parskip=0.05in
 \textwidth=6.5in 
\oddsidemargin=0truecm \evensidemargin=0truecm

\newtheorem{thm}{Theorem}[section]
\newtheorem{claim}[thm]{Claim}
\newtheorem{lem}[thm]{Lemma}
\newtheorem{define}[thm]{Definition}
\newtheorem{cor}[thm]{Corollary}
\newtheorem{obs}[thm]{Observation}
\newtheorem{example}[thm]{Example}
\newtheorem{construct}[thm]{Construction}
\newtheorem{conjecture}[thm]{Conjecture}
\newtheorem{THM}{Theorem}
\newtheorem{question}{Question}
\newtheorem{fact}[thm]{Fact}
\newtheorem{prop}[thm]{Proposition}

\def\F{{\mathbb{F}}}
\def\Q{{\mathbb{Q}}}
\def\Z{{\mathbb{Z}}}
\def\N{{\mathbb{N}}}
\def\R{{\mathbb{R}}}
\def\K{{\mathbb{K}}}
\def\C{{\mathbb{C}}}
\def\A{{\mathbb{A}}}
\def\P{{\mathbb{P}}}
\def\cP{{\cal P}}
\def\cS{{\mathcal S}}
\def\cE{{\mathcal E}}
\def\V{{\mathbf{V}}}
\def\I{{\mathbf{I}}}
\def\bx{{\mathbf x}}
\def\by{{\mathbf y}}
\def\E{{\mathbb E}}

\def\half{ \frac{1}{2}}
\newcommand{\ip}[2]{\langle #1,#2 \rangle}
\def\sumN{\sum_{i=1}^n}
\def\_{\,\,\,\,\,}
\def\prob{{\mathbf{Pr}}}
\newcommand{\entropy}[1]{ {\text{H}_{\infty}\left({#1}\right)} }
\def\modulo{\text{mod}}
\def\omm{ \{0,1\} }
\def\id{ \textit{id} }

\def\D{{\partial}}
\def\gap{\textsf{gap}}
\def\sign{\textsf{sign}}
\def\spar{\textsf{sparse}}
\def\span{\textsf{span}}
\def\Part{\textbf{Part}}
\def\Mon{\textbf{Mon}}
\def\sing{\textbf{sing}}
\def\Und{\textsf{Und}}
\def\Comp{\textsf{Comp}}
\def\rank{\textsf{rank}}
\def\poly{\textsf{poly}}
\def\codim{\textsf{codim}}
\def\cp{\textsf{cp}}
\def\uni{\textsf{Uni}}
\def\ext{\textsf{\bf Ext}}
\def\extt{\textsf{\bf Ext2}}

\def\fplus{{\,+_f\,}}

\newcommand{\epsclose}{\stackrel{\epsilon}{\thicksim} }
\newcommand{\eclose}[1]{\stackrel{{#1}}{\thicksim} }
\newcommand{\eps}{\epsilon}
\newcommand{\Anote}[1]{\begin{quote}{\sf Avi's Note:} {\sl{#1}} \end{quote}}

\begin{document}

\title{Incidence Theorems -- Lecture Notes}
\date{}
\maketitle


\fi

\section{Szemeredi-Trotter in finite fields}\label{sec-stfinite}

We will now see how to use Corollary~\ref{cor-entropy} to give a non trivial bound of $N^{3/2 - \eps}$ for some constant $\eps>0$ on the number of incidences of $N$ points and $N$ lines in $\F^2$. We wish to prove:
\begin{thm}[ST over finite fields \cite{BKT04}]\label{thm-stfinite}
	Let $L$ be a set of $N$ lines in $\F^2$ and let $P$ be a set of $N$ points in $\F^2$. Then, if $p^\alpha < N < p^{2-\alpha}$ for some $\alpha>0$ then $|I(P,L)| \lesssim N^{3/2 - \eps}$, where $\eps>0$ depends only on $\alpha$.
\end{thm}

The proof will be in two steps:
\begin{enumerate}
	\item Reduce the problem to the case when the $N$ points are contained in an $N^{1/2}$ by $N^{1/2}$ grid $A \times B \subset \F^2$.
	\item Prove the required bound over a grid using Corollary~\ref{cor-entropy}.
\end{enumerate}

We will use the following notations: for a point $p \in P$ denote by $L(p) = \{ \ell \in L \,|\, p \in \ell\}$ and for a line $\ell \in L$ denote $P(\ell) = \{ p \in P \,|\, p \in \ell\}$. Suppose $P,L$ are such that $|I(P,L)| \gg N^{3/2 - \eps}$ where $\eps$ will be chosen small enough to derive a contradiction later. We start with throwing away some lines/points to ensure certain regularity conditions. First, remove all lines with at most $N^{1/2 - 2\eps}$ points on them. This can reduce the number of incidences by a negligible fraction. Since there are  still $\gtrsim N^{3/2-\eps}$ incidences we must have at least $N^{1-\eps}$ lines left after this step (otherwise use Cauchy-Schwarz to bound the number of incidences). Next, remove all lines that have {\em at least} $N^{1/2 + 2\eps}$ points on them. Recall that by Cauchy-Schwarz we have a bound of $N^{3/2}$ on the number of incidences and so in this second step we will remove at most $N^{1-2\eps}$ lines, which means  that we will still have at least $\gtrsim N^{1-\eps}$ lines with at least $N^{1/2 - 2\eps}$ points on each. Thus, the total number of incidences will remain at least $N^{3/2 - 3\eps}$. At this point we have that for each line $\ell \in L$,  $$N^{1/2 - 2\eps} \lesssim |P(\ell)| \lesssim N^{1/2 + 2\eps}.$$ We can perform the same procedure on points and obtain that for all $p \in P$, $N^{1/2 - 2\eps} \lesssim |L(p)| \lesssim N^{1/2 + 2\eps}$. Since we are only removing points we will still have after this step the bound $|P(\ell)| \leq N^{1/2 + 2\eps}$ for each line (the lower bound might not hold).

\subsection{Translating the problem to a grid}

To translate our problem to a grid we first find two points $p_0,p_1 \in P$ such that most incidences happen on intersections of a line through $p_0 $ and a line through $p_1$. 
\begin{claim}\label{cla-pregrid}
There exist points $p_0$ and $p_1$ in $P$ such that there exist a subset $P' \subset P$ with $|P'| \geq N^{1 - c\eps}$ and s.t $P' \subset \{ \ell_0 \cap \ell_1 \,|\, \ell_0 \in L(p_0), \ell_1 \in L(p_1) \}$. Here $c >0$ is some absolute constant.
\end{claim}
\begin{proof}
For $p \in P$ define the set of points that lie on some line through $p$ to be $$\Gamma_2(p) = \{ p' \in P \,|\, \exists \ell \in L \text{ s.t } p,p' \in \ell\},$$ (these are vertices of distance two from $p$ on the incidence graph). We are looking for a pair $p_0,p_1$ with large $|\Gamma_2(p_0)\cap \Gamma_2(p_1)|$. To find them we will consider the expected size of 
\begin{eqnarray*}
\E_{p_0,p_1}\left[ |\Gamma_2(p_0)\cap \Gamma_2(p_1)| \right] &=& \frac{1}{N^2} \sum_{p_0,p_1 \in P} \sum_{q \in P} \sum_{\ell_0,\ell_1 \in L(q)} 1_{p_0 \in \ell_0}\cdot 1_{p_1 \in \ell_1} \\
&=& \frac{1}{N^2} \sum_{q \in P} \sum_{\ell_0,\ell_1 \in L(q)} |P(\ell_0)|\cdot |P(\ell_1)| \\
&=& \frac{1}{N^2} \sum_{q \in P} \left( \sum_{\ell \in L(q)} |P(\ell)| \right)^2 \\
&\geq& \frac{1}{N^3} \left( \sum_{q\in P} \sum_{\ell \in L(q)} |P(\ell)|\right)^2 \\
&=& \frac{1}{N^3} \left( \sum_{\ell \in L} |P(\ell)|^2 \right)^2 \\
&\geq& \frac{1}{N^5} \left( \sum_{\ell \in L} |P(\ell)| \right)^4 \\
&\geq& \frac{1}{N^5}\left( N^{3/2 - 2\eps}\right)^4 \geq N^{1 - c\eps},
\end{eqnarray*}
where, in the chain of inequalities, we used Cauchy-Schwarz twice. Thus,  we can pick $p_0,p_1$ so that the expectation is achieved and set $P' = \Gamma_2(p_0)\cap \Gamma_2(p_1)$ to be the required set.
\end{proof}
Notice that, by our assumption on $L(p)$ we have that $|I(P',L)| \gg N^{3/2 - 2\eps}$ and so we have not lost anything by replacing $P$ with $P'$. If we draw a picture of the lines through $p_0$ and the lines through $p_1$ we get a skewed grid that contains the large set $P'$. Our next goal is to `straighten-out' this grid so that the lines through $p_0$ are parallel to the X axis and the lines through $p_1$ are parallel to the Y axis. This will be obtained using a projective transformation sending $p_0$ and $p_1$ to the line at infinity.

\subsection{Projective space over $\F$}
Since we are working over a finite field it makes sense to stop for a minute to define the basic properties of projective space. It will help to keep in mind the mental picture of real projective space in which we place the real plane on a slice $z=1$ in three dimensional space and then project points to the half sphere by passing a line through the origin.

More accurately, the construction of $d$-dimensional projective space over $\F$ is as follows. Take the space $\F^{d+1} \setminus \{0\}$ and identify two non zero vectors $x,y \in \F^{d+1}$ if there exists a non zero $\lambda \in \F$ such that $x = \lambda y$. Call the resulting space $\P\F^d$ or the $d$-dimensional projective space. Points in $\P\F^d$ are given using $d+1$ {\em homogenous coordinates} $x = (x_0:\ldots:x_d)$ with each point having exactly $p-1$ different homogenous coordinates. The `regular' or `affine' $d$-dimensional space $\F^d$ can be embedded into $\P\F^d$ by sending $x = (x_1,\ldots,x_d) \in \F^d$ to $x' = (1:x_1:x_2:\ldots:x_d) \in \P\F^d$. Notice that, since the first coordinate is fixed to one, two different vectors map to two different points (the zero vector goes to $(1:0:\ldots:0)$ which is non zero!). Using this embedding, we call the points with homogenous coordinates having $x_0=0$ {\em points at infinity}. The set of all such points is called the hyperplane at infinity and is another projective space of dimension smaller by one. For example, the points at infinity in $\P\F^2$ form a projective line $\P\F^1$ called the line at infinity. 

To get a feeling for these concepts consider the following example. Let $\ell$ be a line in $\F^2$. Suppose $\ell$ is given by the equation $aX + bY +c = 0$ with $a,b,c \in \F$. Now embed $\F^2$ in $\P\F^2$ using three homogenous coordinates $(X:Y:Z)$ so that the points at infinity are those with $Z=0$. A point $(X,Y)$ on $\ell$ will map to $(X:Y:1)$ and will satisfy the equation $aX + bY + cZ = 0$. Notice that homogenous equations do not care about choice of homogenous coordinates and so it makes sense to look at their common solutions in projective space. Thus, we can identify the line $\ell$ with the line $\ell'$ in projective space given by the homogenous equation $aX + bY + cZ = 0$. The affine points (points that are not at infinity) on $\ell'$ are precisely those that come from points in $\ell$. There is however a new point, at infinity, given by $(-b:a:0)$ (or any of its non zero multiples). At least one of $a,b$ are non zero and so this makes sense. Notice that the coordinates of this point correspond to the direction of $\ell$. This means that, if we take another line $\ell_2$ in the same direction of $\ell$ and embed it into $\P\F^2$ it will also contain the same point at infinity! Thus, lines in the same direction intersect at a fixed point at infinity corresponding to their direction. 

One last thing we need to consider are linear mappings over $\P\F^2$. These can be given by any $3 \times 3$ matrix and act on the points of $\P\F^2$ in the obvious way. Notice that such a mapping may take points not at infinity to the line at infinity and vice versa. Also notice that such mappings map lines to lines.

Let us go back to our set of points $P'$ and the lines $L$. We can embed these into $\P\F^2$ and then perform a liner transformation taking $p_0$ to the point $(1:0:0)$ at infinity and $p_1$ to $(0:1:0)$ at infinity. By our previous discussion one can check that, considering the `affine' points (those with $Z=1$) after the transformation, the lines through $p_0$  are now parallel to the $X$ axis and the lines through $p_1$ are parallel to the Y-axis. There might be some points in $P'$ that were moved to infinity but all of those must lie on a single line passing through $p_0$ and $p_1$ and so there are at most $N^{1/2 + 2\eps}$ of those and we can safely ignore them. This means that, after the projective transformation, most of the set $P'$ is in the `affine' part $(z=1)$ and so we can go back to $\F^2$ (discarding the $z=1$ coordinate) and we now have that the set $P'$ is contained in a grid $A \times B$ with $|A|,|B| \leq N^{1/2 + 2 \eps}$

\subsection{Counting incidences on the grid}
Renaming $N$ to be $N^{1+\eps'}$ for some $\eps' >0 $ and using the projective transformation above we see that Theorem~\ref{thm-stfinite} will follow from the following claim:
\begin{claim}[ST over a grid]
Let $P,L$ be  sets of at most $N$ points/lines and suppose $P \subset A \times B$ with $|A|,|B| \leq N^{1/2}$. If $p^\alpha < N < p^{2-\alpha}$ for sufficiently small $\alpha> 0$, then $|I(P,L)| \leq N^{3/2 - \eps}$.
\end{claim} 
\begin{proof}
Our goal will be to reduce to Corollary~\ref{cor-entropy}. Our grid is given by `rows' $b \in B$ and `columns' $a \in A$. For each $b \in B$ let $R(b) = P \cap (A \times \{b\})$ denote the set of points in $P$ that have Y coordinate equal to $b$. Denote also by $H(b) = \{ \ell \in L \,|\, \ell\cap R(b) \neq \emptyset \}$ the set of lines that pass through some point in $R(b)$.  We can ignore the few lines that are parallel to either the X axis or the Y axis.

The first step is to find two rows $b_0$ and $b_1$ such that many lines pass through both $R(b_0)$ and $R(b_1)$. This will be obtained, again, using a probabilistic argument. Notice that in the inequalities below we use the fact that each line can intersect $R(b)$ in at most one  point for each $b \in B$.
\begin{eqnarray*}
\E_{b_0,b_1}\left[ |H(b_0) \cap H(b_1)| \right] &=& \frac{1}{N} \sum_{b_0,b_1 \in B} \sum_{\ell \in L} \sum_{p \in {R(b_0)}} \sum_{q \in R(b_1)} 1_{p \in \ell} \cdot 1_{q \in \ell} \\
&=& \frac{2}{N} \sum_{\ell \in L} \sum_{p,q \in P} 1_{p \in \ell} \cdot 1_{q \in \ell} \\
&=& \frac{2}{N}\sum_{\ell \in L} |P(\ell)|^2 \\
&\geq& \frac{2}{N^2} \left( \sum_{\ell \in L} |P(\ell)| \right)^2 \\
&\geq& N^{1 - c\eps}.
\end{eqnarray*}

Therefore, we can find two elements in $B$, w.l.o.g take these to be $b=0$ and $b=1$ such that $|H(0) \cap H(1) | \geq N^{1 - c\eps}$ for some constant $c > 0$. Let $L' = H(0) \cap H(1)$ be the set of lines that contain both a point with $b=0$ and a point with $b=1$ in $P$. As before, we could have removed all lines  with less than $N^{1/2 - 2\eps}$ points on them and so we can assume $|I(P,L')| \geq N^{3/2 - c\eps}.$

Since at most $O(N)$ incidences can occur on the lines $b=0$ or $b=1$ we have that $$ |\{ (p,\ell) \in P \times L' \,|\, p \in \ell \text{ and } p \not\in R(0) \cup R(1) \}| \gtrsim N^{3/2 - c \eps}. $$

Consider a point $p = (a,b)$ with $b \not\in \{0,1\}$ that lies on a line $\ell \in L'$. This line passes through two points, say $(x_0,0)$ and $(x_1,1)$ with $x_0,x_1 \in A$ and so we have $(a,b) = (1-b)(x_0,0) + b(x_1,1)$ which means that $(1-b)x_0 + bx_1 \in A$. This gives
$$ |\{ (b,x_0,x_1) \in B \times A \times A \,|\, (1-b)x_0 + bx_1 \in A \}| \gtrsim  N^{3/2 - c\eps}.$$ So, there exists a subset $B' \subset B$ with $|B'| > N^{1/2 - 2c\eps} > N^{1/4}$ such that for all $b \in B'$ we have 
$$ |\{ (x_0,x_1) \in A \times A \,|\, (1-b)x_0 + bx_1 \in A \}| \gtrsim  N^{1 - 2c\eps}.$$ Dividing by $b$, we see that this implies that $E\left(A,\frac{b}{1-b}A\right) \leq N^{1/2 + O(\eps)} = |A|^{1+ O(\eps)}$ for all $b \in B'$ (since many sums fall in the small set $A$). This contradicts Corollary~\ref{cor-entropy} if we take $\eps$ small enough. Notice here that we need use the bound $p^\alpha < N < p^{2-\alpha}$, with $\alpha$ taken to be $\eps/C$ for some large constant $C$, to satisfy the conditions of Corollary~\ref{cor-entropy}.
\end{proof}

\ifediting
\bibliographystyle{alpha}
 \bibliography{incidence}

\end{document}
\fi

\newif\ifediting

\ifediting
\documentclass[11pt]{article}

\usepackage{amsmath,amsthm,amssymb}
\newcommand{\remove}[1]{}
\setlength{\topmargin}{0.3in} \setlength{\headheight}{0in}
\setlength{\headsep}{0in} \setlength{\textheight}{8.0in}
\setlength{\topsep}{0.1in} \setlength{\itemsep}{0.0in}
\parskip=0.05in
 \textwidth=6.5in 
\oddsidemargin=0truecm \evensidemargin=0truecm

\newtheorem{thm}{Theorem}[section]
\newtheorem{claim}[thm]{Claim}
\newtheorem{lem}[thm]{Lemma}
\newtheorem{define}[thm]{Definition}
\newtheorem{cor}[thm]{Corollary}
\newtheorem{obs}[thm]{Observation}
\newtheorem{example}[thm]{Example}
\newtheorem{construct}[thm]{Construction}
\newtheorem{conjecture}[thm]{Conjecture}
\newtheorem{THM}{Theorem}
\newtheorem{question}{Question}
\newtheorem{fact}[thm]{Fact}
\newtheorem{prop}[thm]{Proposition}

\def\F{{\mathbb{F}}}
\def\Q{{\mathbb{Q}}}
\def\Z{{\mathbb{Z}}}
\def\N{{\mathbb{N}}}
\def\R{{\mathbb{R}}}
\def\K{{\mathbb{K}}}
\def\C{{\mathbb{C}}}
\def\A{{\mathbb{A}}}
\def\P{{\mathbb{P}}}
\def\cP{{\cal P}}
\def\cS{{\mathcal S}}
\def\cE{{\mathcal E}}
\def\V{{\mathbf{V}}}
\def\I{{\mathbf{I}}}
\def\bx{{\mathbf x}}
\def\by{{\mathbf y}}
\def\E{{\mathbb E}}

\def\half{ \frac{1}{2}}
\newcommand{\ip}[2]{\langle #1,#2 \rangle}
\def\sumN{\sum_{i=1}^n}
\def\_{\,\,\,\,\,}
\def\prob{{\mathbf{Pr}}}
\newcommand{\entropy}[1]{ {\text{H}_{\infty}\left({#1}\right)} }
\def\modulo{\text{mod}}
\def\omm{ \{0,1\} }
\def\id{ \textit{id} }

\def\D{{\partial}}
\def\gap{\textsf{gap}}
\def\sign{\textsf{sign}}
\def\spar{\textsf{sparse}}
\def\span{\textsf{span}}
\def\Part{\textbf{Part}}
\def\Mon{\textbf{Mon}}
\def\sing{\textbf{sing}}
\def\Und{\textsf{Und}}
\def\Comp{\textsf{Comp}}
\def\rank{\textsf{rank}}
\def\poly{\textsf{poly}}
\def\codim{\textsf{codim}}
\def\cp{\textsf{cp}}
\def\uni{\textsf{Uni}}
\def\ext{\textsf{\bf Ext}}
\def\extt{\textsf{\bf Ext2}}

\def\fplus{{\,+_f\,}}

\newcommand{\epsclose}{\stackrel{\epsilon}{\thicksim} }
\newcommand{\eclose}[1]{\stackrel{{#1}}{\thicksim} }
\newcommand{\eps}{\epsilon}
\newcommand{\Anote}[1]{\begin{quote}{\sf Avi's Note:} {\sl{#1}} \end{quote}}

\begin{document}

\title{Incidence Theorems -- Lecture Notes}
\date{}
\maketitle


\fi

\section{Multi-source extractors}\label{sec-multi}

An extractor (short for {\em randomness extractor}) is an algorithm that transforms `weak' sources of randomness, into strong random bits. For example, suppose $X$ is a random variable distributed uniformly over a set $S \subset \{0,1\}^n$  of size $|S| = 2^k$. Informally, $X$ contains $k$ bits of randomness and so we would hope to use $X$ to generate $k$ (or close to $k$) unbiased random bits. We do not know, however, the set $S$ and so have to construct a single function $f: \{0,1\}^n \mapsto \{0,1\}^k$ such that $f(X)$ will be uniform for all such $X$. It is not hard to see that such a function does not exist (even if we require the output to be only one random bit). There are two different ways around this obstacle and both give rise to interesting questions. One is to allow $f$ to use a small number of auxiliary random bits (independent of $X$). Such a function $f$ is called a seeded-extractor. We will talk more about these later when we discuss applications of the finite field Kakeya problem. Another approach is to assume some structure on the source $X$ (say, that the set $S \subset \{0,1\}^n$ belongs to some `nice' family of sets). This restriction allows, in many interesting cases, for a single {\em deterministic extractor} $f$. 

\subsection{Extractors for constant number of sources: BIW}

One well-studied class of deterministic extractors are for sources belonging to the class of several independent blocks. In this family, the random source is partitioned into blocks $X = (X_1,\ldots,X_t) \in \left(\{0,1\}^n \right)^t$ such that the different blocks are independent (as random variables) and each contains some minimal amount of entropy. The right notion of entropy (and the most commonly used) is {\em min-entropy}. The min-entropy of $X$, denoted $\entropy{X}$ is defined as the maximal $k$ such that $\P[X = x] \leq 2^{-k}$ for all $x$ in the support of $X$. If $X$ is, as above, uniform on a set of size $2^k$ (these are called flat sources) then it has min-entropy $k$. Conversely, one can show that every source with min-entropy $k$ is a convex combination of flat sources of min-entropy $k$ \cite{CG88}. Thus, it is enough to argue about flat sources. A deterministic $(k,\eps)$-extractor for $t$-sources is a function $$f :  \left(\{0,1\}^n \right)^t \mapsto \{0,1\}^m$$ such that for every $t$ independent random variables $X_1, \ldots,X_t \in \{0,1\}^n$, each of min-entropy at least $k$, the output $f(X_1,\ldots,X_t)$ is $\eps$-close to the uniform distribution in statistical distance\footnote{The statistical distance between two distributions is the $\ell_1$ distance of their probability vectors.}. Clearly $m$ is at most $t \cdot k$ and the goal is in general to output as many bits as possible with error $(\eps)$ as small as possible. In what follows, we will mostly talk about extractors with one bit of output since this is usually the hard case (once you have one bit you can usually get more).

A simple probabilistic argument shows that there are deterministic extractors (even for two sources) that works for min-entropy $k \sim \log(n)$. However, explicit constructions are pretty hard to find. Today, the best constructions all use in some way or another tools from incidences in finite fields (or, equivalently the sum product theorem). We will sketch the first construction to make use of these tools. This is a result of Barak, Impagliazzo and Wigderson \cite{BIW06} and was the first explicit extractor for a constant number of sources that worked for any linear min-entropy $k = \Omega(n)$. The idea is as follow: consider three independent sources $X,Y,Z \in \{0,1\}^n$ all with min-entropy $\delta \cdot n$. Suppose $n$ is prime, and identify the three variables with elements in the finite field $\F = GF(2^n)$ that does not contain subfields (there are some subtleties to discuss if $n$ is not prime but we will not go there). Let $W = X + YZ$ be computed over $\F$. Using the Szemeredi-Trotter theorem we can show that $W$ is close to having min-entropy at least $(\delta + \eps)n$ for some small positive $\eps$. We will prove this below with entropy replaced by set-size. Once we know this, we can iterate this construction and take the functions:
$$ f_1(X_1,X_2,X_3) = X_1 + X_2X_3 $$
$$ f_2(X_1,\ldots,X_9) = (X_1 + X_2 X_3) + (X_4 + X_5 X_6)\cdot (X_7 + X_8 X_9)$$
etc.. and prove by induction that after a constant number of steps (depending on $\delta,\eps$) we will get a distribution that is close to uniform (the last step that goes from high min-entropy to close-to-uniform requires a slightly different argument). Let us now prove a set-size variant of the lemma at the heart of this argument (the generalization to min-entropy is left as an exercise). The proof will work over any field $\F$ in which the Szemeredi-Trotter type bound $I(P,L) \lesssim N^{3/2 - \eps}$ holds. As was mentioned before, even though we proved this bound over finite prime fields, the proof holds over any field that does not contain large subfields (which is relevant w.r.t the construction described above).

\begin{lem}
Let $A,B,C \subset \F$ be subsets of size $|\F|^{\alpha} < N < |\F|^{1-\alpha}$ of a field $\F$ in which the Szemeredi-Trotter bound holds. Then $|A + BC| > N^{1+ \eps}$, with $\eps>0$ depending only on $\alpha$.
\end{lem}
\begin{proof}
Let $$S(x) = |\{ (a,b,c) \in A \times B \times C \,|\, a + bc = x \}|$$ denote the `weight' of $x$ in the distribution $A + BC$ (i.e., when we sample three independent samples from $A,B,C$ and compute $a+bc$). We have
\begin{equation}\label{eq-mu1}
	\sum_x S(x) = N^3.
\end{equation}
On the other hand, if we assume in contradiction that  $|A+BC| \leq N^{1+\eps}$, we get
\begin{equation}\label{eq-mu2}
	\sum_x S(x)^2 \geq \frac{\left( \sum_x S(x )\right)^2}{|A+BC|} \geq N^{5 - \eps}
\end{equation}
Now, if we define $T = \{ x \,|\, S(x) > N^{2 - 2\eps} \}$ and using the two inequalities above we get
$$ N^{1-2\eps} \leq |T| \leq N^{1+2\eps}.$$
This implies
\begin{equation}\label{eq-inc}
	|\{ (a,b,c,x) \in A \times B \times C \times T \,|\, a+bc = x \}| \geq N^{3 - 4\eps}.
\end{equation}
This can be viewed as a bound on line/point incidences by defining a set of points $P = C \times T$ and a set of lines $L = \{ \ell_{a,b} \}$ with $\ell_{a,b}$ defined by the equation $a + Xb = Y$ for all $a\in A, b \in B$. The number of lines/points is at most $N^{2+2\eps}$ and the number of incidences is $N^{3 - 4\eps}$. If $\eps$ is sufficiently small this will contradict Szemeredi-Trotter.
\end{proof}

\subsection{Bourgain's two source extractor}

When the number of sources is two (the smallest possible) much less is known. Suppose we want a two-source extractor for min-entropy $k$ that outputs a single bit (with some fixed small $\eps$). A probabilistic argument shows that this can be done with $k \sim \log(n)$. A simple explicit construction exists when $k > n/2$ (take the inner product modulo two) \cite{CG88}. For a long time this was the best known explicit construction. This was changed a few years back when Bourgain \cite{Bour2source} showed how to use the ST theorem to construct an extractor for two sources of min-entropy $k = (1/2 - \eps)n$ for some positive $\eps$. It is an open problem to give an explicit construction of a two-source  extractor for min-entropy significantly less than $n/2$. There are construction of weaker objects called {\em dispersers} for two sources that only output a bit that is non constant (i.e., a bit that is equal to both zero and one with some positive probability).   These constructions work for min-entropy as low as $k = n^{o(1)}$ \cite{BRSW,BKSSW}. These constructions use a whole lot of tools, among which are those that we have developed here. We will now show Bourgain's construction with one bit output (it is possible to extract more bits).
 
We start with the analysis of the  inner product extractor, which works for min-entropy larger than $n/2$. Recall that it is enough to consider two `flat' sources $A,B \subset \{0,1\}^n$. To bound the distance of $\ip{A}{B}$ from the uniform distribution on one bit it is enough to bound the following quantity which we will refer to as the bias
$$ \bias(A,B) = \left| \frac{1}{|A||B|} \sum_{a \in A} \sum_{b \in B} (-1)^{\ip{a}{b}}\right|.$$

To bound the bias we use some Cauchy-Schwarz calculations:
\begin{eqnarray*}
\bias(A,B) &\leq& \frac{1}{|A||B|} \sum_{a\in A} \left| \sum_{b \in B} (-1)^{\ip{a}{b}} \right| \\
&\leq& \frac{1}{|A|^{1/2}|B|} \left( \sum_{a \in A} \left| \sum_{b \in B} (-1)^{\ip{a}{b}}\right|^2\right)^{1/2} \\
&\leq& \frac{1}{|A|^{1/2}|B|} \left( \sum_{a \in \{0,1\}^n} \left| \sum_{b \in B} (-1)^{\ip{a}{b}}\right|^2\right)^{1/2} \\
&=& \frac{1}{|A|^{1/2}|B|} \left( \sum_{a \in \{0,1\}^n}  \sum_{b,b' \in B} (-1)^{\ip{a}{b-b'}} \right)^{1/2} \\
&=& \frac{1}{|A|^{1/2}|B|} \left(  \sum_{b \in B} 2^n\right)^{1/2} \\
&=& \left( \frac{2^n}{|A||B|}\right)^{1/2}.
\end{eqnarray*}
The bias is equal to the difference between the probability that $\ip{A}{B}$ is equal to one and the probability that it is equal to zero. Thus, if $|A||B| \geq C \cdot 2^n$ then the distance of $\ip{A}{B}$ from the uniform distribution will be roughly $1/\sqrt{C}$. This shows that the inner product function is a $(k,\eps)$ extractor for $k \gg \frac{n}{2} + \log(1/\eps)$. 

But how can the above calculation be useful if we wish to handle smaller entropy? Clearly, the inner product function is not enough since we can take $A$ and $B$ to be orthogonal subspaces of dimension $n/2$ each. Can we fix our construction to avoid such bad examples?

The first step is to observe that, in the calculation above, one can replace the set size of $A,B$ with a more refined quantity. For a distribution $\mu$ on some finite set $\Omega$ (i.e., $\mu$ is a function from $\Omega$ to $\R_{\geq 0}$ with sum of values equal one) we will denote the  $\ell_2$-{\em energy} of $\mu$ by $$E(\mu) = \left( \sum_x \mu(x)^2 \right)^{-1}.$$ Notice that if $\mu$ is a uniform distribution on some subset $A$ then $E(\mu) = |A|$. Notice also that our old notations for additive energy $E(A,B) = \frac{|A|^2|B|^2}{|Q(A,B)|}$ satisfies $E(A,B) = E(\mu_{A+B})$, where $\mu_{A+B}$ is the distribution obtained by sampling two independent variables $a \in A$ and $b \in B$ at uniform and then outputting $a+b$. Another interpretation of $E(\mu)$ is as the inverse of the `collision probability' $cp(\mu) = \sum_x \mu(x)^2$ which is the probability of two independent copies of $\mu$ being equal to each other. We can similarly define bias for distributions as
$$ \bias(\mu_1,\mu_2) = \left| \E_{x_1 \sim \mu_1, x_2 \sim \mu_2}\left[ (-1)^{\ip{x_1}{x_2}}\right] \right|. $$

It is straightforward to verify that the calculation above also works if we replace set size with energy. That is:
$$ \bias(\mu_1,\mu_2)  \leq \left( \frac{2^n}{E(\mu_1)E(\mu_2)}\right)^{1/2}.$$ What is more surprising is that the bound is not changed by much if we replace the two distributions $\mu_1,\mu_2$ with the distributions of sums of {\em several independent copies} drawn from the same distributions. To see this observe that
\begin{eqnarray*}
\bias(\mu_1,\mu_2)^2 &=&
\left| \E_{x_1 \sim \mu_1, x_2 \sim \mu_2}\left[ (-1)^{\ip{x_1}{x_2}}\right] \right|^2 \\
 &\leq& \left( \E_{x_1 \sim \mu_1} \left| \E_{ x_2 \sim \mu_2}\left[ (-1)^{\ip{x_1}{x_2}}\right] \right| \right)^2 \\
&\leq& \E_{x_1 \sim \mu_1} \left| \E_{ x_2 \sim \mu_2}\left[ (-1)^{\ip{x_1}{x_2}}\right] \right|^2 \\
&=&   \E_{x_1 \sim \mu_1, x_2,x_3 \sim \mu_2}\left[ (-1)^{\ip{x_1}{x_2+x_3}}\right] \\
&=& \bias(\mu_1, \mu_2 \text{`+'} \mu_2),
\end{eqnarray*}
where the ad hoc notation $\mu_2 \text{`+'} \mu_2$ means summing two independent copies drawn from $\mu_2$ (this is actually the convolution $\mu_2 * \mu_2$). Iterating this calculation four times we can obtain, for example, the following claim
\begin{claim}\label{cla-foursum}
Let $A,B \subset \{0,1\}^n$. Then
$$ \bias(A,B) \leq \bias(4 \cdot A, 4 \cdot B) ^{1/16},$$ where $4\cdot A$ denotes the distribution of sums of four independent uniform variables from $A$ (similarly for $B$).
\end{claim}

Bourgain's approach to constructing a two-source extractor for minentorpy rate $1/2 - \eps$ is as follows: Construct a set $S \subset \{0,1\}^n$ such that for all subsets $A \subset S$ with $|A| > |S|^{1/2 - \eps}$ we have $E(4 \cdot A) \gg 2^{n/2}$. Then define the extractor $f : S \times S \mapsto \{0,1\}$ as $f(x,y) = \ip{x}{y}$. Formally, we will need to identify $S$ with some $\{0,1\}^{n'}$ but this will not be a problem. This will work since, if we take two (flat) sources $A,B \subset S$ of size $|S|^{1/2 - \eps}$ (this corresponds to min-entropy rate $> 1/2 - \eps$) then the bias of their inner product is bounded by $ \left( 2^n / E(4\cdot A) E(4 \cdot B)\right)^{1/32}$ which will be close to zero since $$ E(4\cdot A) E(4 \cdot B)  \gg 2^n $$ (the power of $1/32$ really doesn't change much).

Due to some technical difficulties in working over fields of characteristic two, we will construct the set $S$ over the group $\Z_3^n$ instead of over $\Z_2^n$. To justify this `switch' observe that we can replace $(-1)$ in the summations above with a complex root of unity of order 3, say $\omega = \exp(2\pi i/3)$ and define
$$ \bias_\omega(\mu_1,\mu_2) = \left| \E_{x_1 \sim \mu_1, x_2 \sim \mu_2}\left[ \omega^{\ip{x_1}{x_2}}\right] \right|,$$ where the inner product is over $\Z_3$. Then, the same calculation as above gives
$$ \bias_\omega(\mu_1,\mu_2)  \leq \left( \frac{3^n}{E(\mu_1)E(\mu_2)}\right)^{1/2}$$ as well as
$$ \bias_\omega(A,B) \leq \bias_\omega(4 \cdot A, 4 \cdot B) ^{1/16}.$$ This means that, if we can construct a set $S$ such that every subset $A$ of size $|S|^{1/2-\eps}$ satisfies $E(4 \cdot A) \gg 3^{(1/2+\eps)n}$ we will get that, for all roots of unity of order $3$ the bias $\bias_\omega(A,B)$ is close to zero for all sets $A,B$ in $S$ of size $> |S|^{1/2-\eps}$. It is not hard to show then that the distribution of $\ip{a}{b}$, with $a \in A,b \in B$ is close to the uniform distribution on three elements (so the output of the extractor is not a bit but rather a uniform element in a set of size three).

The construction of $S \subset \Z_3^n$ is as follows. Suppose $n = 2p$, where $p$ is a prime number (there are ways to handle other values of $n$ but this is a technicality). Let $\F$ be a finite field of size $3^p$ so that $\F$ does not have large subfields (i.e., we can use the ST theorem in $\F^2$). Identify $\Z_3^n$ with $\F^2$ by writing each element of $\F$ in some basis of $\F$ over $GF(3)$. So addition in $\F^2$ is the same as coordinate wise addition modulo 3 in $\Z_3^n$. We can now define:
$$ S  = \{ (x,x^2) \,|\, x \in \F \} \subset \F^2 \sim \Z_3^n. $$

We proceed with the analysis. Let $\tilde A \subset S$ be of size $|\tilde A| > |S|^{1/2 - \eps} = 3^{p(1/2 -\eps)}$. Then there is a subset $A \subset \F$ of the same size such that $\tilde A = \{ (a,a^2) \,|\, a \in A \}$. We need to show that $E(4 \cdot \tilde A) = E(4 \cdot A) \geq 3^{p(1+\eps)} \gg 3^{n/2}$. For this purpose, define  for all $x,y \in \F$ the set $$ R_{x,y} = \{ (a_1,a_2,a_3,a_4) \in A^4 \,|\, \sum a_i = x, \,\, \sum a_i^2 = y \}.$$ Notice that $$ E(4 \cdot A) = \frac{|A|^8}{\sum |R_{x,y}|^2}. $$ Thus, if we could show that $$ R = \sum_{x,y} |R_{x,y}|^2 \leq |A|^{6 - 8\eps}$$ we would have that $E(4 \cdot A) \geq |A|^{2 + 8\eps} \geq 3^{p(1+ \eps)},$ for sufficiently small $\eps$, as required.

 We will bound the sum $R$ by partitioning it into two parts. Let $c >0$ be a constant to be chosen later. Define  $T_1 = \{ (x,y) \,|\, |R_{x,y}| \leq |A|^{2 - c\eps} \}$ and $T_2 = \{ (x,y) \,|\, |R_{x,y}| > |A|^{2 - c\eps} \}$. Then $R = R_1 + R_2$, where $R_1$ is the sum of $|R_{x,y}|^2$ over $(x,y) \in T_1$ and $R_2$ is the sum over the `large' terms in $T_2$ (the rest of the terms). $R_1$ is easy to bound since the total number of terms in $T_1$ is at most $|\F|^2 \leq |A|^{2 + 8\eps}$ and each term is at most $|A|^{2 - c\eps}$ and so the total bound is $$ R_1 \leq |A|^{2 + 8 \eps}|A|^{4 - 2c \eps} \ll |A|^{6 - 8\eps} $$ if $c$ is sufficiently large.

 To bound $R_2$ we will bound the size of the set $T_2$ by $|A|^{2 - 8\eps}$. If we can do that than we will be done since we can combine this bound with the trivial bound of $|A|^2$ on each of the $|R_{x,y}|$'s to obtain $R_2 \leq |A|^{6 - 8\eps}$ (to see the trivial bound of $|A|^2$ notice that fixing $a_1,a_2$ allows us to solve for $a_3,a_4$). Suppose in contradiction that $|T_2| > |A|^{2 - 8\eps}$. By definition, for each $(x,y) \in T_2$ there are at least $|A|^{2-c\eps}$ solutions $(a_1,a_2,a_3,a_4) \in A^4$ to the equations $$ a_1 + a_2 + a_3 + a_4 = x, $$ $$ a_1^2 + a_2^2 + a_3^2 + a_4^2 = y.$$ Let $T_3 = \{ (x,(x^2 - y)/2) \,|\, (x,y) \in T_2 \}$ (this is where we need the characteristic to be different than two!) so that $|T_3| = |T_2|$ and such that for each $(x,y) \in T_3$ we have at least $|A|^{2 - c\eps}$ solutions $(a_1,a_2,a_3,a_4) \in A$ to the equations 
 $$ a_1 + a_2 + a_3 + a_4 = x, $$ $$ a_1a_2 + a_1a_3 + a_1a_4 + a_2a_3 + a_2a_4 + a_3a_4 = y.$$ We can now eliminate $a_4$ so that for all $(x,y) \in T_3$ we have at least $|A|^{2 - c\eps}$ solutions $(a_1,a_2,a_3) \in A^3$
to the single equation $$ y = a_1a_2 + a_2a_3 +a_3a_1 - (a_1+a_2+a_3)^2 + (a_1 + a_2 + a_3)\cdot x. $$ Thus we have the bound
\begin{eqnarray*}
	 |\{ (x,y,a_1,a_2,a_3) \in T_3 \times A^3 \,|\, y = a_1a_2 + a_2a_3 +a_3a_1 -\\ (a_1+a_2+a_3)^2 + (a_1 + a_2 + a_3)\cdot x \}| \geq |A|^{4 - (8+c)\eps}.
\end{eqnarray*}
 We can now fix $a_3$ to some value $b \in A$ so that
\begin{eqnarray*}
	|\{ (x,y,a_1,a_2) \in T_3 \times A^2 \,|\, y = a_1a_2 + a_2b +ba_1 -\\ (a_1+a_2+b)^2 + (a_1 + a_2 + b)\cdot x \}| \geq |A|^{3 - (8+c)\eps}.
\end{eqnarray*}
This last quantity can be viewed as the set of incidences of the lines $\ell_{a_1,a_2}, (a_1,a_2) \in A^2$ define as $\ell_{a_1,a_2} = \{ (u,v) | v =  a_1a_2 + a_2b +ba_1 - (a_1+a_2+b)^2 + (a_1 + a_2 + b)\cdot u \}$ and the set of points $T_3$. The number of lines is clearly at most $|A|^2$ and the number of points $|T_3|$ is at most $|A|^{2 +  c\eps}$ since we have $\sum_{x,t} |R_{x,y}| = |A|^4$ and so there can only be at most $|A|^{2 + c\eps}$ summands larger than $|A|^{2 - c\eps}$. Taking $\eps$ to be small enough we will get a contradiction to the ST theorem since the number of points/lines is roughly $|A|^2$ and the number of incidences approaches $|A|^{3}$. This completes the proof.

\ifediting
\bibliographystyle{alpha}
 \bibliography{incidence}

\end{document}
\fi

\chapter{Kakeya sets}\label{chap-3}

\newif\ifediting

\ifediting
\documentclass[11pt]{article}

\usepackage{amsmath,amsthm,amssymb}
\newcommand{\remove}[1]{}
\setlength{\topmargin}{0.3in} \setlength{\headheight}{0in}
\setlength{\headsep}{0in} \setlength{\textheight}{8.0in}
\setlength{\topsep}{0.1in} \setlength{\itemsep}{0.0in}
\parskip=0.05in
 \textwidth=6.5in 
\oddsidemargin=0truecm \evensidemargin=0truecm

\newtheorem{thm}{Theorem}[section]
\newtheorem{claim}[thm]{Claim}
\newtheorem{lem}[thm]{Lemma}
\newtheorem{define}[thm]{Definition}
\newtheorem{cor}[thm]{Corollary}
\newtheorem{obs}[thm]{Observation}
\newtheorem{example}[thm]{Example}
\newtheorem{construct}[thm]{Construction}
\newtheorem{conjecture}[thm]{Conjecture}
\newtheorem{THM}{Theorem}
\newtheorem{question}{Question}
\newtheorem{fact}[thm]{Fact}
\newtheorem{prop}[thm]{Proposition}

\def\F{{\mathbb{F}}}
\def\Q{{\mathbb{Q}}}
\def\Z{{\mathbb{Z}}}
\def\N{{\mathbb{N}}}
\def\R{{\mathbb{R}}}
\def\K{{\mathbb{K}}}
\def\C{{\mathbb{C}}}
\def\A{{\mathbb{A}}}
\def\P{{\mathbb{P}}}
\def\cP{{\cal P}}
\def\cS{{\mathcal S}}
\def\cE{{\mathcal E}}
\def\V{{\mathbf{V}}}
\def\I{{\mathbf{I}}}
\def\bx{{\mathbf x}}
\def\by{{\mathbf y}}
\def\E{{\mathbb E}}

\def\half{ \frac{1}{2}}
\newcommand{\ip}[2]{\langle #1,#2 \rangle}
\def\sumN{\sum_{i=1}^n}
\def\_{\,\,\,\,\,}
\def\prob{{\mathbf{Pr}}}
\newcommand{\entropy}[1]{ {\text{H}_{\infty}\left({#1}\right)} }
\def\modulo{\text{mod}}
\def\omm{ \{0,1\} }
\def\id{ \textit{id} }

\def\D{{\partial}}
\def\gap{\textsf{gap}}
\def\sign{\textsf{sign}}
\def\spar{\textsf{sparse}}
\def\span{\textsf{span}}
\def\Part{\textbf{Part}}
\def\Mon{\textbf{Mon}}
\def\sing{\textbf{sing}}
\def\Und{\textsf{Und}}
\def\Comp{\textsf{Comp}}
\def\rank{\textsf{rank}}
\def\poly{\textsf{poly}}
\def\codim{\textsf{codim}}
\def\cp{\textsf{cp}}
\def\uni{\textsf{Uni}}
\def\ext{\textsf{\bf Ext}}
\def\extt{\textsf{\bf Ext2}}

\def\fplus{{\,+_f\,}}

\newcommand{\epsclose}{\stackrel{\epsilon}{\thicksim} }
\newcommand{\eclose}[1]{\stackrel{{#1}}{\thicksim} }
\newcommand{\eps}{\epsilon}
\newcommand{\Anote}[1]{\begin{quote}{\sf Avi's Note:} {\sl{#1}} \end{quote}}

\begin{document}

\title{Incidence Theorems -- Lecture Notes}
\date{}
\maketitle


\fi

\section{Kakeya sets in $\R^n$}\label{sec-kakeyareal}

The Kakeya problem in $\R^n$ deals with the most efficient way to `pack' many tubes  ($\eps$-neighborhoods of line segments) that point in different directions. As we shall see,  this question reduces  to a discrete question about incidences of line segments pointing in sufficiently separated directions. The starting point is the definition of a Kakeya set.

\begin{define}[Kakeya Set]
A compact set $K \subset \R^n$ is a {\em Kakeya} set if it contains a unit line segment in each direction. More formally, for every $x \in  S^{n-1}$ there exists $y = f(x) \in K$
 such that $\{ y + tx \,|\, t \in [0,1] \} \subset K$.	
\end{define}
  
It is known \cite{Bes28} that Kakeya sets can have measure zero (we will not prove this here). A more refined question has to do with the minimal {\em dimension} of a Kakeya set. For simplicity, we will use the Minkowski dimension (also known as covering/box dimension) but other notions (in particular Hausdorff dimension) are often studied in the literature.  The Minkowski dimension (which we will refer to simply  as  `dimension' from now on) is defined as follows. Let $B_\eps(K)$ denote the minimal number of balls of radius $\eps$ needed to cover the (bounded) set $K \subset \R^n$. The dimension of $K$ is defined as $$ \dim(K) = \limsup_{\eps \rightarrow 0} \frac{\log B_\eps(K)}{\log (1/\eps)}. $$ (technically, this is the {\em upper} Minkowski dimension). Roughly speaking, if $\dim(K) \leq d$ then $K$ can be covered by $\sim (1/\eps)^d$ balls of radius $\eps$, where the $\sim$ notation hides constants that might depend on the dimension $n$. It is a good exercise at this point to verify that this definition of dimension agrees with the usual definition of dimension for subspaces (intersected with a unit ball) and, more generally, algebraic surfaces. For example, the dimension of a line segment of length $L$ is 1 since it can be covered by $L(1/\eps)$ balls of radius $\eps$ and the factor of $L$ disappears in the limit. Also, a set of positive measure in $\R^n$ must have dimension $n$.

The Kakeya conjecture (sometimes called the Euclidean Kakeya conjecture) states that Kakeya sets $K \subset \R^n$ must have  dimension $n$, which is the highest possible. This conjecture is open for $n \geq 3$ (we will prove the $n=2$ case below) and is related to several important questions in analysis, PDE's and number theory. We refer the reader to the excellent survey \cite{Tao01} for more on these applications/connections.

For this section only, we will  think of $n$ as constant and use our asymptotic notations $\sim,\gtrsim,\lesssim$ to suppress constant depending on $n$ (these will disappear in the limit when $\eps \rightarrow 0$). Thus, one can replace the quantity $B_\eps(K)$ with a slightly more convenient quantity having to do with the number of grid points close to $K$.  More formally, let $G_\eps = \eps \Z^n$ denote the $\eps$-grid in $\R^n$. Notice that every point in $\R^n$ is at distance at most $\sqrt{n} \eps$ from some grid point. Let $G_\eps(K)$ denote the number of points in $G_\eps$ that are at distance at most $10\sqrt{n}\cdot\eps$ (the constant $10$ is arbitrary and is there just to give some wiggle room). Thus, in our notations, $G_\eps(K) \sim B_\eps(K)$ and so we can use $G_\eps(K)$ from now on. We will sometimes abuse notations and treat $G_\eps(K)$ as the {\em set} of points of distance at most $10\sqrt{n}\cdot\eps$ from $K$. 
 
Before moving on to the  discretized setting, mentioned above, we will prove the $n=2$ case of the Kakeya conjecture.
\begin{thm}[Davies \cite{Davies}]
Let $K \subset \R^2$ be a Kakeya set. Then $\dim(K) = 2$.	
\end{thm}
\begin{proof}
Let $K'$ be the $\eps$-neighborhood of $K$ (i.e., all points of distance at most $\eps$ from $K$) and notice that $G_\eps(K) \sim G_\eps(K')$. We will show $G_\eps(K') \gtrsim \frac{1}{\eps^2} \cdot \frac{1}{\log (1/\eps)}$, which will prove the theorem. 

Consider $\sim 1/\eps$ tubes of width $\eps$ with one endpoint at the origin and with the other endpoints spread along the  first quadrant part of the unit circle in $\R^2$. That is, take $\ell_j$ to be the line segment connecting the origin with $(\cos(\eps j \pi/2),\sin(\eps j \pi/2)$ and take $T_j $ to be its $\eps$-neighborhood. Since $K$
 is a Kakeya set, we can `shift' each of the tubes $T_j$ (without changing its direction) so that they are contained in $K'$. Suppose we have already done that and that for all $j$, $T_j \subset K$. 

Notice that, for each $j$ we have $G_\eps(T_j) \gtrsim 1/\eps$ and that, for $i \neq j$ we have $G_\eps(T_i \cap T_j) \lesssim \frac{1}{\eps|i-j|}$ (this is why we took tubes only in the first quadrant). Using Cauchy-Schwarz we get
\begin{eqnarray*}
\frac{1}{\eps^2} &\sim& \sum_j G_\eps(T_j) 
= \sum_{x \in G_\eps(K')} \sum_{j} 1_{x \in T_j} \\
&\lesssim& \left( G_{\eps}(K')\right)^{1/2} \cdot \left( \sum_{x \in G_\eps(K')} \left( \sum_j 1_{x \in T_j}\right)^2\right)^{1/2} \\
&\sim& \left( G_{\eps}(K')\right)^{1/2} \cdot\left( \sum_{i,j} G_\eps(T_i \cap T_j)\right)^{1/2} \\
&\lesssim& \left( G_{\eps}(K')\right)^{1/2} \cdot \left( \frac{1}{\eps^2} + \frac{1}{\eps}\sum_{i \neq j} \frac{1}{|i-j|}\right)^{1/2} \\
&\lesssim& \left( G_{\eps}(K')\right)^{1/2} \cdot \left( \frac{\log (1/\eps)}{\eps^2}\right)^{1/2}.
\end{eqnarray*}
Rearranging we get the bound 
$$ G_\eps(K') \gtrsim \frac{1}{\eps^2 \log(1/\eps)} $$ which gives $\dim(K)=2$, when $\eps$ goes to zero in the definition of dimension.
\end{proof}

\subsection{The $n/2$ bound}

We will now see a proof that gives a lower bound of $n/2$ on the dimension of Kakeya sets in $\R^n$. This will also set up some of the notations for the next part which will use additive combinatorics to get a better bound of the form $(4/7)n$. 

Similarly to the set of tubes $T_j$ used above we will now need an $\eps$-separated set of directions $\Omega \subset S^{n-1}$. Since we are ignoring constants depending on $n$ we can easily find such a set with $|\Omega| \sim (1/\eps)^{n-1}$. Thus, if $K$ is a Kakeya set we have that for all $w \in \Omega$ there exists $a_w \in \R^n$ such that the segment $\ell_w = \{ a_w + tw \,|\, t \in [0,1] \} \subset K$. Let us denote by $b_w = a_w + w$ the second `endpoint' of the line segment in direction $w$. For each $w \in \Omega$ let $a_w',b_w'$ be the grid points (in $G_\eps$) closest to $a_w,b_w$. Consider the line segment $\ell_w'$ connecting $a_w'$ to $b_w'$. Since $\ell_w'$ is obtained from $\ell_w$ by moving its endpoints by at most $\lesssim \eps$ we have that the set of directions $\Omega'$ of line segments $\ell_w'$ has size at least $\gtrsim |\Omega| \gtrsim (1/\eps)^{n-1}$. Let $A = \{a_w'\,|\, w' \in \Omega\}$ and $B = \{b_w' \,|\, w' \in \Omega\}$. Then both $|A|,|B|$ are at most $G_\eps(K)$ (since $a_w,b_w$ are in $K$ and $a_w',b_w'$ are the closest  grid points to them). On the other hand, we have $ |B - A| \geq |\Omega'|$ (since the differences between $a_w'$ and $b_w'$ cover all directions in $\Omega'$). Since $|B-A| \leq |A||B|$ we have $ |A||B| \gtrsim (1/\eps)^{n-1}$ which implies $G_\eps(K) \gtrsim (1/\eps)^{(n-1)/2}$. This means that $\dim(K) \geq (n-1)/2$. 

To go from $(n-1)/2$ we use a tensoring argument: Observe that, if $K$ is a Kakeya set in $\R^n$ then, for all $t \in \N$, $K^t \subset \R^{nt}$ is also a Kakeya set. It is also simple to verify that $\dim(K^t) = t\dim(K)$ and so, using our previous bound on $K^t$ we get that $t\dim(K) \geq (nt-1)/2$. Dividing by $t$ and taking $t$ to infinity we get that $\dim(K) \geq n/2$.

In \cite{Wolff99}, Wolff proved an even stronger bound of $(n+2)/2$ for general $n$. We will not see this proof here and focus on later developments, starting with the work of Bourgain \cite{Bou99}, giving (increasingly higher) bounds of the form $\alpha n$ for $\alpha > 1/2$. These results use ideas and tools from additive combinatorics. 

\subsection{Additive combinatorics methods}

The proof of the $n/2$ bound we saw above uses only the `endpoints' of the line segments (after moving them slightly so that they are on a grid). Not using other information cannot go beyond $n/2$ as there are sets of points $|A|,|B| $ on the grid with $|B-A| \sim  |A||B|$ and this is all we used in the proof. To go beyond this barrier we will also use, as a starting point, the {\em mid points} of the segments. That is, the points $(a_w + b_w)/2$. We saw that we can shift $a_w,b_w$ by at most $\lesssim \eps$ so that they are on the grid $G_\eps$. It is easy to see that a similar shifting argument can also put all three points $a_w',b_w'$ and the mid point $c_w' = (a_w' + b_w')/2$ on the grid $G_\eps$. The distance we need to shift the endpoints will grow by at most a constant factor which we do not care about. Define as before the sets $A,B$ to contain the points $a_w',b_w'$ with $w' \in \Omega'$ (the new set of directions we obtain after the shifting). As before we have $|\Omega'| \gtrsim (1/\eps)^{n-1}.$ We can also place all points in $A,B$ in an $O(\eps)$-neighborhood of $G_\eps(K)$ and so we have $|A|,|B| \lesssim G_\eps(K)$. Let us denote $N = G_\eps(K)$

We still know that $|A-B| \geq |\Omega'| \gtrsim (1/\eps)^{n-1}$ is large. Now we can also incorporate the midpoints to claim that, in some sense, the sumset $A+B$ is small! To see this consider all sums of the form $(a_w' + b_w')/2$ with $w' \in \Omega'$. These sums will all fall in a set of size $\lesssim N = G_\eps(K)$ and so, if we assume the dimension of $K$ is at most $d$, this set will be of size at most $\lesssim (1/\eps)^d$. Thus, the intuition is that, since the difference set is large, the sumset cannot be too small and so we will get a contradiction if the dimension of $K$ is smaller than some bound. There is, however, a serious difficulty. The sumset $\{a_w' + b_w' \,|\, w' \in \Omega' \}$ whose size we want to argue about (we can discard the $1/2$) is not really the sumset $A+B$ but rather a sub-sumset determined by some fixed family of pairs (indexed by $\Omega'$). To see the way around it we remind ourselves of the Balog-Szemeredi-Gowers theorem which says that, if a dense family of pairs in a sumset $A+B$ is small, then there are large subsets $A',B'$ with small sumsets. Thus, if the family of pairs $(a_w',b_w')$ with $w \in \Omega'$ is dense, in the sense that $|\Omega'| \geq (|A||B|)^{1-\delta}$ (for some small constant $\delta$) then we can hope to save the situation in some way. This is indeed the case if we are shooting for a $(1/2 + \delta')n$ type of bound on $\dim(K)$. To see this, notice that $|A|,|B| \leq N$ and, if we assume in contradiction that $N \ll (1/\eps)^{(1/2 + \delta)n}$ we get that $|\Omega'| \sim (1/\eps)^{n-1} \sim N^{2 - \delta}$ as required.

In \cite{Bou99} Bourgain carries out the above plan with the aid of a modified version of the BSG theorem tailored for this situation.  Bourgain's proof was simplified considerably by Katz and Tao \cite{KT99,KT02} who also improved  the constant $\alpha$ from Bourgain's original $13/25$ to almost $0.596..$. We will see below the simplified proof which gives $4/7$. Let us now state the general reduction from the Kakeya dimension question to a simple to state question in additive combinatorics. This more general reduction will allow us to use more points (not just the midpoints) which will be useful in simplifying the proofs. 

\begin{define}[$\SD(R,\beta)$]
Let $A,B \subset H$ be finite subsets of an abelian group $H$ with $|A|,|B| \leq N$. Let $\Gamma \subset A \times B$. Let $R \subset \N$ and suppose that for all $r \in R$ we have $|\{ a + rb \,|\, (a,b) \in \Gamma\}| \leq N$. We say that the statement $\SD(R,\beta)$ holds over $H$ if for every pair of sets $A,B$ as above, we have $|\{ a - b \,|\, (a,b) \in \Gamma \}| \leq N^{\beta}$.
\end{define}
 
\begin{lem}
Suppose $\SD(R,\beta)$ holds over $\R^n$ for $R = \{1,2,\ldots,r\}$ and $\beta > 1$. Then for all Kakeya sets $K \subset \R^n$ we have $\dim(K) \geq n/\beta$.
\end{lem}
\begin{proof}
Let $K \subset \R^n$ be a Kakeya set. We will treat $r$ as a constant (as $\eps$ will go to zero). Consider, as before an $\eps$-seperated set of directions $\Omega \subset S^{n-1}$ of size $\sim (1/\eps)^{n-1}$ and let $a_w,b_w \in K$ be the endpoints of a unit line segment in direction $w \in \Omega$ that is contained in $K$. Fix $\eps > 0$ to be sufficiently small and let $N = G_\eps(K)$. We can move each pair $(a_w,b_w)$ by at most $O(\eps)$ to new points $(a_w',b_w')$ on the grid $G_\eps$ so that all combinations $a_w' + jb_w'$ for all $j \in R$ fall in a set of size $\sim G_\eps(K)$ (similarly to what we did for sums). Since the line segments $\ell_w$ were moved by $O(\eps)$ we have that  the new set of directions $w' = a_w' - b_w'$, denoted $\Omega'$, is also of size at least $\gtrsim (1/\eps)^{n-1}$. Therefore, for all $j \in R$ we have $$ |\{ a_w' + jb_w' \,|\, w' \in \Omega' \} | \lesssim N.$$ Using the $\SD(R,\beta)$ assumption we get that $$ (1/\eps)^{n-1} \lesssim|\{ a_w' - b_w' \,|\, w' \in \Omega' \}| \lesssim N^{\beta}$$ which gives the required bound (after a tensoring argument). 
\end{proof}

Thus, in order to prove the Kakeya conjecture it suffices to show that $\SD(R,1)$ holds for some fixed set $R \subset \N$. 

\subsection{The $4n/7$ bound}

We will now prove that $\SD(\{1,2\},7/4)$ holds over any abelian group $H$  which will give a $4n/7$ bound on $\dim(K)$ (we will assume that the order of $1$ in $H$ is larger than 2). Since we want to bound the size of the set $\{ a - b \,|\, (a,b) \in \Gamma \}$ we may assume w.l.o.g that the difference $a-b$ is distinct for each $(a,b) \in \Gamma$ (since removing edges will only decrease the bound $N$ on the other sets in the definition). We are thus interested in bounding the number of edges in $\Gamma$ or $|\Gamma|$.

The main ingredient in the proof is the notion of a `gadget' which we now define. A `gadget' will be a substructure in the graph $\Gamma$ with certain restriction on linear combinations on edges. More formally, a gadget $G$ is defined as a 4-tuple $G = (V_A,V_B,E,C)$ with
\begin{itemize}
	\item $V_A = (a_1,\ldots,a_s)$, $V_B = (b_1,\ldots,b_\ell)$ two sets of formal variables.
	\item $E$ a subset of $V_A \times V_B$ (we call these `edges').
	\item $C$ a set of constraints of the form $a_i + rb_j = a_{i'} + r'b_{j'}$ with $i,j \in [s]$, $i',j' \in [r]$ and $r,r'$ integers.
\end{itemize}

An example of a simple gadget is $G_1 = (V_A,V_B,E,C)$ with:
\begin{eqnarray*}
 V_A = \{a_1,a_2\}, V_B = \{b_1,b_2\},\\ E = \{(a_1,b_1), (a_2,b_2) \}, C = \{ a_1 + 2b_1 = a_2 + 2b_2\}.
\end{eqnarray*}
We say that a gadget $G$ {\em appears} in the graph $\Gamma = A\times B$ (with $A,B$ subsets of the abelian group $H$) if we can map $V_A,V_B$ to subsets of $A,B$ such that the set of edges $E$ is contained in the set of edges induced by $\Gamma$ and the constraints in $C$ are satisfied. For example, if we take $A = B = \{1,2,3,4,5\}$ and $\Gamma = A \times B$ then the gadget $G_1$ above appears in $\Gamma$ by taking $a_1=1,a_2 = 3,b_1=5,b_2=4$ (since $1+2\cdot5 = 3 + 2\cdot 4$).

We can also count the number of times a gadget appears in $\Gamma$ in the obvious way as the number of different ways to map $V_A,V_B$ into subsets of $A,B$ so that the edges/constraints are satisfied. For example, we will show that, if we take $\Gamma \subset A \times B$ such that for all edges $(a,b) \in \Gamma$ we have $a + 2b \in H' \subset H$ then $G_1$ will appear in $\Gamma$ at least $|\Gamma|^2/|H'|$ times. This fact follows from a Cauchy-Schwarz calculation that is given by the following lemma.
\begin{lem}\label{lem-cs}
Let $W$ be a finite set and let $f : W \mapsto Z$ be a mapping to some other finite set $Z$. Then
$$ \left| \{ (v,u) \in W^2 \,|\, f(v)=f(u) \}\right| \geq |W|^2/|Z|. $$
\end{lem}
\begin{proof}
The size of the set is given by the sum
\begin{eqnarray*}
\sum_{u,v \in W}1_{f(u) = f(v)} &=& \sum_{u,v \in W} \sum_{z \in Z} 1_{f(u)=z} \cdot 1_{f(v)=z} \\
&=& \sum_{z \in Z} \left( \sum_{u \in W} 1_{f(u) = z}\right)^2 \\
&\geq& \frac{1}{|Z|} \cdot \left( \sum_{z \in Z} \sum_{u \in W} 1_{f(u)=z} \right)^2 \\
&=& \frac{|W|^2}{|Z|}.
\end{eqnarray*}
\end{proof}
 
We can apply this lemma to the gadget $G_1$ as follows. Take the function $f: \Gamma \mapsto H'$ to be $f(a,b) = a+2b$. We then get the promised bound $|\Gamma|^2/|H'|$.

 Recall that in the definition of $\SD(R,\beta)$ we have a bound $N$ on the sizes  $|A|,|B|$ as well as on the sizes of each of the sets $\{ a + rb \,|\, (a,b) \in \Gamma\}$. Once we have a gadget $G$ and we can give a lower bound on the number of times it appears in $\Gamma$ the next step is to give a corresponding upper bound on the number of times $G$ appears in $\Gamma$ in terms of $N$. This will be done by showing that we can `encode' each gadget using a few elements, each in a set of size at most $N$. For example, the gadget $G_1$ can be encoded as $(a_1,a_2,a_1+2b_1)$ since, from this triple, we can recover both $b_1$ and $b_2$ (using the fact that $a_1+2b_1 = a_2 + 2b_2$). Since all three elements in this triple are in a set of size at most $N$ we get that there can be at most $N^3$ appearances of $G_1$ in $\Gamma$. Combining this with the lower bound obtained from Lemma~\ref{lem-cs} we get $ |\Gamma|^2/N \leq N^3$ or $|\Gamma| \leq N^2$. This bound on $\Gamma$ is not very interesting and we proved it just to give an idea of the proof technique. To prove the claimed $4/7$ bound we need to get a bound of $|\Gamma| \leq N^{7/4}$ which will require a more elaborate gadget.

\subsubsection*{A more elaborate gadget}

Consider the gadget $G_{4/7}$ given by 
$$ V_A = \{a_1,a_2\}$$
$$ V_B = \{b_1,b_2,b_3\}$$
$$ E = \{ (a_1,b_1),(a_1,b_2),(a_2,b_2),(a_2,b_3) \} $$
$$ C = \{ a_1 + 2b_1 = a_2 + 2b_3 \}.$$

Let $\Gamma \subset A \times B$ be as in the definition of $\SD(\{1,2\},\beta)$ so that $|A|,|B| \leq N$ and $|\{ a + rb \,|\, (a,b) \in \Gamma\}| \leq N$ for $r=1,2$. 
Our first step is to give a lower bound on the number of appearances of $G_{4/7}$ in $\Gamma$. For this purpose consider the set $$M = \{((a,b),(a',b')) \in \Gamma^2 \,|\, a = a' \}$$ (the set of paths of length two). Using Lemma~\ref{lem-cs} we have $|M| \geq |\Gamma|^2/N$. Let $f : M \mapsto H^3$ be defined as $$f((a,b),(a',b')) = (b',a+2b). $$  Notice that each collision of $f$ gives an appearance of the gadget $G_{4/7}$. Since the image of $f$ is contained in a set of size $N^2$,  Lemma~\ref{lem-cs} gives at least $$|M|^2/N^2 \geq |\Gamma|^4/N^4$$ collisions/appearances of $G_{4/7}$. 

We now give an upper bound using the `encoding' argument. Here we will use the fact, mentioned above, that w.l.o.g the differences on the edges of $\Gamma$ are distinct. This will be useful since, knowing the different $a-b$ on some edge identifies this edge and so also identifies its two endpoints. 

Let $G' = (a_1,a_2,b_1,b_2,b_3)$ be an appearance of $G_{4/7}$ (so that all edges/constraints are satisfied). We will show that $G'$ can be recovered from the triple $(b_3,a_1+b_2,a_1+b_1)$. Since we have the bound $|\{ a + b \,|\, (a,b) \in \Gamma\}| \leq N$ we know that there are at most $N^3$ such triples which will give the same upper bound on the number of appearances of $G_{4/7}$. We now describe the decoding. The first step is to decode $$a_2 - b_1 = (a_1 + b_1) - 2b_3$$ (using the constraint $a_1 + 2b_1 = a_2 + 2b_3$). Then we  compute $$ b_1 - b_2 = (a_1 + b_1) - (a_1 + b_2)$$. Using these two we can compute $$ a_2 - b_2 = (a_2 - b_1) + (b_1 - b_2).$$ Now, using the distinctness of differences we can recover $a_2,b_2$ (since it is an edge in $\Gamma$) and from them the rest of the vertices in $G'$. Putting the two bounds together gives the required $|\Gamma| \leq N^{7/4}$.

\ifediting
\bibliographystyle{alpha}
 \bibliography{incidence}

\end{document}
\fi

\newif\ifediting

\ifediting
\documentclass[11pt]{article}

\usepackage{amsmath,amsthm,amssymb}
\newcommand{\remove}[1]{}
\setlength{\topmargin}{0.3in} \setlength{\headheight}{0in}
\setlength{\headsep}{0in} \setlength{\textheight}{8.0in}
\setlength{\topsep}{0.1in} \setlength{\itemsep}{0.0in}
\parskip=0.05in
 \textwidth=6.5in 
\oddsidemargin=0truecm \evensidemargin=0truecm

\newtheorem{thm}{Theorem}[section]
\newtheorem{claim}[thm]{Claim}
\newtheorem{lem}[thm]{Lemma}
\newtheorem{define}[thm]{Definition}
\newtheorem{cor}[thm]{Corollary}
\newtheorem{obs}[thm]{Observation}
\newtheorem{example}[thm]{Example}
\newtheorem{construct}[thm]{Construction}
\newtheorem{conjecture}[thm]{Conjecture}
\newtheorem{THM}{Theorem}
\newtheorem{question}{Question}
\newtheorem{fact}[thm]{Fact}
\newtheorem{prop}[thm]{Proposition}

\def\F{{\mathbb{F}}}
\def\Q{{\mathbb{Q}}}
\def\Z{{\mathbb{Z}}}
\def\N{{\mathbb{N}}}
\def\R{{\mathbb{R}}}
\def\K{{\mathbb{K}}}
\def\C{{\mathbb{C}}}
\def\A{{\mathbb{A}}}
\def\P{{\mathbb{P}}}
\def\cP{{\cal P}}
\def\cS{{\mathcal S}}
\def\cE{{\mathcal E}}
\def\V{{\mathbf{V}}}
\def\I{{\mathbf{I}}}
\def\bx{{\mathbf x}}
\def\by{{\mathbf y}}
\def\E{{\mathbb E}}

\def\half{ \frac{1}{2}}
\newcommand{\ip}[2]{\langle #1,#2 \rangle}
\def\sumN{\sum_{i=1}^n}
\def\_{\,\,\,\,\,}
\def\prob{{\mathbf{Pr}}}
\newcommand{\entropy}[1]{ {\text{H}_{\infty}\left({#1}\right)} }
\def\modulo{\text{mod}}
\def\omm{ \{0,1\} }
\def\id{ \textit{id} }

\def\D{{\partial}}
\def\gap{\textsf{gap}}
\def\sign{\textsf{sign}}
\def\spar{\textsf{sparse}}
\def\span{\textsf{span}}
\def\Part{\textbf{Part}}
\def\Mon{\textbf{Mon}}
\def\sing{\textbf{sing}}
\def\Und{\textsf{Und}}
\def\Comp{\textsf{Comp}}
\def\rank{\textsf{rank}}
\def\poly{\textsf{poly}}
\def\codim{\textsf{codim}}
\def\cp{\textsf{cp}}
\def\uni{\textsf{Uni}}
\def\ext{\textsf{\bf Ext}}
\def\extt{\textsf{\bf Ext2}}

\def\fplus{{\,+_f\,}}

\newcommand{\epsclose}{\stackrel{\epsilon}{\thicksim} }
\newcommand{\eclose}[1]{\stackrel{{#1}}{\thicksim} }
\newcommand{\eps}{\epsilon}
\newcommand{\Anote}[1]{\begin{quote}{\sf Avi's Note:} {\sl{#1}} \end{quote}}

\begin{document}

\title{Incidence Theorems -- Lecture Notes}
\date{}
\maketitle


\fi

\section{Kakeya sets in finite fields}\label{sec-kakeyafin}

In his influential survey on the Kakeya problem, Wolff \cite{Wolff99} defined the finite field analog of the problem. Below,  $\F$ will denote a finite field of size $q$ (not necessarily prime).
\begin{define}
 A Kakeya set $K \subset \F^n$ is a set containing a line in every direction. More formally, for all $x \in \F^n$ there exists $y \in \F^n$ such that $\{ y + tx \,|\, t \in \F \} \subset K$.
\end{define}

Wolff asked whether a bound of the form $|K| \geq C_n \cdot q^n $ holds for all Kakeya sets $K$, with $C_n$ a constant depending only on $n$. Here, one should think of $n$ as fixed and the field size $q$ goes to infinity (thinking of $q \sim 1/\eps$ helps). The proofs we saw in the previous section, using additive combinatorics, can be carried out also over finite fields. For example, using the $\SD(\{1,2\},7/4)$ statement (over the abelian group $\F^n$) one gets a bound of $|K| \geq C_n \cdot q^{4n/7}$. In \cite{Dvir08}, the polynomial method was used to give an answer to Wolff's question. Initially, a proof of $C_n q^{n-1}$ was shown and then, using an observation of Alon and Tao, the tight exponent $C_n q^n$ was also obtained. We will see both the original proof and the improvement, which is another nice example of the usefulness of working in projective space.

\subsection{Proof of the finite field Kakeya conjecture}

To start, we define Nikodym sets which are closely related to Kakeya sets.
\begin{define}
A Nikodym set $K \subset \F^n$ is a set for which, through every point not in $K$ there is a line that intersects $K$ in all points but one. More formally, if for all $y \not\in K$ there exists $x$ such that $\{y + tx \,|\, t \in \F^* \} \subset K$.
\end{define}
This definition seems `stronger' than the Kakeya definition. However, the two definitions are related by a factor of $q$.
\begin{claim}
If there exists a Kakeya set $K$ of size $T$ in $\F^n$ then there exists a Nikodym set $M$ in $\F^n$ of size at most $qT$. In fact, one can take $M = \{ tx \,|\, t \in \F, x \in K \}$.
\end{claim}
\begin{proof}
For each $x \in \F^n$, there is $y \in \F^n$ such that $\{ y + tx \,|\, t \in \F\} \subset K$. This means that $\{ sy + stx \,|\, s,t \in \F\} \subset M$. Fixing $t = 1/s$ and going over all $s \neq 0$ we get $\{ sy + x \,|\, s \in \F^* \} \subset M$ and so $M$ is a Nikodym set.
\end{proof}

What can we do with a small Nikodym set $M$? Suppose we have a polynomial $f(x_1,\ldots,x_n)$  of degree $d$ and we know the values of $f$ on all points of $M$. If the degree of $f$ is less than $q-1$ we could use these values to recover the values of $f$ everywhere! To see this, suppose we wish to find the value of $f$ at a point $x \not\in M$. Let $y \in \F^n$ be such that the punctured line $\ell = \{ x + ty \,|\, t \in \F^* \}$ is contained in $M$. Restricting $f$ to the line $\ell$ we get a polynomial of degree at most $d$ and we know its values in $q-1 > d$ points. Therefore, we can recover the coefficients of the restricted polynomial and compute its value at the missing point $x$. But this means that the number of points in $M$ must be larger than the number of coefficients in a degree $q-2$ polynomial. Since, otherwise, we could find a non zero polynomial of degree $q-2$ that vanishes everywhere in $M$ and is not identically zero. This will be a contradiction since, using the above decoding procedure, we would get that the polynomial is zero everywhere.

The last step of this argument requires proving that a multivariate polynomial that is not identically zero, has a non zero value at some point in $\F^n$. This is known as the Schwartz-Zippel Lemma:
\begin{lem}
Let $f \in \F[x_1,\ldots,x_n]$ be a non zero polynomial of degree $d$. Then there are at most $dq^{n-1}$ points in $\F^n$ where $f$ vanishes.
\end{lem} 
\begin{proof}
By induction on $n$. The $n=1$ case is the fundamental theorem of algebra. For larger $n$, write $f(x_1,\ldots,x_n) = \sum_{j=1}^{r}g_j(x_2,\ldots,x_n)x_1^j$, such that w.l.o.g $g_r(x_2,\ldots,x_n)$ is non zero of degree $d-r$. By induction, there are at most $(d-r)q^{n-2}$ zeros of $g_r$. For each one of them, the restricted polynomial $f$ (which is now a polynomial in the single variable $x_1$) might vanish identically and so have $q$ zeros. For the rest of the assignments (at most $q^{n-1}$) to $g_r$, $f$ will remain a non zero univariate polynomial of degree $r$ and so can have at most $r$ zeros. Combining, we get at most $(d-r)q^{n-1} + rq^{n-1} = dq^{n-1}$ zeros for $f$.
\end{proof}

We can now give the proof of the $q^{n-1}$ bound on Kakeya sets.
\begin{thm}
For every Kakeya set $K \subset \F^n$ we have $|K| > (1/n!) q^{n-1}$.
\end{thm}
\begin{proof}
Let $M$ be a Nikodym set of size $q|K|$. If $|K| \ll (1/n!) q^{n-1}$ then $ |M| \ll (1/n!)q^n$ and we can find a polynomial $f$ of degree $d \leq q-2$ that vanishes on $M$ and is not identically zero. For every point $x \not\in M$ consider the restriction of $f$ to the line $\ell$ passing through $x$ that has $q-1$ points in $M$. The restriction of $f$ to this line is a degree $d$ polynomial and so, since $d < q-1$ we get that $f$ must vanish everywhere, contradicting the Schwartz-Zippel lemma.
\end{proof}

Using a tensoring argument, as we saw for Kakeya sets over the reals, one can amplify this bound to $C_{n,\eps}q^{n-\eps}$ for all $\eps > 0$. There is, however, a clever way to get rid of this $\eps$ completely. This has to do with working over projective space. 

Recall that the $n$ dimensional projective space $\P\F^n$ is defined formally as the set of $n+1$ dimensional non zero vectors with two vectors identified if they are a constant multiple of each other. We embed the affine space $\F^n$ in $\P\F^n$ by adding a coordinate $x_0 = 1$ so that the points at infinity are given by the hyperplane $x_0=0$. Recall that a line in $\F^n$ in direction $y$ will hit the point at infinity with coordinates $(0,y_1,\ldots,y_n)$ (since multiplying by a constant doesn't change the point the choice of $y$ is also up to a constant). There is a way to extend the polynomial method to work over the projective space. In projective space, we only consider homogeneous polynomials, those in which every monomial has the same degree. The set of zeros of a homogeneous polynomial is well defined in $\P\F^n$ since $f(ax_1,\ldots,ax_n) = a^d f(x_1,\ldots,x_n)$ for all non zero $a \in \F$. When we embed $\F^n$ into $\P\F^n$ in the above described manner, we can accompany this with an embedding of $\F[x_1,\ldots,x_n]$ into the set of homogeneous polynomials in variables $x_0,x_1,\ldots,x_n$. This is done by sending $f(x_1,\ldots,x_n)$ of degree $d$ into $$f^h(x_0,x_1,\ldots,x_n) = x_0^df(x_1/x_0,\ldots,x_n/x_0)$$ or, in other words, multiplying each monomial of $f$ of degree $d-r$ with $x_0^r$ so that the resulting polynomial is homogeneous of the same degree of $f$. Notice that, setting $x_0=1$ in $f^h$ we get $f$ back (thus, $f^h$ is consistent with the embedding of points in $\F^n$ into $\P\F^n$). Also notice that, setting $x_0 = 0$ in $f^h$ we get back the homogeneous part of $f$ of highest degree. This is the restriction of $f^h$ to the hyperplane at infinity. 

Suppose $K \subset \F^n$ is a Kakeya set and embed $\F^n $ into $\P\F^n$ using $x_0=1$. Let $K'$ be the embedding of $K$ (which has the same size as $K$). Saying that $K$ contains a line in every direction is the same as saying that, through each point at infinity $(0,y_1,\ldots,y_n)$ there is a line that has $q$ points in $K'$. Suppose now we had a polynomial $f$ of degree $d \leq q-1$ that vanished on $K$ and consider $f^h$ as above. Using the restrictions to all these lines, we get that $f^h$ must vanish at all points at infinity. This means that the homogeneous part of highest degree of $f$ (which is the same as $f^h(0,x_1,\ldots,x_n)$) vanishes identically. But this is a contradiction since we assumed that $f$ is non zero and so it must have a homogeneous part of highest degree which is non zero.

Using the above argument we get 
\begin{thm}
	For all Kakeya sets $K \subset \F^n$ we have $|K| \geq (1/n!)q^n$
\end{thm}
  
In the finite field setting we also might care about the constant in front of the $q^n$ (this doesn't appear in the real case since we are taking a limit). There is a better bound of $|K| \geq (1/2^n)q^n$ on Kakeya sets proved in \cite{DKSS09} which uses a more sophisticated polynomial argument with zeros of high multiplicities.
  
\subsection{A construction of small Kakeya sets}\label{sec-upper}
We now turn to describing the smallest known Kakeya sets which are
of size
$$|K| \leq \frac{q^n}{2^{n-1}} + O(q^{n-1}),$$ which is, asymptotically as
$q$ tends to infinity, to within a factor of 2 of the lower bound obtained in \cite{DKSS09}. The construction for the case $n=2$ was
given by \cite{MT04} and the generalization for larger $n$ was
observed by the author for odd characteristic and by \cite{SaSu}
for even characteristic. We give here the construction for odd
characteristic.

We will only worry about lines in directions $b =
(b_1,\ldots,b_n)$ with $b_n=1$. The rest of the lines can be added
using an additional $q^{n-1}$ points, which is swallowed by the
low order term. Our set is defined as follows:
\[ K = \left. \left\{ \left(v_1^2/4 + v_1 \cdot t, \ldots,
v_{n-1}^2/4 + v_{n-1} \cdot t, t\right) \, \right| \,
v_1,\ldots,v_{n-1},t \in \F \right \}. \] Let $b =
(b_1,\ldots,b_{n-1},1)$ be some direction. Then $K$ clearly
contains the line in direction $b$ through the point
$(b_1^2/4,\ldots,b_{n-1}^2/4,0)$. We now turn to showing that $|K|
\leq \frac{q^n}{2^{n-1}}$. Notice that the sum of the first
coordinate of $K$ and the square of the last one is equal to
\[ v_1^2/4 + v_1 \cdot t + t^2 = (v_1/2 + t)^2 \]
and so is a square in $\F$. Since $\F$ has odd characteristic it
contains at most  $\approx q/2$ squares. Let $x_1,\ldots,x_n$
denote the coordinates of the set $K$. Fixing the last coordinate
we get that the first coordinate $x_1$ can take at most $\approx
q/2$ values. The same holds for $x_2,\ldots,x_{n-1}$ and so we get
a bound of $\approx \frac{q^n}{2^{n-1}}$ on the size of $K$.

\ifediting
\bibliographystyle{alpha}
 \bibliography{incidence}

\end{document}
\fi

\newif\ifediting

\ifediting
\documentclass[11pt]{article}

\usepackage{amsmath,amsthm,amssymb}
\newcommand{\remove}[1]{}
\setlength{\topmargin}{0.3in} \setlength{\headheight}{0in}
\setlength{\headsep}{0in} \setlength{\textheight}{8.0in}
\setlength{\topsep}{0.1in} \setlength{\itemsep}{0.0in}
\parskip=0.05in
 \textwidth=6.5in 
\oddsidemargin=0truecm \evensidemargin=0truecm

\newtheorem{thm}{Theorem}[section]
\newtheorem{claim}[thm]{Claim}
\newtheorem{lem}[thm]{Lemma}
\newtheorem{define}[thm]{Definition}
\newtheorem{cor}[thm]{Corollary}
\newtheorem{obs}[thm]{Observation}
\newtheorem{example}[thm]{Example}
\newtheorem{construct}[thm]{Construction}
\newtheorem{conjecture}[thm]{Conjecture}
\newtheorem{THM}{Theorem}
\newtheorem{question}{Question}
\newtheorem{fact}[thm]{Fact}
\newtheorem{prop}[thm]{Proposition}

\def\F{{\mathbb{F}}}
\def\Q{{\mathbb{Q}}}
\def\Z{{\mathbb{Z}}}
\def\N{{\mathbb{N}}}
\def\R{{\mathbb{R}}}
\def\K{{\mathbb{K}}}
\def\C{{\mathbb{C}}}
\def\A{{\mathbb{A}}}
\def\P{{\mathbb{P}}}
\def\cP{{\cal P}}
\def\cS{{\mathcal S}}
\def\cE{{\mathcal E}}
\def\V{{\mathbf{V}}}
\def\I{{\mathbf{I}}}
\def\bx{{\mathbf x}}
\def\by{{\mathbf y}}
\def\E{{\mathbb E}}

\def\half{ \frac{1}{2}}
\newcommand{\ip}[2]{\langle #1,#2 \rangle}
\def\sumN{\sum_{i=1}^n}
\def\_{\,\,\,\,\,}
\def\prob{{\mathbf{Pr}}}
\newcommand{\entropy}[1]{ {\text{H}_{\infty}\left({#1}\right)} }
\def\modulo{\text{mod}}
\def\omm{ \{0,1\} }
\def\id{ \textit{id} }

\def\D{{\partial}}
\def\gap{\textsf{gap}}
\def\sign{\textsf{sign}}
\def\spar{\textsf{sparse}}
\def\span{\textsf{span}}
\def\Part{\textbf{Part}}
\def\Mon{\textbf{Mon}}
\def\sing{\textbf{sing}}
\def\Und{\textsf{Und}}
\def\Comp{\textsf{Comp}}
\def\rank{\textsf{rank}}
\def\poly{\textsf{poly}}
\def\codim{\textsf{codim}}
\def\cp{\textsf{cp}}
\def\uni{\textsf{Uni}}
\def\ext{\textsf{\bf Ext}}
\def\extt{\textsf{\bf Ext2}}

\def\fplus{{\,+_f\,}}

\newcommand{\epsclose}{\stackrel{\epsilon}{\thicksim} }
\newcommand{\eclose}[1]{\stackrel{{#1}}{\thicksim} }
\newcommand{\eps}{\epsilon}
\newcommand{\Anote}[1]{\begin{quote}{\sf Avi's Note:} {\sl{#1}} \end{quote}}

\begin{document}

\title{Incidence Theorems -- Lecture Notes}
\date{}
\maketitle


\fi

\section{Randomness Mergers from Kakeya sets.}\label{sec-mergers}

In CS, the interest in the finite field Kakeya problem originated in the work of Lu, Reingold, Vadhan and Wigderson \cite{LRVW03}. Motivated by extractor constructions, the following question was raised: Suppose $X_1,\ldots,X_k$ are random variables each distributed over $\F^n$, where $\F$ is a finite field of order $q$. We do not assume that the $X_i$'s are independent and are guaranteed only that one of them is uniformly distributed over $\F^n$. The question is, what can we say about the entropy of  a random linear combination of $X_1,\ldots,X_k$? To make things simpler, suppose we only have two variables $X,Y \in \F^n$ such that $X$ is uniform on $\F^n$ and  $Y$ could depend on $X$. Let $Z = aX + bY$, where $a,b \in \F$ are both chosen uniformly at random and independently of $X,Y$ and of each other. How `random' is $Z$? 

The connection between this question and the finite field Kakeya problem is as follows: Suppose we had a small Kakeya set $K \subset \F^n$ and take $M = \{ ax \,|\, x \in K, a \in \F \}$ to be the corresponding Nikodym set (see previous section) that is of comparable size to $|K|$. We know that for each $x \in \F^n$ there exists $y = y(x) \in \F^n$ such that $\{ y(x) + tx \,|\, t \in \F\} \subset K$. This means that $\{ stx + sy(x) \,|\, s,t \in \F \} \subset M$. Renaming $st=a$ and $s=b$ we get that $\{ ax + by(x) \,|\, a \in \F, b \in \F^*\} \subset M$. What this means is that, given $X$, one could set $Y = Y(X)$ such that all linear combinations $aX + bY$ with $b$ non zero hit the small set $M$. This means that the output will land in $M$ with high probability (at least $1 - 1/q$) which would imply that $Z = aX + bY$ has low entropy (e.g when using min-entropy). Thus, to answer the question of \cite{LRVW03} we must (in the least) solve the finite field Kakeya conjecture! This problem is even more challenging, since it involves entropy and randomness (which we still need to define properly). Luckily, the polynomial method is sufficiently robust to handle even this harder scenario.

We start with some definitions.  The {\sf statistical distance} between two
distributions $P$ and $Q$ on a finite domain $\Omega$ is defined
as
\[\mathop{\max}_{S\subseteq \Omega} \left|P(S) - Q(S)\right|
.\] We say that $P$ is $\eps$-{\sf close} to $Q$ if the
statistical distance between $P$ and $Q$ is at most $\eps$. The {\sf min-entropy} of a random variable $X$ is
defined as
\[  \entropy{X} \triangleq \min_{x \in \text{supp}(X)}
\log\left(\frac{1}{\prob[X=x]}\right) \] (all logarithms are taken
to the base 2). Intuitively, having min-entropy at least $k$ means having at least $k$ bits of entropy. We say that a random variable $X$ is $\eps$-close to having min-entropy $k$ if there exists another random variable $X'$ such that $X'$ has min-entropy $\geq k$ and $X$ is $\eps$-close to $X'$.

Notice that a r.v $X$ distributed over $\F^n$ can have min-entropy between zero and $n\log(q)$. If $X$ has min-entropy $\beta n \log(q)$ we call $\beta$ the min-entropy {\em rate} of $X$. The following lemma is very useful and allows us to move from min-entropy to set size:

\begin{lem}\label{lem-guv}
Say $X$ is distributed over a finite set $\Omega$ and $X$ is not $\eps$-close to having min-entropy at least $k$. Then there exists a set $T \subset \Omega$ with $|T| \leq 2^k$ such that $\Pr[X \in T] \geq \eps$. 
\end{lem}
\begin{proof}
Take $T = \{ a \in \Omega \, |\, \Pr[X = a] \geq 2^{-k} \}.$ Clearly, $|T| \leq 2^k$ since the sum of probabilities $\Pr[X = a]$ cannot exceed one. If $\Pr[X \in T] < \eps$ we could change the distribution of $X$ slightly by moving the probability mass from $T$ to other values so that the resulting r.v $X'$ will have min-entropy $\geq k$ and will be $\eps$-close to $X$.
\end{proof}

We start by analyzing the case of two random variables:
\begin{thm}
Let $X,Y$ be  two (not necessarily independent) random variables distributed over $\F^n$ and suppose one of them is uniformly distributed. Let $a,b \in \F$ be chosen independently at random and let $Z = aX + bY$. Let $\alpha >0$ be any real number such that $q > n^{10/\alpha}$. Then  $Z$ is $\eps$-close to having min entropy rate $1-\alpha$ with $\eps = q^{-\alpha/10}$.
\end{thm}
\begin{proof}
By symmetry we may assume w.l.o.g that $X$ is uniform.
If $Z$ is not $\eps$-close to having min-entropy rate $1-\alpha$ then, by Lemma~\ref{lem-guv} there is a set $T \subset \F^n$ of size $|T| \leq q^{(1-\alpha)n}$ such that $\Pr[Z \in T] \geq \eps$. Using the polynomial method, we will find a non zero polynomial $f \in \F[x_1,\ldots,x_n]$ of low degree that vanishes on $T$. Let $d$ be the required degree. We need $d$ to satisfy $$ {n+d \choose d} > q^{1-\alpha}n. $$ Using the inequality ${n+d \choose d} \geq (d/n)^n$ and the bound $q > n^{10/\alpha}$ we see that it is enough to take $d = q^{1-\alpha/5}$. 

For each $x \in \F^n$ let $$ p_x = \Pr[ Z \in T \,|\, X = x].$$ and let $$ G = \{ x \in \F^n \,|\, p_x \geq \eps/2 \}. $$ Since $\Pr[Z \in T] \geq \eps$ we have that $\Pr[X \in G] \geq \eps/2$ (this follows from a simple averaging argument). Since $X$ is uniform this implies $|G| \geq (\eps/2)q^n$. We will now show that $f$ vanishes on all points in $G$. 

Fix some $x \in G$. We know that $$ \Pr[aX + bY \in T\,|\, X = x] \geq \eps/2. $$ Thus, we can fix $Y = y$ to some specific value so that the same inequality still holds. That is, there is some $y \in F^n$ such that $$ \Pr[ax + by \in T] \geq \eps/2.$$ Notice that in the last probability the randomness is only over the choice of $a,b$ and that $x,y$ are both fixed. Let $g(a,b) = f(ax + by)$ be the restriction of $f$ to the plane spanned by $x,y$. By the above calculation we get that $g$ has at least $(\eps/2)q^2$ zeros. We know that $g$ can have at most $dq$ zeros (see Schwartz-Zippel lemma from the previous section) and so, if $d > (\eps/2)q$ (which holds in our choice of parameters) we would get that $g(a,b)$ is identically zero. Thus, we have that $g(1,0)=f(x)=0$ and so we conclude that $f$ vanishes on all of $G$.

Now, since $f$ can have at most $dq^{n-1}$ zeros (by Schwartz-Zippel) we get $$ (\eps/2)q^n \leq |G| \leq dq^{n-1} $$ which is a contradiction for the choice of $\eps$ given in the theorem. This concludes the proof.
\end{proof}
 
Looking at things more broadly, a procedure such as the one described above is called a {\em merger}. Mergers allow us to combine several (dependent) random variables, one of which is uniform, into a single variable that has high min-entropy. Mergers are allowed to use a short random `seed' (given above by $a,b \in \F$) and one can show that without this seed the task is impossible. Above we analyzed a simple merger for two sources. Mergers for many sources are important in constructions of seeded-extractors which are procedures that can extract randomness from arbitrary distributions of low min-entropy and that use an additional short random seed. One can generalize the construction above to work with many sources (taking independent coefficients $a_1,\ldots,a_k$ and outputting $\sum_{i}a_iX_i$). This is problematic, however, since the length of the seed grows linearly with the number of sources. One can however, pick the coefficients in a correlated way and get a merger with shorter seed. This is done by passing a curve of degree $k$ through the $k$ points $X_1,\ldots,X_k$ and outputting a random point on this curve. The analysis given above, using the polynomial method, generalizes to this setting as well (see \cite{DW08,DKSS09}).

\ifediting
\bibliographystyle{alpha}
 \bibliography{incidence}

\end{document}
\fi

\chapter{Sylvester-Gallai type problems}\label{chap-4}

\newif\ifediting

\ifediting
\documentclass[11pt]{article}

\usepackage{amsmath,amsthm,amssymb}
\newcommand{\remove}[1]{}
\setlength{\topmargin}{0.3in} \setlength{\headheight}{0in}
\setlength{\headsep}{0in} \setlength{\textheight}{8.0in}
\setlength{\topsep}{0.1in} \setlength{\itemsep}{0.0in}
\parskip=0.05in
 \textwidth=6.5in 
\oddsidemargin=0truecm \evensidemargin=0truecm

\newtheorem{thm}{Theorem}[section]
\newtheorem{claim}[thm]{Claim}
\newtheorem{lem}[thm]{Lemma}
\newtheorem{define}[thm]{Definition}
\newtheorem{cor}[thm]{Corollary}
\newtheorem{obs}[thm]{Observation}
\newtheorem{example}[thm]{Example}
\newtheorem{construct}[thm]{Construction}
\newtheorem{conjecture}[thm]{Conjecture}
\newtheorem{THM}{Theorem}
\newtheorem{question}{Question}
\newtheorem{fact}[thm]{Fact}
\newtheorem{prop}[thm]{Proposition}

\def\F{{\mathbb{F}}}
\def\Q{{\mathbb{Q}}}
\def\Z{{\mathbb{Z}}}
\def\N{{\mathbb{N}}}
\def\R{{\mathbb{R}}}
\def\K{{\mathbb{K}}}
\def\C{{\mathbb{C}}}
\def\A{{\mathbb{A}}}
\def\P{{\mathbb{P}}}
\def\cP{{\cal P}}
\def\cR{{\cal R}}
\def\L{{\cal L}}
\def\B{{\mathcal B}}
\def\S{{\mathbf S}}
\def\cE{{\mathcal E}}
\def\V{{\mathbf{V}}}
\def\I{{\mathbf{I}}}
\def\bx{{\mathbf x}}
\def\by{{\mathbf y}}
\def\E{{\mathbb E}}
\def\SD{{\mathbf{SD}}}
\def\SC{{\mathbf{SC}}}

\def\half{ \frac{1}{2}}
\newcommand{\ip}[2]{\langle #1,#2 \rangle}
\def\sumN{\sum_{i=1}^n}
\def\_{\,\,\,\,\,}
\def\prob{{\mathbf{Pr}}}
\newcommand{\entropy}[1]{ {\text{H}_{\infty}\left({#1}\right)} }
\def\modulo{\text{mod}}
\def\omm{ \{0,1\} }
\def\id{ \textit{id} }

\def\D{{\partial}}
\def\gap{\textsf{gap}}
\def\sign{\textsf{sign}}
\def\bias{\textsf{bias}}
\def\spar{\textsf{sparse}}
\def\span{\textsf{span}}
\def\Part{\textbf{Part}}
\def\Mon{\textbf{Mon}}
\def\sing{\textbf{sing}}
\def\Und{\textsf{Und}}
\def\Comp{\textsf{Comp}}
\def\rank{\textsf{rank}}
\def\poly{\textsf{poly}}
\def\codim{\textsf{codim}}
\def\cp{\textsf{cp}}
\def\uni{\textsf{Uni}}
\def\ext{\textsf{\bf Ext}}
\def\extt{\textsf{\bf Ext2}}
\def\cros{{\textbf{cr}}}
\def\dist{{\textbf{dist}}}
\def\SL{{\textbf{SL}}}
\def\stab{{\textrm{Stab}}}
\def\Dec{{\textrm{Dec}}}

\def\fplus{{\,+_f\,}}
 
\newcommand{\epsclose}{\stackrel{\epsilon}{\thicksim} }
\newcommand{\eclose}[1]{\stackrel{{#1}}{\thicksim} }
\newcommand{\eps}{\epsilon}

\begin{document}

\title{Incidence Theorems -- Lecture Notes}
\date{}
\maketitle


\fi

\section{Sylvester-Gallai type theorems over the reals}\label{sec-sgreal}

The Sylvester-Gallai (SG) theorem states that, in any configuration of $n$ points in the real plane, not all on the same line, there exists a line passing through exactly two of the points. Another way of stating it is as saying that, if in a configuration of points, every pair of points is collinear with a third point, then all points must lie on the same line. This theorem has an extremely simple proof: Suppose the points $v_1,\ldots,v_n$ are not on a line, and let $v_i$ be a point such that the distance between $v_i$ and some line, say $\ell_{v_1,v_2}$, is minimal among all such distances (i.e., between a point and a line defined by the set of points). Now, the line $\ell_{v_1,v_2}$ contains a third point $v_3$. One can draw a picture and see that one of the distances $\dist( v_1, \ell_{v_i,v_3})$, $\dist( v_2, \ell_{v_i,v_3})$,  $\dist( v_3, \ell_{v_i,v_2})$ or  $\dist( v_3, \ell_{v_i,v_1})$ is smaller than $\dist(v_i, \ell_{v_1,v_2})$. \qed

Over the complex numbers this theorem is no longer true! There are configurations of points that lie in a two dimensional plane and with the property that every pair is collinear with a third point. The complex SG theorem, proved by Kelly in \cite{Kel86}, says that this is the highest dimension possible and that every such configuration is contained in some two dimensional affine (complex) plane. The proof of Kelly's theorem  originally used deep tools from Algebraic Geometry but recently an elementary proof was found by Elkies, Pretorius and Swanepoel \cite{EPS06}.

A nice way to think about the SG theorem (and the way which leads to interesting generalizations) is as translating local information (about collinear triples) into global information (all points being on a line). We will now study a more relaxed version of this question when the local information is partial. We start with some definitions. The affine dimension of a set of points $\dim(v_1,\ldots,v_n)$ is the dimension of the smallest affine subspace containing them. Given $v_1,\ldots,v_n$ we call a line passing through at least two points in the set {\em special} if it contains at least three points in the configuration. Otherwise we call the line an {\em ordinary} line. So, the standard SG theorem says that, in every configuration of dimension  at least 2 (or 3 over the complex numbers) there is at least one ordinary line.

\begin{define}[$\delta$-SG configuration]
\label{def: d SG} Let $\delta \in [0,1]$. The $n$ distinct points $v_1,\ldots,v_n \in \C^d$ is called a {\em $\delta$-SG
configuration} if for every $i \in [n]$, there exists a family of special lines $L_i$ all passing through $v_i$ and at least
$\delta n$ of the points $v_1,\ldots,v_n$ are on the lines in $L_i$. (Note that each collection $L_i$ may cover a different
subset of the $n$ points.)
\end{define}

We will now prove the following theorem:
\begin{thm}[Quantitative SG theorem]\label{thm-deltaSG}
Let $\delta \in (0,1]$. Let $v_1,\ldots,v_n \in \C^d$ be a $\delta$-SG configuration. Then
$$\dim\{v_1,\ldots,v_n\} < O(1 / \delta^2).$$
\end{thm}
 
This theorem, proven in \cite{BDWY11}, does not imply Kelly's theorem since, for $\delta=1$ we do not get the constant $2$ (in the original paper the constant 10 is arrived at). A more recent work  \cite{DSW12} improves the techniques in the proof of Theorem~\ref{thm-deltaSG} to give a quantitatively better bound of $O(1/\delta)$, which also gives the constant 2 for $\delta=1$ over the complex numbers  (which gives a new proof of Kelly's theorem).

\subsection{Rank of design matrices}
The proof of Theorem~\ref{thm-deltaSG} is by reduction to a question about the rank of matrices with certain restrictions on their zero/non zero patterns. These are called design matrices:

\begin{define}[Design matrix]\label{def-designmatrix}
Let $A$ be an $m \times n$ matrix over some field. For $i \in [m]$ let $R_i \subset [n]$ denote the set of indices of all non zero entries in the $i$'th row of $A$. Similarly, let $C_j \subset [m]$, $j \in [n]$,
denote the set of non zero indices in the $j$'th column. We say that $A$ is a {\em $(q,k,t)$-design matrix} if
\begin{enumerate}
\item For all $i \in [m]$, $|R_i| \leq q$.
\item For all $j \in [n]$, $|C_j| \geq k$.
\item For all $j_1 \neq j_2 \in [n]$, $|C_{j_1} \cap C_{j_2} | \leq t$.
\end{enumerate}
\end{define}

The reason for studying these matrices in connection with SG configurations will become clear later. For now, let us state the main result we will need to prove:

\begin{thm}[Rank of design matrices]\label{thm-rankdesign}
Let $A$ be an $m \times n$ complex  $(q,k,t)$-design matrix. Then
$$\rank(A) \geq n - \left(\frac{q\cdot t \cdot n}{2k}\right)^2 .$$
\end{thm}

For a stronger (optimal in some settings) form of this theorem we refer the reader to \cite{DSW12}. We will prove Theorem~\ref{thm-rankdesign} in Section~\ref{sec-rankdesign} and will continue now with the proof of Theorem~\ref{thm-deltaSG}.

\subsection{Proof of Theorem~\ref{thm-deltaSG} using the rank bound}

Let $V$ be the $n \times d$ matrix whose $i$'th row is the vector $v_i$. Assume w.l.o.g. that $v_1 = 0$. Thus
$$\dim\{v_1,\ldots,v_n\} = \rank(V) .$$

The overview of the proof is as follows. We will first build an $m \times n$ matrix $A$ that will satisfy $A \cdot V = 0$. Then,
we will argue that the rank of $A$ is large because it is a design matrix. This will show that the rank of $V$ is small.

Consider a special line $\ell$ which passes through three points $v_i,v_j,v_k$. This gives a linear dependency among the three
vectors $v_i,v_j,v_k$ (we identify a point with its vector of coordinates in the standard basis). In other words, this gives a
vector $a = (a_1,\ldots,a_n)$ which is non zero only in the three coordinates $i,j,k$ and such that $a \cdot V = 0$. If $a$ is
not unique, choose an arbitrary vector $a$ with these properties.
Our strategy is to pick a family of collinear triples among the points in our configuration and to build the matrix $A$ from rows
corresponding to these triples in the above manner.

We will need the following combinatorial lemma. 

\begin{lem}\label{lem-triples}
Let $r \geq 3$. 
Then there exists a set  $T \subset [r]^3$ of $r^2 - r$ triples that satisfies the following properties:
\begin{enumerate}
\item Each triple $(t_1,t_2,t_3) \in T$ is of three distinct elements.
\item For each $i \in [r]$ there are exactly $3(r-1)$ triples in $T$ containing $i$ as an element.
\item For every pair $i,j \in [r]$ of distinct elements there are at most $6$ triples in $T$ which contain both $i$ and $j$ as elements.
\end{enumerate}
\end{lem}
\begin{proof}
This follows from a result of Hilton \cite{Hil73} on diagonal Latin squares and we will omit it (see \cite{BDWY11} for more details.)
\end{proof}

Let $\L$ denote the set of all special lines in the configuration (i.e., all lines containing at least three points). Then each
$L_i$ is a subset of $\L$ containing lines passing through $v_i$.  For each $\ell \in \L$ let $V_\ell$ denote the set of points
in the configuration which lie on the line $\ell$. Then $|V_\ell| \geq 3$ and we can assign to it a family of triples $T_\ell
\subset V_\ell^3$, given by Lemma~\ref{lem-triples} (we identify $V_\ell$ with $[r]$, where $r = |V_\ell|$ in some arbitrary
way).

We now construct the matrix $A$ by going over all lines $\ell \in \L$ and for each triple in $T_\ell$ adding as a row of $A$ the
vector with three non zero coefficients $a = (a_1,\ldots,a_n)$ described above (so that $a$ is the linear dependency between the
three points in the triple).

Since the matrix $A$ satisfies $A \cdot V = 0$ by construction, we only have to argue that $A$ is a design matrix and bound its
rank.
\begin{claim}
The matrix $A$ is a $(3,3k,6)$-design matrix,  where $k \triangleq \lfloor \delta n \rfloor - 1$.
\end{claim}
\begin{proof}

By construction, each row of $A$ has exactly $3$ non zero entries. The number of non zero entries in column $i$ of $A$
corresponds to the number of triples we used that contain the point $v_i$. These can come from all special lines containing
$v_i$. Suppose there are $s$ special lines containing $v_i$ and let $r_1,\ldots,r_s$ denote the number of points on each of those
lines. Then, since the lines through $v_i$ have only the point $v_i$ in common, we have that
\[ \sum_{j=1}^s (r_j - 1) \geq k. \]
The properties of the families of triples $T_\ell$ guarantee that there are $3(r_j-1)$ triples containing $v_i$ coming from the $j$'th line. Therefore there are at least $3k$ triples in total containing $v_i$.

The size of the  intersection of columns $i_1$ and $i_2$ is equal to the number of triples containing the points $v_{i_1},v_{i_2}$ that were used in the construction of $A$. These triples can only come from one special line (the line containing these two points) and so, by Lemma~\ref{lem-triples}, there can be at most $6$ of those.
\end{proof}

Applying Theorem~\ref{thm-rankdesign} we get that
\begin{eqnarray*}
 \rank(A) & \geq & n - \left(\frac{3 \cdot 6 \cdot n}{2 \cdot 3k}\right)^2
\geq n - \left(\frac{3 \cdot n}{\delta n  - 2}\right)^2 \\
 & \geq & n - \left(\frac{3 \cdot n \cdot 13 }{ 11 \cdot \delta n} \right)^2 >   n - O(1/\delta^2).
\end{eqnarray*}
Which completes the proof. \qed

\subsection{Extensions to other fields}

We will discuss Sylvester-Gallai type problems over small finite fields in Section~\ref{sec-sgfinite}. For now, let us see that Theorem~\ref{thm-rankdesign} extends to any field of characteristic zero (or very large finite characteristic). Since the reduction to the $\delta$-SG bound was field independent, we can also extend Theorem~\ref{thm-deltaSG} to these fields.

The argument is quite generic and relies on Hilbert's Nullstellensatz.
\begin{define}[$T$-matrix]
Let $m,n$ be integers and let $T \subset [m] \times [n]$. We call an $m \times n$ matrix $A$ a {\em $T$-matrix} if all entries of $A$ with indices in $T$ are non zero and all entries with indices outside $T$ are zero.
\end{define}

\begin{thm}[Effective Hilbert's Nullstellensatz \cite{Kol88}]
Let $g_1,\ldots,g_s \in \Z[y_1,\ldots,y_t]$ be degree $d$ polynomials with coefficients in $\{0,1\}$  and let
$$Z \triangleq \{ y \in \C^t \,|\, g_i(y)=0 \,\, \forall i \in [s] \}.$$
Suppose $h \in \Z[z_1,\ldots,z_t]$ is another polynomial with coefficients in $\{0,1\}$ which vanishes on $Z$.
Then there exist positive integers $C,D$ and polynomials $f_1, \ldots, f_s \in \Z[y_1,\ldots,y_t]$ such that
\[ \sum_{i=1}^s f_i \cdot g_i \equiv C \cdot h^D. \]
Furthermore, one can bound $C,D$ and the maximal absolute value of the coefficients of the $f_i$'s by an explicit function $H_0(d,t,s)$.
\end{thm}

\begin{thm}\label{thm-finitefield}
Let $m,n,r$ be integers and let $T \subset [m] \times [n]$. Suppose that all complex $T$-matrices have rank at least $r$. Let $\F$ be a field of either characteristic zero or of finite large enough characteristic $p > P_0(n,m)$, where $P_0$ is some explicit function of $n$ and $m$. Then, the rank of all $T$-matrices over $\F$  is at least $r$.
\end{thm}
\begin{proof}
Let $g_1,\ldots,g_s \in \C[\{x_{ij} \ | \  i \in [m], j \in [n]\} ]$ be the determinants of all $r \times r$ sub-matrices of an $m \times n$ matrix of variables $X = (x_{ij}$). The statement ``all $T$-matrices have rank at least $r$'' can be phrased as ``if $x_{ij}=0$ for all $(i,j) \not\in T$ and $g_k(X)=0$ for all $k \in [s]$ then $\prod_{(i,j) \in T} x_{ij}=0$.'' That is, if all entries outside $T$ are zero and $X$ has rank smaller than $r$ then it must have at least one zero entry also inside $T$. From Nullstellensatz we know that there are integers $\alpha, \lambda > 0$ and polynomials $f_1,\ldots,f_s$ and $h_{ij}, (i,j) \not \in T$, with integer coefficients such that
\begin{equation}\label{eq-rankpoly}
\alpha \cdot \left( \prod_{(i,j) \in T} x_{ij} \right)^\lambda \equiv \sum_{(i,j) \not\in T} x_{ij} \cdot h_{ij}(X)  + \sum_{k=1}^s f_i(X)\cdot g_i(X).
\end{equation}
This identity implies the high rank of $T$-matrices also over any field $\F$ in which $\alpha \neq 0$. Since we have a bound on $\alpha$ in terms of $n$ and $m$ the result follows.
\end{proof}

\ifediting
\bibliographystyle{alpha}
 \bibliography{incidence}

\end{document}
\fi

\newif\ifediting

\ifediting
\documentclass[11pt]{article}

\usepackage{amsmath,amsthm,amssymb}
\newcommand{\remove}[1]{}
\setlength{\topmargin}{0.3in} \setlength{\headheight}{0in}
\setlength{\headsep}{0in} \setlength{\textheight}{8.0in}
\setlength{\topsep}{0.1in} \setlength{\itemsep}{0.0in}
\parskip=0.05in
 \textwidth=6.5in 
\oddsidemargin=0truecm \evensidemargin=0truecm

\newtheorem{thm}{Theorem}[section]
\newtheorem{claim}[thm]{Claim}
\newtheorem{lem}[thm]{Lemma}
\newtheorem{define}[thm]{Definition}
\newtheorem{cor}[thm]{Corollary}
\newtheorem{obs}[thm]{Observation}
\newtheorem{example}[thm]{Example}
\newtheorem{construct}[thm]{Construction}
\newtheorem{conjecture}[thm]{Conjecture}
\newtheorem{THM}{Theorem}
\newtheorem{question}{Question}
\newtheorem{fact}[thm]{Fact}
\newtheorem{prop}[thm]{Proposition}

\def\F{{\mathbb{F}}}
\def\Q{{\mathbb{Q}}}
\def\Z{{\mathbb{Z}}}
\def\N{{\mathbb{N}}}
\def\R{{\mathbb{R}}}
\def\K{{\mathbb{K}}}
\def\C{{\mathbb{C}}}
\def\A{{\mathbb{A}}}
\def\P{{\mathbb{P}}}
\def\cP{{\cal P}}
\def\cS{{\mathcal S}}
\def\cE{{\mathcal E}}
\def\V{{\mathbf{V}}}
\def\I{{\mathbf{I}}}
\def\bx{{\mathbf x}}
\def\by{{\mathbf y}}
\def\E{{\mathbb E}}

\def\half{ \frac{1}{2}}
\newcommand{\ip}[2]{\langle #1,#2 \rangle}
\def\sumN{\sum_{i=1}^n}
\def\_{\,\,\,\,\,}
\def\prob{{\mathbf{Pr}}}
\newcommand{\entropy}[1]{ {\text{H}_{\infty}\left({#1}\right)} }
\def\modulo{\text{mod}}
\def\omm{ \{0,1\} }
\def\id{ \textit{id} }

\def\D{{\partial}}
\def\gap{\textsf{gap}}
\def\sign{\textsf{sign}}
\def\spar{\textsf{sparse}}
\def\span{\textsf{span}}
\def\Part{\textbf{Part}}
\def\Mon{\textbf{Mon}}
\def\sing{\textbf{sing}}
\def\Und{\textsf{Und}}
\def\Comp{\textsf{Comp}}
\def\rank{\textsf{rank}}
\def\poly{\textsf{poly}}
\def\codim{\textsf{codim}}
\def\cp{\textsf{cp}}
\def\uni{\textsf{Uni}}
\def\ext{\textsf{\bf Ext}}
\def\extt{\textsf{\bf Ext2}}

\def\fplus{{\,+_f\,}}

\newcommand{\epsclose}{\stackrel{\epsilon}{\thicksim} }
\newcommand{\eclose}[1]{\stackrel{{#1}}{\thicksim} }
\newcommand{\eps}{\epsilon}
\newcommand{\Anote}[1]{\begin{quote}{\sf Avi's Note:} {\sl{#1}} \end{quote}}

\begin{document}

\title{Incidence Theorems -- Lecture Notes}
\date{}
\maketitle


\fi

\section{Rank lower bound for design matrices}\label{sec-rankdesign}

We will now prove Theorem~\ref{thm-rankdesign}. First, we discuss a simpler case:

\subsection{The bounded entries case}

When the ratios between different entries of the matrix are bounded in absolute value (say, they are all in $[1/c,c]$ for some  positive constant $c$), the proof is quite easy. Observe that, in this case the $n \times n$ matrix $M = A^* A$ ($A^*$ is the conjugate transpose) is a Hermitian matrix with diagonal elements which are all at least $k/c^2$ in absolute value and the off diagonal are at most $tc^2$ in absolute value.  Also notice that it is enough to give a lower bound on the rank of $M$ since it is equal to the rank of $A$. This bound will follow from the following simple lemma (see  \cite{Alo09} for more on this lemma) which provides a bound on the rank of matrices whose diagonal entries are much larger than the off-diagonal ones.

\begin{lem}\label{lem-diagdom}
Let $A = (a_{ij})$ be an $n \times n$ complex Hermitian matrix and let $0 < \ell < L$ be integers.
Suppose that $a_{ii} \geq L$ for all $i \in [n]$ and that $|a_{ij}| \leq \ell$ for all $i \neq j$.
Then
\[ \rank(A) \geq \frac{n}{1+ n\cdot (\ell/L)^2} \geq n - (n\ell/L)^2. \]
\end{lem}
\begin{proof}
We can assume w.l.o.g. that $a_{ii} = L$ for all $i$.
If not, then we can make the inequality into an equality by multiplying the $i$'th row and column by $(L/a_{ii})^{1/2} < 1$ without changing the rank or breaking the symmetry. Let $r = \rank(A)$ and let $\lambda_1,\ldots,\lambda_r$ denote the non-zero  eigenvalues of $A$ (counting multiplicities). Since $A$ is Hermitian we have that the $\lambda_i$'s are real. We have
\begin{eqnarray*}
 n^2 \cdot L^2 &=& tr(A)^2
 = \left( \sum_{i=1}^r \lambda_i \right)^2
 \leq  r \cdot \sum_{i=1}^r \lambda_i^2
 = r \cdot \sum_{i,j=1}^n |a_{ij}|^2 \\
 &\leq& r \cdot ( n\cdot L^2 + n^2 \cdot \ell^2).
\end{eqnarray*}
Rearranging, we get the required bound. The second inequality in the statement of the lemma follows from the fact that $1/(1+x) \geq 1 - x$ for all $x$.
\end{proof}

Plugging in the parameters in the above lemma we get a rank bound of $n - O(nt/k)^2$ on the rank of $M$ which is what we wanted to prove (in this case, the parameter $q$ is not used). To handle the general case, when the entries are not bounded, we will use a technique called {\em matrix scaling}.

\paragraph{Remark: }Notice that it would suffice if the entries {\em in each row} were bounded (i.e., in a range $[1/c,c]$) since then we could scale each row to get a bounded ratio matrix. Recalling our application to the SG theorem, in which each row corresponded to a collinear triple, we can see that such an {\em unbalanced} triple (where some ratio is very large) can only come from  a collinear triple $v_1,v_2,v_3$ such that the distance between $v_1,v_2$ is, say, much larger than the distance between $v_2,v_3$. Hence, if we know, for some reason, that such triples do not exist in our configuration we can just apply the above argument without a need for further work.

\subsection{Matrix scaling} 

We now define matrix scaling:

\begin{define}\label{def-scaling}[Matrix scaling]
Let $A$ be an $m \times n$ complex matrix. Let $\rho \in \C^{m}, \gamma \in \C^n$ be two complex vectors
with all entries non-zero.
We denote by $$ \SC(A,\rho,\gamma)$$ the matrix obtained from $A$ by multiplying the $(i,j)$'th element of $A$ by $\rho_i \cdot \gamma_j$. We say that two matrices $A,B$ of the same dimensions are a scaling of each other if there exist non-zero vectors $\rho,\gamma$ such that $B = \SC(A,\rho,\gamma)$.
It is easy to check that this is an equivalence relation. We refer to the elements of the vector $\rho$ as the {\em row scaling coefficients} and to the elements of $\gamma$ as the {\em column scaling coefficients}. Notice that two matrices which are a scaling of each other have the same rank and the same pattern of zero and non-zero entries.
\end{define}

Matrix scaling originated in a paper of Sinkhorn \cite{Sinkhorn} and has been
widely studied since (see \cite{LSW98} for more background). The goal is to find a scaling that satisfies certain conditions on the row/column sums. For example, given a square matrix (say, with non negative entries), we would like to find a scaling that makes the matrix doubly stochastic (i.e., with row sums equal column sums equal one). Sinkhorn showed that, if all entries are positive (no zeros) this is possible. The proof was using an iterative algorithm: keep normalizing the row sums and the column sums in alternating steps. This will converge to a scaling that gives a doubly stochastic matrix (for a more efficient variant see \cite{LSW98}). If the matrix contains zeros things are a bit trickier. Take for example the $2 \times 2$ matrix $$ \left(\begin{array}{cc}
1 & 1\\
0 & 1\\
\end{array}\right).
$$ It is clear that there is no scaling of this matrix that makes it doubly stochastic. However, we can `almost' achieve this by making the row/columns sums arbitrarily close to 1. Sometimes, this approximate scaling is good enough, as we shall see in our application. Clearly, we need some condition on the pattern of zeros and non zeros of the matrix (or at least that no row/columns in zero!). The following definition will give a necessary condition that will suffice for our purposes (a more general condition which is both necessary and sufficient is known. See \cite{BDWY11,RS89}).

\begin{define}[Non-zero diagonal]
Let $A$ be an $n \times n$ real matrix. We say that $A$ has a non-zero diagonal if all of its diagonal entries are non-zero. If $A$ is an $nk \times n$ matrix we say that $A$ has non-zero diagonal if its rows can be reordered so that for each $i =0,\ldots,k-1$ the rows $in+1, \ldots,in + n $ form an $n \times n$ matrix with non-zero diagonal (i.e., $A$ is, up to ordering, a concatenation of square non-zero diagonal matrices).
\end{define}

The following is a special case of a theorem from \cite{RS89} that gives sufficient conditions for finding a scaling of a matrix which has certain row and column sums.

\begin{thm}[Matrix scaling theorem]\label{thm-scaling}
Let $A$ be an $nk \times n$ real matrix with non-negative entries and non-zero diagonal. Then, for every $\eps >0$, there
exists a scaling $A'$ of $A$ such that the sum of each row of $A'$ is at most  $1+\eps$ and the sum of each column of $A'$ is at
least $k - \eps$. Moreover, the scaling coefficients used to obtain $A'$ are all positive real numbers.
\end{thm}

The proof of the theorem uses convex programming techniques. One defines an appropriate function and shows that, at the points at which it is maximized, the vanishing of the partial derivatives gives the required bounds on the row/column sums. We will prove this theorem in Section~\ref{sec-proofscaling}.

We will need the following easy corollary of the above theorem.

\begin{cor}[$\ell_2^2$-scaling]\label{cor-l2scale}
Let $A = (a_{ij})$ be an $nk \times n$ complex matrix with non-zero diagonal. Then, for every $\eps > 0$, there exists a scaling $A'$ of $A$ such that for every $i \in [nk]$ $$ \sum_{j \in [n]} |a_{ij}|^2 \leq 1+\eps$$ and for every $j \in [n]$ $$ \sum_{i \in [m]} |a_{ij}|^2 \geq k - \eps.$$
\end{cor}
\begin{proof}
Let $B = (b_{ij}) = ( |a_{ij}|^2 )$. Then $B$ is a real non-negative matrix with non-zero diagonal. Applying Theorem~\ref{thm-scaling} we get that for all $\eps>0$ there exists a scaling $B' = \SC(B, \rho, \gamma)$, with $\rho,\gamma$ positive real vectors,  which has row sums at most $1+\eps$ and column sums at least  $k - \eps$. Letting $\rho'_i = \sqrt{\rho_i}$ and $\gamma'_i = \sqrt{\gamma_i}$ we get a scaling $\SC(A,\rho',\gamma')$ of $A$ with the required properties.
\end{proof}

\subsection{Proof of Theorem~\ref{thm-rankdesign}}
To prove the theorem we will first find a scaling of $A$ so that the norms (squared) of the columns are large and such that each entry is small.

Our first step is to find an $nk \times n$ matrix $B$ with non-zero diagonal that will be composed from rows of $A$ s.t. each row is repeated with multiplicity between $0$ and $q$. To achieve this  we will describe an algorithm that builds the matrix $B$ iteratively by concatenating to it rows from $A$. The algorithm will
\emph{mark} entries of $A$ as it continues to add rows. Keeping track of these marks will help us decide which rows to add next. Initially all the entries of $A$ are \emph{unmarked}. The algorithm proceeds in $k$ steps. At  step $i$ ($i$ goes from $1$ to $k$) the algorithm picks $n$ rows from $A$ and adds them to $B$. These $n$ rows are chosen as follows: For every $j \in \{1,\ldots,n\}$ pick a row that has an unmarked non-zero entry in the $j$'th column and mark this non-zero entry. The reason why such a row exists at all steps is that each column contains at least $k$ non-zero entries, and in each step we mark at most one non-zero entry in each column.
\begin{claim}\label{cla-matrixB}
The matrix $B$ obtained by the algorithm has non-zero diagonal and each row of $A$ is added to $B$ at most $q$ times.
\end{claim}
\begin{proof}
The $n$ rows added at each of the $k$ steps form an $n\times n$ matrix with non-zero diagonal. The bound on the number of times each row is added to $B$ follows from the fact that each row has at most $q$ non-zero entries and each time we add a row to $B$ we mark one of its non-zero entries.
\end{proof}

Since the matrix $B$ is obtained from the rows of $A$ its rank is at most the rank of $A$. Fix some $\eps>0$ (which will later tend to zero).
Applying Corollary~\ref{cor-l2scale} we get a scaling $B'$ of $B$ such that the $\ell_2$-norm of each row is at most $\sqrt{1+\eps}$ and the $\ell_2$-norm of each column is at least $ \sqrt{k -\eps}$.

Our final step is to argue about the rank of $B'$ (which is at most the rank of $A$). To this end, consider the matrix $$ M = (B')^* \cdot B',$$
where $(B')^*$ is $B'$ transposed conjugate.
Then $M = (m_{ij})$ is an $n \times n$ Hermitian matrix.
The diagonal entries of $M$ are exactly the squares of the $\ell_2$-norm of the columns of $B'$.
Therefore, $$ m_{ii} \geq  (k-\eps)$$ for all $i \in [n]$.

We now upper bound the off-diagonal entries.
The off-diagonal entries of $M$ are the inner products of different columns of $B'$.
The intersection of the support of each pair of different columns is at most $tq$ since, in $A$ the intersections were at most $t$, and each row in $A$ is repeated at most $q$ times in $B$ (which has the same support as $B'$).
The norm of each row is at most $\sqrt{1+\eps}$.
For every two real numbers $\alpha,\beta$ so that $\alpha^2 + \beta^2 \leq 1+\eps$
we have $|\alpha \cdot \beta| \leq 1/2 + \eps'$, where $\eps'$ tends to zero as $\eps$ tends to zero.
Therefore
$$ |m_{ij}| \leq tq\cdot (1/2+\eps')$$ for all $i \neq j \in [n]$. Applying Lemma~\ref{lem-diagdom} we get that
$$\rank(A) = \rank(A') \geq n - \left(\frac{q\cdot t(1/2+\eps') \cdot n}{k-\eps}\right)^2.$$
Since this holds for all $\eps > 0$
it holds also for $\eps = 0$, which gives the required bound on the rank of $A$. \qed

\subsection{Proof of the matrix scaling theorem}\label{sec-proofscaling}

For simplicity we will prove the theorem for a square $n \times n$ matrix $A$ (it is easy to modify the proof to fit the more general case). We wish to find a scaling of $A$ with row/columns sums approaching 1. We will call  a scaling with row/column sums {\em exactly} 1, a `good' scaling. Since the constant 1 is arbitrary we will call a scaling `good' even if its row/column sums are equal to some other constant. We will first give a  condition on the pattern of zeros/non-zeros of $A$ that guarantees $A$ has a good scaling. Then we will argue about matrices with non-zero diagonal (which might not have a good scaling as we saw above) and obtain approximate scalings for all $\eps>0$.

We first set up some convenient notations. Let $$ S = \text{supp}(A) = \{ (i,j) \in [n]^2 \,|\, A_{i,j} \neq 0 \}.$$ For $s = (i,j) \in [n]^2$ we will denote $A_s = A_{i,j}$. For each $s = (i,j)\in [n]^2$ let $e_s \in \R^{2n}$ be a vector with $1$ at positions $i$ and $n+j$ and zeros everywhere else. We think of vectors in $\R^{2n}$ as divided into two parts: the first $n$ coordinates corresponding to rows of $A$ and the last $n$ coordinates to columns of $A$. Thus, the vector $e_s$ corresponds to the $(i,j)$'s entry of $A$. We also define $u \in \R^{2n}$ to be the vector with all entries equal to $1/n$. Also, let $\bar 1 \in \R^{2n}$ be the all $1$ vector and notice that $\bar 1 \cdot u = \bar 1 \cdot e_s = 2$ for all $s \in [n]^2$ (where $x \cdot y$ denotes the standard inner product).

With these notations in place we can state our scaling problem in a nicer form. Since we are looking for positive scaling coefficients, it is convenient to treat all of them as exponential functions. Thus, we will find row/column coefficients $\rho_1,\ldots,\rho_n,\gamma_1,\ldots,\gamma_n$ and the scaling defined by them will multiply the rows/columns by $\exp(\rho_i),\exp(\gamma_j)$ (which are always positive). Solving for these exponents will make the problem  easier to analyze as we shall now see.

\begin{claim}\label{good-scaling}
There exists a good scaling of $A$ iff there exists $x \in \R^n$ such that 
\begin{equation}\label{sumexp}
	\sum_{s \in S} \left( A_s \exp(x \cdot e_s) \right)   e_s = u.
\end{equation}
\end{claim}
\begin{proof}
Consider the $i$'s position in $u$. The above equality implies $$\sum_{j \in [n]}A_{i,j} \exp(x_i) \exp(x_{n+j}) = 1/n$$ and so, in the scaling with row coefficients $x_1,\ldots,x_n$ and column coefficients $x_{n+1},\ldots,x_{2n}$, the row sums are all $1/n$. Similarly, the column sums are also $1/n$ and so every $x$ satisfying (\ref{sumexp}) gives a good scaling. Conversely, give a good scaling with positive coefficients, we can take logarithms and find an $x$ solving (\ref{sumexp}).
\end{proof}

Notice that (\ref{sumexp}) implies that $u$ must be in the convex hull of the vectors $e_s, s \in S$. To see this, take the inner product with $(1/2)\bar 1$ -- this implies that the sum of coefficients in the linear combination of the $e_s$'s is one. Thus, this is also a necessary condition for having a good scaling. What we will show below is that this is also a sufficient condition for approximate scaling and that, if $u$ is in the {\em interior} of the convex hull (i.e., if there is a convex combination of the $e_s$'s with all coefficients non zero) then there is a `good' scaling. Below, we will prove the following lemma, which is the hardest part of the proof.

\begin{lem}\label{lem-cond}
For all vectors $v \in \R^{2n}$ such that $\bar 1 \cdot v = \bar 1 \cdot e_s$ for all $s \in S$ the following holds: If $v$ is in the interior of the convex hull of the vectors $e_s, s \in S$ then there exists an $x \in \R^{2n}$ satisfying 	\begin{equation}\label{sumexpv}
		\sum_{s \in S} \left( A_s \exp(x \cdot e_s) \right)   e_s = v.
	\end{equation} If $v = u = (1/n)\bar 1$, this condition implies that $A$ has a good scaling.
\end{lem}

We defer the proof for the subsection below and continue with the proof of Matrix-Scaling theorem. Since $A$ has a non zero diagonal, the vector $u$ is in the convex hull of $e_s, s\in S$ (just take the combination of $e_s$ with $s =(i,i)$ and coefficients $1/n$). However, it might not be in the {\em interior} and so we cannot hope to find a good scaling. However, for every $\eps > 0$ there is a vector $u'$ of distance at most $\eps$ from $u$ that is in the interior. Thus, we can find a scaling of $A$ with row/column sums equal to the entries of $u'$. This implies the existence of an approximate scaling for every $\eps$.

\subsection*{Proof of Lemma~\ref{lem-cond}}
The idea is to define a convex function $f(x)$ such that, if $f$ has a minimizer (w.l.o.g it is a global minimizer since $f$ is convex), the vanishing of the gradient at that point implies the equality in (\ref{sumexpv}). Then we will show that $f$ obtains a minimum by showing that its value grows to infinity when $||x||$ does.

The function we will use is
$$ f(x) =  \ln\left( \sum_{s \in S} A_s \exp(x \cdot e_s)\right)e_s - x \cdot v. $$
A simple application of Cauchy-Schwarz shows that $f$ is indeed a convex function. It is straightforward to verify that the gradient of $f$ (the vector of $2n$ partial derivatives) is
$$ \nabla f(x) = \frac{ \sum_{s \in S} A_s \exp(x \cdot e_s) e_s}{\sum_{s \in S} A_s \exp(x \cdot e_s))} - v $$
and so, if $\nabla f(x) = 0$ then $x$ satisfies (\ref{sumexpv}) up to scaling of $v$ (which can be scaled back to one by a rescaling of $x$).

We now show that $f$ goes to infinity when $||x||$ does. This is not precisely true: let $$ F = \{ y \in \R^{2n} \,|\, y\cdot (e_s - e_{s'}) = 0\, \forall s,s' \in S \} $$ be the subspace of vectors $y$ for which $y \cdot e_s$ is constant for all $s \in S$. Observe that $f(x+y) = f(x)$ for all $x \in \R^{2n}$ and all $y \in F$. Let $E = F^{\bot}$ be the dual subspace to $F$. Hence, we can think of $f$ as a function on $E$ and use the fact that $f$ has a minimizer on $\R^{2n}$ iff it has one on $E$. The fact that $f$ has a minimizer on $E$ will follow from the following claim:

\begin{claim}
There exists a constant $C_1 \in \R$ and a positive constant $C_2 \in R$ such that for all $x \in E$ we have $$ f(x) > C_1 + C_2||x||.$$
\end{claim}
\begin{proof}
Notice that, since $E \cap F = 0$ we have that for all nonzero $x \in E$ the following quantity
$$ \Delta(x) = \max_{s \in S} x \cdot e_s  - \min_{s \in S} x \cdot e_s $$
is positive. Let 
$$ \alpha = \min_{x \in E, ||x||=1} \Delta(x) > 0.$$ Notice that, for all $x$ we have $\Delta(x) \geq \alpha||x||$.

To prove the claim, fix some $x \in E$. Let $s_m,s_M \in S$ be such that $\Delta(x) = x \cdot e_{s_M} - x \cdot e_{s_m}$. Since $v$ is in the interior of the convex hull of $e_s, s \in S$ there are strictly positive coefficients $\lambda_s \in \R, s\in S$ such that $v = \sum_{s \in S}\lambda_s e_s$ and $\sum_{s \in S}\lambda_s = 1$. The following calculation completes the proof:
\begin{eqnarray*}
f(x) &\geq& \ln( A_{s_M}\exp(x \cdot e_{s_M})) - x \cdot v \\
&=& C_1 + x \cdot e_{s_M} - x \cdot v\\
&=& C_1 + \sum_{s \in S} \lambda_s ( x \cdot e_{s_M} - x \cdot e_s) \\
&\geq& C_1 + \lambda_{s_m} (x \cdot e_{s_M} - x \cdot  e_{s_m} ) \\
&=& C_1 + \lambda_{s_m} \Delta(x) \\
&\geq& C_1 + \lambda_{s_m}\alpha ||x||.
\end{eqnarray*}	
\end{proof}
This completes the proof of the lemma.

\ifediting
\bibliographystyle{alpha}
 \bibliography{incidence}

\end{document}
\fi

\newif\ifediting

\ifediting
\documentclass[11pt]{article}

\usepackage{amsmath,amsthm,amssymb}
\newcommand{\remove}[1]{}
\setlength{\topmargin}{0.3in} \setlength{\headheight}{0in}
\setlength{\headsep}{0in} \setlength{\textheight}{8.0in}
\setlength{\topsep}{0.1in} \setlength{\itemsep}{0.0in}
\parskip=0.05in
 \textwidth=6.5in 
\oddsidemargin=0truecm \evensidemargin=0truecm

\newtheorem{thm}{Theorem}[section]
\newtheorem{claim}[thm]{Claim}
\newtheorem{lem}[thm]{Lemma}
\newtheorem{define}[thm]{Definition}
\newtheorem{cor}[thm]{Corollary}
\newtheorem{obs}[thm]{Observation}
\newtheorem{example}[thm]{Example}
\newtheorem{construct}[thm]{Construction}
\newtheorem{conjecture}[thm]{Conjecture}
\newtheorem{THM}{Theorem}
\newtheorem{question}{Question}
\newtheorem{fact}[thm]{Fact}
\newtheorem{prop}[thm]{Proposition}

\def\F{{\mathbb{F}}}
\def\Q{{\mathbb{Q}}}
\def\Z{{\mathbb{Z}}}
\def\N{{\mathbb{N}}}
\def\R{{\mathbb{R}}}
\def\K{{\mathbb{K}}}
\def\C{{\mathbb{C}}}
\def\A{{\mathbb{A}}}
\def\P{{\mathbb{P}}}
\def\cP{{\cal P}}
\def\cS{{\mathcal S}}
\def\cE{{\mathcal E}}
\def\V{{\mathbf{V}}}
\def\I{{\mathbf{I}}}
\def\bx{{\mathbf x}}
\def\by{{\mathbf y}}
\def\E{{\mathbb E}}

\def\half{ \frac{1}{2}}
\newcommand{\ip}[2]{\langle #1,#2 \rangle}
\def\sumN{\sum_{i=1}^n}
\def\_{\,\,\,\,\,}
\def\prob{{\mathbf{Pr}}}
\newcommand{\entropy}[1]{ {\text{H}_{\infty}\left({#1}\right)} }
\def\modulo{\text{mod}}
\def\omm{ \{0,1\} }
\def\id{ \textit{id} }

\def\D{{\partial}}
\def\gap{\textsf{gap}}
\def\sign{\textsf{sign}}
\def\spar{\textsf{sparse}}
\def\span{\textsf{span}}
\def\Part{\textbf{Part}}
\def\Mon{\textbf{Mon}}
\def\sing{\textbf{sing}}
\def\Und{\textsf{Und}}
\def\Comp{\textsf{Comp}}
\def\rank{\textsf{rank}}
\def\poly{\textsf{poly}}
\def\codim{\textsf{codim}}
\def\cp{\textsf{cp}}
\def\uni{\textsf{Uni}}
\def\ext{\textsf{\bf Ext}}
\def\extt{\textsf{\bf Ext2}}

\def\fplus{{\,+_f\,}}

\newcommand{\epsclose}{\stackrel{\epsilon}{\thicksim} }
\newcommand{\eclose}[1]{\stackrel{{#1}}{\thicksim} }
\newcommand{\eps}{\epsilon}
\newcommand{\Anote}[1]{\begin{quote}{\sf Avi's Note:} {\sl{#1}} \end{quote}}

\begin{document}

\title{Incidence Theorems -- Lecture Notes}
\date{}
\maketitle


\fi

\section{Sylvester-Gallai over finite fields}\label{sec-sgfinite}

Let $\F$ denote a finite field of $q$ elements. We can extend our definition of  SG configurations (and $\delta$-SG) to finite fields. To simplify matters we will replace the collinearity condition with linear dependence (i.e., we will assume we have many dependent triples). This will require us to assume that no two points are multiples of each other (or they will be dependent with any third point). We will call a set of points $v_1,\ldots,v_n \in \F^n$ a {\em proper} set if no two points are a constant multiple of each other and the zero point is not in the set (so a proper set is a subset of projective space). 

\begin{define}[SG configuration]
Let $V = \{ v_1,\ldots,v_n \} \subset \F^d$ be a proper set of points. $V$ is  called an SG configuration if for every $i \neq j \in [n]$, there exists $k \in [n] \setminus \{i,j\}$ with $v_i,v_j,v_k$ linearly dependent. $V$ is a $\delta$-SG
 configuration, with $\delta \in [0,1]$ if for each $i$ there are at least $\delta n$ values of $j$ for which there exists $k$ s.t $v_i,v_j,v_k$ are linearly dependent.
\end{define}

To simplify the presentation we will restrict ourselves to the case $\delta=1$, where {\em every} pair is in some dependent triple. The results we will prove can all be generalized  to the more general case, when $\delta$ can be any constant, in a straightforward way. We will mention along the way which changes need to be done to handle arbitrary $\delta$.

Using this definition, one can ask the same question as before: `what is the smallest dimension of a SG configuration?'.  To see that the answer has to be different from the real/complex case notice that the set $V = \F^d$ (taking one representative from each line through the origin) is an SG configuration and so we could have $\dim(V) \geq \log_q n$ (with $n = |V| \sim q^{d-1}$). Also, if $\F$ has characteristic $p$ we can take the set $V = \F_p^d$ (modulo constant multiples) and get $\dim(V) \gtrsim \log_p n$. We will prove two bounds: The first is a generic upper bound of  $\dim(V) \leq O( \log_2 n )$ which holds over {\em any} field \cite{GKST, DvirShpilka06}. The second result will be a bound of the form $\dim(V) \leq O(\log_p n) + \poly(p)$  over prime fields of size $p$ \cite{BDSS11}. This second bound is asymptotically tight, as the $V = \F_p^d$ example shows, for any constant $p$. When $p$ is a growing function of $n$  a bound of the form $O(\log_p n)$ is conjectured to exist\footnote{We are not aware of any results for large finite fields of small characteristic, other than the general $\log_2 n$ bound.}.
 
Another way of stating these two bounds is as saying that, if $V \subset \F^d$ is an SG configuration of dimensions $d = \dim(V)$ then $|V| \geq 2^{\Omega(d)}$ (over any field) or $|V| \geq p^{\Omega(d)}$, when $\F$ is a prime field of size $p < d^{o(1)}$. Thus, the size of the smallest SG configuration of dimensions $d$ grows exponentially with $d$ (with the basis of the exponent being larger for fields of larger characteristic). 

\subsection{The $O(\log_2 n)$ bound}

We will prove this bound in two stages. First we will prove it over $\F_2$ and then see how to handle arbitrary fields in a similar way. Let $V = \{v_1,\ldots,v_n\}$ be an SG configuration in $\F_2^d$. W.l.o.g the dimension of $V$ is equal to $d$ and so, we can perform a linear change of basis so that $v_1,\ldots,v_d$ are the standard basis vectors $e_1,\ldots,e_d$ (with $e_i$ having one in coordinate $i$ and zero elsewhere). We will only use the SG property for $v_1,\ldots,v_d$ (i.e., the fact that for each $i \in [d]$ and each $j$ there is a $k$ s.t $v_i,v_j,v_k$ are dependent). Therefore, the bound $O(\log_2 n)$ will hold also for this special case (and, in this case, it is tight over any field. See below). Observe that, when $v_i = e_i$, if the triple $v_i,v_j,v_k$ is dependent then we have $e_i = v_j + v_k$ or, in other words, $v_j,v_k$ differ only in the $i$'th coordinate. Let $\B = \{0,1\}^d$ be the boolean cube with edges going between vectors that differ in exactly one coordinate. Consider $V$ as a subset of $\B$ and let us try to estimate the number of edges of $\B$ that connect two elements of $V$. For each $e_i$ we have at least $\Omega(n)$ edges in `direction' $i$ (i.e. pairs that differ in the $i$ coordinate alone). This follows from the SG property and the discussion above. Thus, in total we have at least $\Omega(n \cdot d)$ edges inside $V$. We now use the following lemma, which is known as the `edge isoperimetric inequality for the hyper cube' (the bound we prove is not the best possible but it will suffice for our purposes).
\begin{lem}
Let $S \subset \B = \{0,1\}^d$ be some subset of the boolean hypercube. Then there are at most $|S|\log_2|S|$ edges going between elements of $S$. 
\end{lem}
\begin{proof}
The proof is by induction on $d$, where the case $d=1$ is trivial. Write $S_0$ for the set elements on $S$ with 1st bit equal to zero and let $S_1 = S \setminus S_0$. Let $E(S)$ denote the number of edges in $S$ and let $E(S_0),E(S_1)$ be defined similarly. We can think of $S_0,S_1$ as subsets of the $d-1$ dimensional cube and so, by induction, both are bounded by  $|S_0|\log_2 |S_0| $ and  $|S_1|\log_2 |S_1| $ respectively. Observe that the edges in $S$ are divided into three disjoint sets: the edges in $S_0$, the edges in $S_1$ and the edges between $S_0$ and $S_1$. This last set of edges has size at most $\min\{|S_0|,|S_1|\}$ since each element in $S_0$ can have at most one neighbor in $S_1$ and vice versa. We thus have $$ E(S) \leq |S_0|\log_2 |S_0| + |S_1|\log_2 |S_1| + \min\{|S_0|,|S_1|\}.$$ Let $m = |S|$ and consider the function $$ f(x) = x\log_2 x + (m-x)\log_2 (m-x) + x$$ in the range $0 \leq x \leq m/2$ (we think of $x$ as being equal to $\min\{|S_0|,|S_1|\}$). Using some basic calculus we see that $f(x)$ is maximized at the end points at which it is equal to $f(0) = f(m/2) = m\log_2 m$. This implies $E(S) \leq m \log_2 m$ as was required.
\end{proof}

Using the lemma we get
$$ n \cdot d \leq O( n \cdot \log_2 n) $$ or $d = \dim(V) \leq O(\log_2 n)$. 

Now consider an arbitrary field $\F$ (not necessarily finite) and suppose $V \subset \F^d$ with $d = \dim(V)$. We will show how to use $V$ to find a subset of the boolean cube $\B$ that has roughly $n \cdot d$ edges. Suppose $e_i,v_j,v_k$ is a dependent triple as before. Now, there exist non zero field coefficients $a,b$ such that $e_i = av_j + bv_k$. This, however, does not imply that $v_j,v_k$ differ in only the $i$'th coordinate. To be able to derive such a conclusion we would like to have $a = -b$. This would be true if we knew that  $v_i$ and $v_j$ have the same value in some coordinate other than $i$. To make this happen (in most triples) we will normalize each $v_i$ so that its first non zero coordinate is 1. More formally, let $f(v) \in [n]$ be the minimal $\ell \in [n]$ such that the $\ell$'th coordinate of $v$ is non zero. We can  multiply each $v_i$ by a constant so that $(v_i)_{f_{(v_i)}} = 1$ for all $i$ (clearly this keeps the SG property intact). Now, for each $i$ we have a set $M_i$ of $\sim n$ pairs $v_j,v_k$ so that $e_i$ is spanned by $v_j,v_k$. Call a pair  $(v_j,v_k) \in M_i$ `good'  if both $f(v_j)$ and $f(v_k)$ are not $i$. If $e_i$ is spanned by a good pair $(v_j,v_k) \in M_i$ then, we must have $f(v_j) = f(v_k)$ and so, by the above, $e_i = av_j - av_k$ and so $v_j,v_k$ differ only in the $i$'th coordinate. 
\begin{claim}
There are at least $\Omega(n\cdot d)$ good pairs (in all of $M_1,\ldots,M_d$ together). 
\end{claim}
\begin{proof}
The total number of pairs is $\Omega(n \cdot d)$ and so we only have to bound the number of `bad' pairs. Each vector $v_j$ can be `responsible' for a bad pair in only one of the $M_i$'s, namely in $M_{f(v_j)}$. Therefore, the total number of bad pairs is bounded by $O(n)$. This complete the proof (assuming $d$ is larger than some absolute constant).
\end{proof}

We will now reduce to the binary case by mapping each field element randomly to either 0 or 1. Every good pair will remain good with probability at least $1/2$ and so (using expectations) we can find a set $V' \subset \B$ (which might be smaller than $V$) that has at least $\Omega(n \cdot d)$ edges in $\B$. The bound now follows from the isoperimetric inequality.

\paragraph{Remark 1:}The same proof as above works for $\delta$-SG configurations using the lower bound $\Omega(\delta d n)$ on the number of edges inside $V$ and gives a $\Omega(\delta^{-1} \log_2 n)$ upper bound on the dimension.

\paragraph{Remark 2:}To see that the $O(\log_2 n)$ bound is tight for the special case we considered (when we only use dependent triples containing $e_1,\ldots,e_d$) take $V = \{0,1\}^d \subset \F^d$ where $\F$ is any field. For every $e_i$ and every $v \in V$ one of the vectors $v+e_i$ or $v - e_i$ is in $V$ and so we have an SG configuration.

\subsection{The $O(\log_p n)$ bound over prime fields}

For the rest of this section $\F$ will denote a finite field of prime size $p$.
The  example $V = \{0,1\}^d$ shows that, to prove the stronger bound of $O(\log_p n)$ we must go beyond the isoperimetric inequality. The new ideas in the proof will come from additive combinatorics. The SG property can be translated into bounds on the additive growth of the set $V$ (up to some scaling) and these bounds will be exactly those encountered in the Balog-Szemeredi-Gowers theorem, encountered in Section~\ref{sec-bsg}. We rephrase this theorem here in a slightly different form (whose proof is a simple reduction to the one we saw).

\begin{thm}\label{thm-bsg2}[Balog-Szemeredi-Gowers]
Let $A \subset G$ be a set of size $N$ in an abelian group $G$. Suppose that $$ |\{ (a_1,a_2) \in A^2 \,|\, a_1 + a_2 \in A \}| \geq N^2/K.$$
Then, there exists a subset $A' \subset A$ with $|A'| \geq N/K^c$ and with $|A'+A'| \leq K^c N$, where $c >0$ is some absolute constant.
\end{thm}

Another ingredient we will need is the following important result of Ruzsa:
\begin{thm}[Ruzsa \cite{Ruzsa}] 
Let $A \subset \F^d$ be such that $|A + A| \leq K|A|$. Then, there exists a subspace $W \subset \F^d$ containing $A$ with $|W| \leq K^c p^{K^c}|A|$, where $c$ is an absolute constant. This implies $\dim(A) \leq \log_p |W| \leq \log_p |A| + K^{c'}$ for some other constant $c'$. 
\end{thm}

We will prove Ruzsa's theorem below and continue with our proof of the $O(\log_p n)$ bound for SG configurations. The first part of the proof will use the above two theorems to find a large subset of $V$ of small dimension. 

\begin{lem}[Small dim subset]
Let $V = \{v_1,\ldots,v_n\} \subset \F^d$ be an SG configuration. Then there exists a subset $V' \subset V$ with $|V'| \geq n/p^c$ such that $\dim(V') \leq \log_p n + p^c$ for some constant $c >0$.
\end{lem}
\begin{proof}
Let $A = \{ \lambda v_i \,|\, i \in [n], \lambda \in \F^* \}$ be the set of size $(p-1)n$ containing all non zero constant multiples of elements from $V$ (recall that no two elements of $V$ are constant multiples of each other). Every dependent triple $v_i,v_j,v_k$ in $V$ with, say, $a_iv_i + a_jv_j = a_kv_k$ implies that the sum of the two elements $a_iv_i,a_jv_j$ (both in $A$) is also in $A$. Using the SG property and this observation we get that, for each $a \in A$ there are at least $\Omega(|A|/p)$ elements $a' \in A$ such that $a + a' \in A$. Using the BSG theorem we get that there is a subset $A' \subset A$ of size $|A'| \geq |A|/p^c$ such that $|A' + A'| \leq p^c|A|$. Ruzsa's theorem now implies that $$\dim(A') \leq \log_p |A'| + p^{c'} \leq \log_p n + \poly(p).$$ We can now take $V' \subset V$ to be set of all elements that have some multiple in $A'$. The size of $V'$ is at least $|A'|/p$ and its dimension is bounded by that of $A'$. This completes the proof.
\end{proof}

We will now show how to `grow' the set $V'$ so that it contains the entire set $V$ without increasing its dimension by much. Let $V'$ be a subset of $V$ given by the Lemma. W.l.o.g we may assume that $\span(V') \cap V = V'$ (otherwise replace $V'$ with its span in $V$). Let $w \in V \setminus V'$ be some element not in $V'$ (and so also not in the span of $V'$). Using the SG property we know that for each $v' \in V'$ there is some $u  \in V$ such that $w,v',u$ are dependent. Since $w$ is not spanned by $V'$ we cannot have $u \in V'$. Thus, we can define a function $f : V' \mapsto V \setminus V'$ such that for all $v' \in V$ we have $w,v',f(v')$ dependent. Observe that if $f(v') = f(v'') = u$ then both $v'$ and $v''$ are in the span of $w,u$ which has size at most $p^2$. This implies that the set of images $f(V')$ has size at least $|V'|/p^2$. Now, let $V'' = \span( V' \cup \{w\}) \cap V$. I.e., add $w$ to $V' $ and take the span of the resulting set inside $V$. Clearly $\dim(V'' )= \dim(V') + 1$. We also add all the elements of $f(V')$ to the set $V''$ since they are spanned by $w$ and some element of $V'$. This means that $|V''| \geq |V'|(1 + 1/p^2)$. Continuing in this manner $\poly(p)$ times we will eventually add all the elements of $V$ and the dimension will grow by an additive factor of $\poly(p)$. This implies that $\dim(V) \leq \log_p n + \poly(p)$ as was required.

\subsubsection{Proof of Ruzsa's theorem}

Let $k\cdot A = A+ \ldots +A$ $k$ times. For $g \in G$ we denote by $A + g = \{a+g\,|\, a \in A\}$. We can assume w.l.o.g that $A = -A = \{-a \,|\, a \in A\}$ since otherwise we can replace $A$ with $A \cup -A$, which will also not grow in addition using Ruzsa calculus.

Consider a maximal integer $r$ such that there exist $b_1,\ldots,b_r \in 3 \cdot A$ for which the $r$ sets $A + b_i, i \in [r]$ do not intersect each other. Notice that each of these $r$ sets is contained in the set $4 \cdot A$, which, using Ruzsa calculus, has size at most $\leq K^c |A|$. Thus, $r \leq K^c$. By construction, we have that for every $b \in 3\cdot A$ the set $b + A$ intersects  $b_i + A$ for some $i \in [r]$. Such an intersection implies that $b \in A - A + b_i = 2 \cdot A + b_i$. We can thus conclude that $$3 \cdot A \subset \cup_{i \in [r]}\left(2\cdot A + b_i\right).$$ Iterating, this means that $$k \cdot A \subset 2\cdot A + \span(b_1,\ldots,b_r)$$ for all $k$. Thus, the span of $A$ has size at most $$ |\span(A)| \leq |2\cdot A| \cdot p^r \leq K^c p^{K^c} \cdot |A|. $$ This concludes the proof.

\ifediting
\bibliographystyle{alpha}
 \bibliography{incidence}

\end{document}
\fi
 
\newif\ifediting

\ifediting
\documentclass[11pt]{article}

\usepackage{amsmath,amsthm,amssymb}
\newcommand{\remove}[1]{}
\setlength{\topmargin}{0.3in} \setlength{\headheight}{0in}
\setlength{\headsep}{0in} \setlength{\textheight}{8.0in}
\setlength{\topsep}{0.1in} \setlength{\itemsep}{0.0in}
\parskip=0.05in
 \textwidth=6.5in 
\oddsidemargin=0truecm \evensidemargin=0truecm

\newtheorem{thm}{Theorem}[section]
\newtheorem{claim}[thm]{Claim}
\newtheorem{lem}[thm]{Lemma}
\newtheorem{define}[thm]{Definition}
\newtheorem{cor}[thm]{Corollary}
\newtheorem{obs}[thm]{Observation}
\newtheorem{example}[thm]{Example}
\newtheorem{construct}[thm]{Construction}
\newtheorem{conjecture}[thm]{Conjecture}
\newtheorem{THM}{Theorem}
\newtheorem{question}{Question}
\newtheorem{fact}[thm]{Fact}
\newtheorem{prop}[thm]{Proposition}

\def\F{{\mathbb{F}}}
\def\Q{{\mathbb{Q}}}
\def\Z{{\mathbb{Z}}}
\def\N{{\mathbb{N}}}
\def\R{{\mathbb{R}}}
\def\K{{\mathbb{K}}}
\def\C{{\mathbb{C}}}
\def\A{{\mathbb{A}}}
\def\P{{\mathbb{P}}}
\def\cP{{\cal P}}
\def\cS{{\mathcal S}}
\def\cE{{\mathcal E}}
\def\V{{\mathbf{V}}}
\def\I{{\mathbf{I}}}
\def\bx{{\mathbf x}}
\def\by{{\mathbf y}}
\def\E{{\mathbb E}}

\def\half{ \frac{1}{2}}
\newcommand{\ip}[2]{\langle #1,#2 \rangle}
\def\sumN{\sum_{i=1}^n}
\def\_{\,\,\,\,\,}
\def\prob{{\mathbf{Pr}}}
\newcommand{\entropy}[1]{ {\text{H}_{\infty}\left({#1}\right)} }
\def\modulo{\text{mod}}
\def\omm{ \{0,1\} }
\def\id{ \textit{id} }

\def\D{{\partial}}
\def\gap{\textsf{gap}}
\def\sign{\textsf{sign}}
\def\spar{\textsf{sparse}}
\def\span{\textsf{span}}
\def\Part{\textbf{Part}}
\def\Mon{\textbf{Mon}}
\def\sing{\textbf{sing}}
\def\Und{\textsf{Und}}
\def\Comp{\textsf{Comp}}
\def\rank{\textsf{rank}}
\def\poly{\textsf{poly}}
\def\codim{\textsf{codim}}
\def\cp{\textsf{cp}}
\def\uni{\textsf{Uni}}
\def\ext{\textsf{\bf Ext}}
\def\extt{\textsf{\bf Ext2}}

\def\fplus{{\,+_f\,}}

\newcommand{\epsclose}{\stackrel{\epsilon}{\thicksim} }
\newcommand{\eclose}[1]{\stackrel{{#1}}{\thicksim} }
\newcommand{\eps}{\epsilon}
\newcommand{\Anote}[1]{\begin{quote}{\sf Avi's Note:} {\sl{#1}} \end{quote}}

\begin{document}

\title{Incidence Theorems -- Lecture Notes}
\date{}
\maketitle


\fi

\section{Locally Correctable Codes}\label{sec-lcc}

We will now see how the question of bounding the dimension of $\delta$-SG configurations comes up naturally in the context of error correction.

\subsection{Error Correcting Codes}

We start by defining Error Correcting Codes (ECCs). We will focus on linear ECCs since these are the most well studied. One way to view  an ECC is as a subspace $C \subset \F^n$, where $\F$ is some finite field. We say that the code has minimum distance $D$ if for all $x \neq y \in C$ the Hamming distance (the number of different coordinates) between $x$ and $y$, denoted $\Delta(x,y)$, is at least $D$. We will sometimes also refer to the normalized minimum distance as the minimum distance divided by $n$. The {\em rate} of the code $C$ is defined as $r(C) = \dim(C)/n$. Codes with high rate and high distance can be used to transmit messages in the presence of errors. More precisely, suppose $\dim(C)=d$ and let $E_C : \F^d \mapsto \F^n$ be a linear mapping whose image is $C$ (thus, $E_C$ is a bijection on its image). To send a message $x \in \F^d$ we send its encoding $y = E_C(x)$ instead. Now, suppose that the transmission is noisy and that the actual received string was not $y$, but some $y' \in \F^n$ with $\Delta(y,y') < D/2$. The receiver could then determine $y$ (and from it, $x$) uniquely from $y'$ since there could not be two distinct $y_1,y_2 \in C$ with both $\Delta(y_1,y')$ and $\Delta(y_2,y')$ smaller than $D/2$ (this would imply $\Delta(y_1,y_2) < D$). 

A nice example of an error correcting code is the Reed Solomon Code: We take $\F$ to be a finite field of size at least $n$ and fix some $n$ distinct field elements $a_1,\ldots,a_n \in \F$. Fixing $d$ to be any integer between $1$ and $n$ we let $C$ be the following subspace: $$ C = \{ (f(a_1),\ldots,f(a_n)) \in \F^n \,|\, f \in \F[T] \text{ has } \deg(f) \leq d \}.$$ That is, $C$ is the subspace of $n$-tuples that are the evaluations of degree $\leq d$ univariate polynomials on $n$ distinct points in the field. The minimal distance of this code can be readily computed since two polynomials of degree $\leq d$ can agree on at most $d$ places. This means that the distance between two distinct vectors in $C$ is at least $n-d$. Thus, if we take, for example, $d =  n/2$ we will get a code of rate $\sim 1/2$ and (normalized) minimum distance $1/2$. Such codes, with constant rate and distance, are sometimes called `good' codes (for obvious reasons). Obtaining good codes over smaller alphabet (idealy, for applications, with $|\F|=2$) can also be obtained using additional ideas. It is important to note here that taking $C$ to be a {\em random} subspace of dimension $d$ will result, with high probability, with a good code. However, this type of construction will not give us any efficient way to perform the decoding (other than going over all elements of $C$). 

Coding theory is a vast area of research spanning engineering, computer science and mathematics and we will not attempt to give a full introduction here. The basic questions on existence/constructions of ECCs of the form described above have, to the most part, satisfactory (if not complete) answers. Our focus will be a specific kind of ECCs -- Locally Correctable Codes (LCCs)-- that are very poorly understood and tightly related to questions regarding SG configurations. LCC's are variants of Locally Decodable Codes (LDC), first defined and studied in a paper by  Katz and Trevisan \cite{KTldc} and much of the discussion below appeared in that seminal work.

\subsection{Locally Correctable Codes}
In the usual ECC setting, the decoder takes a received word $y' \in \F^n$, runs some sophisticated algorithm on $y'$ and returns the unique $y \in C$ which minimizes $\Delta(y,y')$. For example, in  Reed Solomon codes, given a noisy list of values of a polynomial of low degree, we want to interpolate the unique polynomial that agrees with this list in the largest number of places. This type of decoding algorithm is usually very `global', meaning that if one wanted to compute even one coordinate in the `corrected' $y$, they would still need to compute the entire $y$ (and then output a single coordinate). Locally Correctable Codes allow the receiver to recover $y$ from $y'$ in a more local way: The decoder can, given an index $i \in [n]$, recover the $i$'th coordinate of the unique closest $y \in C$, looking at a small random sample of positions in $y'$. Such a decoding procedure cannot always be correct (since the few places we look at might all contain errors), but it could be correct w.h.p over the choices of the coordinates we choose to read. To make the connection to SG configurations clearer we will define a code $C \subset \F^n$ of dimension $d$ as an ordered list of vectors $V = (v_1,\ldots,v_n) \in (\F^d)^n$ (possibly containing repetitions), each corresponding to a single coordinate in $[n]$. Given such $V$ we  take as our code the subspace 
\begin{equation}\label{eq-CV}
	C_V = \{ (\ip{x}{v_1},\ldots,\ip{x}{v_n}) \in \F^n \,|\, x \in \F^d \}.
\end{equation}
Notice that, in this way of writing things, if some $v_i$ is in the span of some other set of vectors $\{ v_{j_1}, \ldots,v_{j_r} \}$  in the list $V$ then the $i$'th position of any $y \in C_V$ can be recovered from the positions $y_{j_1},\ldots,y_{j_r}$. Simply write $$ v_i = \sum_{\ell=1}^r a_\ell v_{j_\ell} $$ and then we have $$ y_i = \ip{x}{v_i} = \sum_{\ell=1}^r a_\ell \ip{x}{v_{j_\ell}} = \sum_{\ell=1}^r a_{\ell}y_{j_\ell}. $$
 
We now give a formal definition of LCCs. We will allow the base field to be any field (even infinite). 
\begin{define}\label{def-lcc}
An $(r,\delta)-LCC$ of dimension $d$ is an ordered list of vectors $V = (v_1,\ldots,v_n) \in (\F^d)^n$ such that $\dim(v_1,\ldots,v_n) = d$ and with the following property, called the {\em LCC property}: for each $i \in [n]$ and every set $S \subset [n]$ of size at most $\delta n$ there exists a set $R \subset [n]\setminus S$ with $|R|\leq r$ such that $v_i \in \span\{v_j \,|\, j \in R\}$. The parameter $r$ is called the {\em query complexity} of $V$.
\end{define} 

It is not immediately obvious why this definition is the right one. True, for every set of `errors' $S \subset [n]$ of size at most $\delta n$ there are $r$ positions outside this set (so the values there are correct) which determine the $i$'th coordinate. So, if we `knew' where the errors were, we could locally correct any coordinate\footnote{This kind of decoding is sometimes interesting in its own right.}. But what if we do not know where the errors are? The following simple and useful lemma will help us resolve this issue. To state the lemma we will require the following definition:

\begin{define}[$r$-Matching]
Let $\Omega$ be some finite set. A family of subsets $M = \{R_1,\ldots,R_k\}$ with each $R_i \subset \Omega$ is called an $r$-Matching in $\Omega$ if 
\begin{itemize}
	\item For all $i \in [k]$, $1 \leq |R_i| \leq r$.
	\item For all $i \neq j \in [k]$,  $R_i \cap R_j = \emptyset$.
\end{itemize}
We denote the {\em size} of the $r$-matching  $M$ by $|M| = k$. We call $M$ a {\em regular} $r$-Matching if all sets $R_i$ are of size exactly $r$. When $r$ is obvious from the context we will sometimes omit it and refer to $M$ simply as a {\em matching}.
\end{define}
	
\begin{lem}
Let $V = (v_1,\ldots,v_n) \in (\F^d)^n$ be an $(r,\delta)-LCC$. Then, for each $i \in [n]$ there exists an $r$-Matching $M_i = \{R_{i,1},\ldots,R_{i,k}\}$ in $[n]$ with $|M_i| = k \geq (\delta/r)n$ such that for every $i \in [n], j \in [k]$ we have $v_i \in \span\{ v_\ell \,|\, \ell \in R_{i,j} \}$. 
\end{lem}
\begin{proof}
For each $i$ we can construct the $r$-Matching $M_i$ iteratively: As long as $|M_i| < (\delta/r)n$ the $r$-tuples in $M_i$ can cover at most $\delta n$ other coordinates and so there has to be an $r$-tuple we can add that is disjoint from all of them.
\end{proof}

In other words, for every  $i \in [n]$ there is a large (at least $\Omega(n)$ if both $\delta$ and $r$ are constants) family of small disjoint sets of  coordinates that determine the $i$'th coordinate. Now, if the fraction of errors is at most $\delta' \ll \delta/r$ only a small fraction of the sets in each $M_i$ will contain some corrupted coordinate and so, picking a random set in $M_i$ will not contain errors w.h.p, allowing for correction of the $i$'th position. This is stated precisely by the following lemma, which justifies Definition~\ref{def-lcc} and connects it with our intuitive description of LCCs.

\begin{lem}\label{lem-qmatching}
Let $V = (v_1,\ldots,v_n) \in (\F^d)^n$ be an $(r,\delta)-LCC$ and let $C_V$ be defined as in (\ref{eq-CV}). Let $\delta' = \eps\delta/r $. Then, there exists an efficient randomized algorithm\footnote{We assume our algorithm can perform field operations at unit cost.} $\Dec : \F^n \times [n] \mapsto \F$ with the following properties
\begin{itemize}
	\item For all $y \in C_V$ and all $y' \in \F^n$ with $\Delta(y,y') \leq \delta' n$ we have $ \Pr[ \Dec(y',i) = y_i] \geq 1 - \eps$, where the probability is over the internal coin tosses of $\Dec$. 
	\item For all $y' \in \F^n$ and all $i \in [n]$, the invocation of $\Dec(y,i)$ reads at most $r$ positions in the input $y'$. 
\end{itemize}
\end{lem}
\begin{proof}
$\Dec(y',i)$ will simply pick a random $j \in [k]$, where $k$ is the size of the $r$-Matchings $M_i = \{ R_{i,1},\ldots,R_{i,k} \}$ given by Lemma~\ref{lem-qmatching}, and compute $y_i$ from the coordinates $\{ y_\ell \,|\, \ell \in R_{i,j} \}$ using the fact that $v_i \in \span\{ v_\ell \,|\, \ell \in R_{i,j} \}$. Since the distance $\Delta(y,y')$ is at most $\delta'n$ there could be at most $\delta' n = \eps \cdot (\delta/r)n \leq \eps |M_i|$ sets $R_{i,j} \in M_i$ that contain a coordinate in which $y$ and $y'$ differ. Thus, with probability at least $1 - \eps$ the decoding will succeed. 
\end{proof}

When $r$ is not constant (say $r = \log n$) LCCs are still interesting but the loss of $1/r$ in the decoding distance is no longer acceptable. A more restrictive  definition of LCCs can be made along the lines of the last lemma, requiring that there exists a decoding procedure $\Dec(y,i)$ that returns $y_i$ with high probability in the presence of $\delta n$ errors, reading only $r$ positions. Here, we opt for the cleaner statement given in  Definition~\ref{def-lcc}. Notice that, for the purpose of proving upper bounds on the dimension of an LCC $V$, our definition is {\em more} general and upper bounds for our definition will imply upper bounds for the stronger definition.

\subsection{Random codes are not locally correctable}
 
The property of being able to decode symbols of the codeword locally is very appealing for real life coding applications. However, for such codes to be used in practice their dimension, which determines the amount of information they can encode, cannot be too small. There is a huge gap between the known upper and lower bounds on the dimension of LCCs with small $r$. This is surprising considering the good understanding we have of `regular' ECCs (without local correction). A partial explanation for this discrepancy is the fact that a random code (of reasonable dimension) is {\em not} an LCC. More formally, suppose $|\F|=q$ and pick the list $V = (v_1,\ldots,v_n)$ at random (i.e., pick each $v_i$ i.i.d in $\F^d$). The probability of any $r+1$ of the chosen vectors to be dependent is at most $q^{r-d}$ (this bound is the probability that the last vector is in the span of the previous $r$). This probability is exponentially small when $r \ll d$ and so, unless $|V| \geq q^{\Omega(d)}$ (i.e., $\dim(V) = O(\log_q n)$) we will not see even {\em one} dependent $r+1$-tuple. Since an LCC must have a quadratic number of  dependent $r+1$-tuples, this shows that a random code is not an LCC\footnote{Random codes in the regime $\dim(V) \sim \log_q n$ are studied in \cite{KS10} and can be shown to have local decoding properties.}. Thus, the construction of LCCs with high dimension and low query complexity is morally different than the construction of regular ECCs. Constructing an ECC amounts to finding a structured example of an object that exists almost everywhere. Constructing an LCC is a task of finding a very rare object with extremely delicate local properties.

\subsection{2-Query LCCs and SG configurations}

The case $r=1$ is not very interesting since it is easy to show that the best (and only) $(1,\delta)$-LCC will have dimension at most $1/\delta$ (since every coordinate must repeat, up to constant multiples, at last $\delta n$ times). The case $r=2$ is already much more interesting and, in this case, we have a pretty good understanding of the parameters obtainable by LCCs. This case is also where the connection to SG configuration will become clear. Suppose $\F$ is finite field of size $q$. A trivial construction of a $(2,\delta)$-LCC with constant $\delta$ is to take $V$ to contain all vectors in $\F^n$. Not accidentally, this is also the trivial construction of an SG configuration. If it is not clear by now, let us state the following easy Lemma:
\begin{lem}
Suppose $V = \{v_1,\ldots,v_n\} \subset \F^d$ is a $\delta$-SG configuration (see Section~\ref{sec-sgreal}). Then the list  $V = (v_1,\ldots,v_n)$ is a $(2,\delta/3)$-LCC.
\end{lem} 
\begin{proof}
Suppose $S \subset [n]$ is some set of size $|S| \leq (\delta/3)n$ and fix some $i \in [n]$. We need to show that there is a pair $j,k \in [n]\setminus S$ such that $v_i$ is spanned by $v_j,v_k$. Using the $\delta$-SG property we know that there is a set $T \subset [n]\setminus S$ of size $|T| \geq (\delta/2)n$ such that for each $j \in T$ there is some $k = k(j) \in [n]$ such that $v_i,v_j,v_k$ are linearly dependent (recall that in $V$ not two vectors are a constant multiple of each other). If there is some $j \in T$ for which $k(j) \not\in S$ we are done since $v_j,v_{k(j)}$ span $v_i$ and both are outside $S$. If for all $j \in T$ we have $k(j) \in S$ then there must be a collision (since $|T| > |S|$) of the form $k(j) = k(j')$ and then both $v_j,v_{j'}$ are in the span of $v_i,v_{k(j)}$. Since $v_j,v_{j'}$ are independent, they must also span $v_i$ and, again, we are done since both are outside $S$.
\end{proof}

Hence, every $\delta$-SG configuration (over any field, not only finite fields) gives a 2-query LCC of distance $\sim \delta$. But is the opposite also true? Can we take any $(2,\delta)$-LCC and convert it to a $\delta'$-SG configuration with $\delta' \sim \delta$? I do not know the answer to this question but suspect that it might be true. The main (and only) difficulty is that an LCC $V$ is a list that can have repetitions. It seems reasonable to conjecture that repetitions should not be useful in creating a good LCC and that, perhaps with some careful combinatorial work, they can be eliminated (at some negligible cost to the other parameters). Even though we do not know of a black-box reduction from 2-query LCCs to SG configurations, all known upper bounds on the dimension of SG configurations extend (with significantly more work) to 2-query LCCs. Suppose $V$ is a $(2,\delta)$-LCC over a field $\F$. The known bound can be summarized as follows:
\begin{itemize}
	\item Over any field $\dim(V) \leq O((1/\delta)\log_2 n)$ \cite{GKST,DvirShpilka06}. 
	\item Over prime fields $\F_p$ we have $\dim(V) \leq \poly(p/\delta) + O( (1/\delta)\log_p n)$ \cite{BDSS11}.
	\item Over fields of characteristic zero (or characteristic $\gg \exp\exp(n)$) we have $\dim(V) \leq \poly(1/\delta)$ \cite{BDWY11}. 
\end{itemize}
The proofs of all three bounds use the same basic ideas used for SG configurations with an extra layer of arguments added to handle repetitions in $V$. Since these arguments are quite cumbersome and taylor made for each proof, it would be very desirable to find a clean black-box way of getting rid of repetitions for any $LCC$ (even with more than two queries).

\subsection{Constructions using polynomials}
 
When $r > 2$ our knowledge is quite limited. The best constructions  are those coming from multivariate polynomials or {\em Reed-Muller} codes \cite{Reed,Muller}. We have already encountered these when we discussed the polynomial method over finite fields. Let $\F$ be a field of size $q$ and let $n = q^s$ for some $s$. We will construct an $(r,\delta)$-LCC $V = (v_1,\ldots,v_n)$ in $\F^n$ by describing the subspace $C_V \subset \F^n$ (see Eq.\ref{eq-CV}). Let $\F^{(e)}[z_1,\ldots,z_m]$ be the set of polynomials in $m$ variables of degree at most $e$. We identify the set of coordinates $[n]$ with the set $\F^m$ using some fixed one-to-one map $\tau : [n] \mapsto \F^m$ and define
$$ C_V = \{ (f(\tau(1)), \ldots, f(\tau(n))) \in \F^n \,|\, f \in \F^{(q-2)}[z_1,\ldots,z_m] \}. $$ That is, a codeword in $C$ is the vector of evaluations of a polynomial of degree $ \leq q-2$ on the entire space $\F^m$.\footnote{One can easily come up with the explicit vectors in the list $V = (v_1,\ldots,v_n)$  and this is a good exercise.} We will now argue that this code is a $(q-1,\delta)$-LCC (with some constant $\delta$). Consider an index $i \in [n]$ and its associated point $\tau(i)$. On every line in $\F^m$ passing through $\tau(i)$, the values of a degree $q-2$ polynomial in $q-1$ places on the line determine the rest of the values on the line. Thus, the value of a codeword at coordinate $\tau(i)$ can be determined from any $q-1$-tuple of coordinates corresponding to the points on any line through $\tau(i)$. Since the space $\F^m$ can be covered completely by lines passing through $\tau(i)$ we can find such a line outside any set $S \subset \F^m$ with $|S| < q^{m-1}$. Thus, we can take $\delta = 1/q$. We can make $\delta$ independent of $q$ by reducing the degree of the polynomials from $q-2$ to , say, $q/10$. Then, we only need to find a line passing through $\tau(i)$ with at least $q/10+1$ points outside $S$. A simple probabilistic argument shows that such a line exists if $|S| < n/10$. The dimension of the code $C_V$ described above is equal to the number of coefficients in a degree $\sim q$ polynomial in $m$ variables, where $m = \log_q n$. When $r \sim q$, the number of queries, is a constant and $n$ tends to infinity, this dimension is roughly $m^{O(r)} = (\log_q n)^{O(r)}$. This is a power of $\log_q n$ that depends linearly on the number of queries. Thus, these codes are quite far from being applicable in practice when we wish to have dimension close to $n$ (or at least polynomial in $n$). 

\subsection{General upper bounds on the dimension of LCCs}

When $r > 2$, the upper bounds on the dimension of $r$-query LCCs are quite weak. The following bound works for any LCC and degrades quite quickly with the number of queries:
\begin{thm}[Katz-Trevisan \cite{KTldc}]
Let $V = (v_1,\ldots,v_n)$ be an $(r,\delta)$-LCC in $\F^n$. Then, when $n$ goes to infinity and $r$ and $\delta$ are fixed, we have $$ \dim(V) \leq  O\left(n^{\frac{r-1}{r}} \cdot \log n \right) .$$
\end{thm}
\begin{proof}
First, we use Lemma~\ref{lem-qmatching} to find  $r$-matchings $M_1,\ldots,M_n$ in $[n]$ of the form $M_i = \{ R_{i,1}, \ldots,R_{i,k} \}$ with $k \geq (\delta/r)n$. Recall also that each of  sets $\{ v_\ell \,|\, \ell \in R_{i,j} \}$ spans the vector $v_i$. We will use a probabilistic argument to find a set $T \subset [n]$ of small size such that $T$ will contain at least one set $R_{i,j}$ for each of the $i$'s. This will imply that $\dim(V) \leq T$ since, the vectors $v_i, i \in T$ span all of the other vectors in $V$.

Consider the following random choice of $V$: take every element in $[n]$ to be in $T$ independently with probability $\mu = \log n \cdot n^{- \frac{1}{r}}$. Then, with  probability higher than $3/4$, we will have $|T| \leq O\left(n^{\frac{r-1}{r}}\cdot \log n\right)$. We now show that w.h.p $T$ will contain at least one set from each of the matchings. The probability for a single set $R_{i,j}$ of size at most $r$ to {\em not} be contained in $T$ is at most $$ \Pr[ R_{i,j} \not\subset T] \leq 1 - \mu^r. $$ Since the sets in each $M_i$ are disjoint we have that $$ \Pr[ \forall j \in [t], \,\, R_{i,j} \not\subset T] \leq (1-\mu^r)^k.$$ Plugging in $\mu$ and the bound for $k$ we get that this probability is smaller than $1/20n$ and so, the probability that there exists an $M_i$ with all sets not in $T$ is at most $1/20$ by a union bound. This means that there exists a choice of $T$ of the appropriate size that contains a set from each matching. This completes the proof.
\end{proof}
 
If this bound were tight then there could be a chance to use LCCs with a constant number of queries in practice. However, most people believe that this bound is not tight and some conjecture that the polynomial constructions achieve the best possible parameters. The best known general upper bound was proved by Woodruff \cite{Woodruff} and gives $\dim(V) \leq \tilde O\left( n^{\frac{\lceil r/2 \rceil - 1}{\lceil r/2 \rceil}}\right)$ which is equal to $\sqrt n$ for $r=3,4$. Any improvement to either this upper bound or the polynomial constructions will be extremely interesting. Notice that the polynomial constructions we saw work only over finite (and smaller than $n$) characteristic. When the characteristic is zero (or larger than $n$) there are no known constructions of constant query LCCs with dimension tending to infinity (one can always take a code of dimension $1/\delta$ with $r=1$). A tempting conjecture, which might be a good starting point for progress is that there are no $3$-query LCCs over fields of characteristic zero.

\subsection{LCCs as low rank sparse matrices}

A nice way to think about the LCC question is to translate it to a question on the rank of matrices with a certain zero/non zero pattern. The following definition defines the particular pattern that arises in this setting.

\begin{define}[LCC-matrix]	
Let $A$ be an $nk \times n$ matrix over $\F$ and let $A_1,\ldots,A_n$ be $k \times n$ matrices so that $A$ is the concatenation of the blocks $A_1,\ldots,A_n$ placed on top of each other (so $A_\ell$ contains the rows of $A$ numbered $k(\ell-1) + 1, \ldots,k \ell$). We say that $A$ is a $(k,r)$-LCC matrix if, for each $i \in [n]$ the block $A_i$ satisfies the following conditions:
\begin{itemize}
	\item Each row of $A_i$ has support size at most $r+1$.
	\item All rows in $A_i$ are non zero in position $i$.
	\item The supports of two distinct rows in $A_i$ intersect only in position $i$.
\end{itemize}
\end{define}

The connection between LCC's and LCC-matrices will be clear from the following lemma:
\begin{lem}
Let $V = (v_1,\ldots,v_n) \in (\F^d)^n$ be a $(r,\delta)$-LCC with $\dim(V)=d$. Then, for $k = (\delta/r)n$, there exists a $(k,r)$-LCC matrix $A$ with  $n$ columns and with $\rank(M) \leq n - d$. Conversely, suppose there exists  a $(k,r)$-LCC matrix $A$ with $n$ columns and with $\rank(M) \leq n - d$. Then there exists an $(r,\delta)$-LCC $V = (v_1,\ldots,v_n)$ of dimension $\dim(V) \geq d$ with $\delta = k/n$.
\end{lem}
\begin{proof}
For the first direction let $B$ be the $n \times d$ matrix whose $i$'th row is the vector $v_i$. We will construct a $(k,r)$-LCC matrix $A$ such that $A \cdot B = 0$. This will prove that the rank of $A$ is at most $n-d$ since the rank of $B$ is $\dim(V)=d$. Let $M_1,\ldots,M_n$ be the $r$-Matchings given by Lemma~\ref{lem-qmatching}. Each $M_i$ can be used to define a block $A_i$ by adding a row for each set $R_{i,j} \in M_i$. We would like this row to have support $\{i\} \cup R_{i,j}$ and to have this row in the (left) kernel of $B$. This is possible since we know that $v_i$ is in the span of $\{ v_\ell \,|\, \ell \in R_{i,j} \}$. Thus, each block $A_i$ will have the required properties and we are done.

For the other direction, if $\rank(A) \leq n-d$ then there is a rank $d$ matrix $B$ of dimensions $n \times d$ so that $A \cdot B=0$. Let $V = (v_1,\ldots,v_n)$ be so that $v_i$ is the $i$th row of $B$. The structure of $A_i$ means that $v_i$ belongs to the span of $k$ disjoint sets $\{ v_\ell \,|\, \ell \in R_{i,j} \}$ with $R_{i,j}$ being the support of the $j$'th row of $A_i$, removing $\{i\}$. If we take any set $S \subset [n]$ of size at most $k = (k/n)n$ there will be some set $R_{i,j}$ that has empty intersection with $S$ and so $V$ will satisfy the LCC property.
\end{proof}
 
Thus, proving upper bounds on the dimension of LCCs is equivalent to proving lower bounds on the rank of LCC-matrices. Over fields of characteristic zero one can try using the results we saw on the rank of design matrices. This will work if the LCC matrix obtained from the code happens to be a design matrix.  The only families of codes we know, those based on Reed Muller codes, clearly satisfy this requirement (since every two points define a single line). Thus, we can use the bound on the rank of design matrices to show that there are no constant query LCCs over fields of characteristic zero that `look like' Reed muller codes or, more generally, whose decoding $r$-tuples satisfy a design condition (i.e., that every pair of  coordinates belongs to a small number of $r$-tuples used by the decoder). Clearly, one can construct artificial examples of LCCs whose decoding structure is {\em not} a design (simply repeat each coordinate twice). However, it is not out of the question to try and show that every LCC can be `modified' in some way to give a design-based LCC with comparable query complexity and dimension.

We conclude this section by mentioning a weaker type of local codes called {\em Locally Decodable Codes} (LDCs). These codes only require that the local decoding will be done for some basis $v_1,\ldots,v_d$ of the span of $V$. This type of decoding does not correct every symbol of the codeword but rather only symbols of the {\em message}. When $r >2$ there are constructions of LDCs that significantly outperform polynomial codes. See for example the excellent survey  \cite{Y_now}.

\ifediting
\bibliographystyle{alpha}
 \bibliography{incidence}

\end{document}
\fi

\cleardoublepage
\phantomsection
\addcontentsline{toc}{chapter}{References}
\bibliographystyle{alpha}
 \bibliography{incidence}

\end{document}